\documentclass{article}

\usepackage{blindtext} 

\usepackage[sc]{mathpazo} 
\usepackage[T1]{fontenc} 
\usepackage{microtype} 

\usepackage[english]{babel} 

\usepackage[hmarginratio=1:1,top=32mm,columnsep=20pt]{geometry} 
\usepackage[hang, small,labelfont=bf,up,textfont=it,up]{caption} 
\usepackage{booktabs} 

\usepackage{lettrine} 

\usepackage{enumitem} 
\setlist[itemize]{noitemsep} 

\usepackage{abstract} 

\usepackage{titlesec} 
\renewcommand\thesection{\Roman{section}} 
\renewcommand\thesubsection{\roman{subsection}} 
\titleformat{\section}[block]{\large\scshape\centering}{\thesection.}{1em}{} 
\titleformat{\subsection}[block]{\large}{\thesubsection.}{1em}{} 

\usepackage{fancyhdr} 
\usepackage{titling} 

\usepackage{hyperref} 

\usepackage{amssymb,amsbsy,amsmath,amsfonts}
\usepackage{graphicx}
\usepackage{subcaption}
\usepackage{biblio}
\usepackage{color}

\DeclareMathOperator*{\argmin}{\arg\!\min}
\usepackage{algorithm}
\usepackage{algpseudocode}

\pretitle{\begin{center}\Huge\bfseries} 
\posttitle{\end{center}} 
\title{Goal-oriented error control of stochastic system approximations using metric-based anisotropic adaptations} 
\author{%
\textsc{J. ~Van Langenhove} \\[1ex] 
\normalsize Sorbonne Universit\'es, UPMC Univ Paris 06,\\ \normalsize UMR 7190, CNRS, Institut Jean le Rond d'Alembert, \\ \normalsize F-75005 Paris, France \\ 
\and 
\textsc{D. Lucor} \\[1ex] 
\normalsize LIMSI, CNRS, Universit\'e Paris-Saclay, \\  \normalsize Rue John von Neumann,\\  \normalsize F-91405 Orsay cedex, France \\ 
\and 
\textsc{F. Alauzet} \\[1ex] 
\normalsize INRIA Saclay, Gamma3 Project,\\ \normalsize F-91126 Palaiseau, France \\ 
\and 
\textsc{A. Belme}\thanks{Corresponding author: \href{mailto:belme@dalembert.upmc.fr}{belme@dalembert.upmc.fr} } \\[1ex] 
\normalsize Sorbonne Universit\'es, UPMC Univ Paris 06,\\ \normalsize UMR 7190, CNRS, Institut Jean le Rond d'Alembert, \\ \normalsize F-75005 Paris, France \\ 
}
\date{} 
\renewcommand{\maketitlehookd}{%
\begin{abstract}
\noindent  The simulation of complex nonlinear engineering systems such as compressible fluid flows may be targeted to make more efficient and accurate the approximation of a specific (scalar) quantity of interest of the system. Putting aside modeling error and parametric uncertainty, this may be achieved by combining goal-oriented error estimates and adaptive anisotropic spatial mesh refinements. To this end, an elegant and efficient framework is the one of (Riemannian) metric-based adaptation where a goal-based {\em a priori} error estimation is used as indicator for adaptivity. \\ 
This work proposes a novel extension of this approach to the case of aforementioned system approximations bearing a stochastic component. In this case, an optimisation problem leading to the best control of the distinct sources of errors is formulated in the continuous framework of the Riemannian metric space. Algorithmic developments are also presented in order to quantify and adaptively adjust the error components in the deterministic and stochastic approximation spaces. \\
The capability of the proposed method is tested on various problems including a supersonic scramjet inlet subject to geometrical and operational parametric uncertainties. It is demonstrated to accurately capture discontinuous features of stochastic compressible flows impacting pressure-related quantities of interest, while balancing computational budget and refinements in both spaces.
\end{abstract}
}

\newcommand\BigR{\mathbb{R}}

\newcommand\BigE{\mathbb{E}}

\newcommand\Ccal{\mathcal{C}}

\newcommand\Ocal{\mathcal{O}}

\newcommand\bx{\mathbf{x}}

 \newcommand\bxi{\boldsymbol{\xi}}

\begin{document}

\maketitle

\section{Introduction}

In order to predict the performance of engineering systems involving fluid flow, the use of CFD simulations is widespread. Still, challenges remain to be addressed in order to 
improve the reliability of predictions of nonlinear flows in complex geometries (e.g. compressible airflow around an aircraft). Putting aside modeling error, sometimes leading to unpredictable outlying fluctuations \cite{Van_Langenhove_IJUQ2016}, a key element of the computational pipeline is the mesh; it discretizes the geometry, acts as the support of the numerical method and  needs to be tailored to the underlying non-linear equations governing the system physics. A poor mesh will fail to capture both geometrical and physical complexities and will thus yield inaccurate results \cite{loseille2008these}.
Critical flow-dominating features (boundary layers, shock waves, contact discontinuities, fluid density changes etc.) are very often concentrated in small regions and are inherently anisotropic. Efficient validation and certification of the numerical results require {\em optimal anisotropic} mesh refinements, i.e. 
a spatial discretization, for a given mesh size, leads to the best possible solution accuracy.\\
One is often solely interested in a specific Quantity of Interest (QoI) that is a functional of the solution field (e.g. the drag or lift of an airfoil). In this case goal-oriented error estimates for functional outputs may be used to drive the mesh adaptation \cite{loseille2010fully,Park-2002,Venditti-2003,DWIGHT20082845,ResminiIJNMF2017}.
These error estimates require the solution of an adjoint problem which takes into account the QoI. While there is a cost associated to solving this additional adjoint problem, goal-based mesh adaptation allows for a faster convergence in the QoI. For a more complete overview of mesh adaptation the reader is referred to \cite{Fidkowski2011,alauzet2016decade} and references therein.\\

Another important aspect of the computational framework is the impact of the  model parameters that are often {\em approximated} from calibration of experiments or other data sources, e.g. \cite{POETTE2012319}. A major task at hand, Uncertainty Quantification (UQ), consists in propagating this model parameter uncertainty throughout subsequent calculations to quantify the statistical behavior of the output QoI \cite{Sullivan_book}. Sometimes, the parametric exploration of a given model - out of its reliability range - cause a deterioration of the prediction, leading to unpredictable outlying fluctuations that need robust handling \cite{Van_Langenhove_IJUQ2016}. In other situations, the model response to parametric change is less biased, but nevertheless remains strongly nonlinear and requires adaptation. 
For instance the propagation of changing operating conditions and variations in the geometry of compressible flows  is well documented \cite{bijl2013uncertainty,ABGRALL2017507}, but particularly challenging, e.g. \cite{PoetteJCP2012}.
 This is because such systems are described by nonlinear conservation laws containing uncertain parameters and resulting in stochastic solutions with discontinuities in {\em both} physical and parametric dimensions.
As each sample in the parameter space corresponds to a costly CFD simulation, it becomes mandatory to lower as much as possible the number of samples needed in order to reach a prescribed QoI statistical accuracy. \\
The mainstay for UQ is the Monte Carlo (MC) method. However, despite its robustness, ease of implementation and parallelization, MC methods quickly become infeasible due to their high computational cost.
Other approaches rely on generalized Polynomial Chaos (gPC) \cite{xiu2002wiener} approximations which consists in constructing a parametrized polynomial approximation of the QoI response. 
For instance, the implementation using a Stochastic Galerkin (SG) method is a so called \emph{intrusive} method, which requires modifications to legacy codes and results in a large system of coupled equations to be solved. The nonintrusive variants, allow for the use of legacy solvers as black boxes similar to the approach used by MC methods. They take the form of a pseudospectral projection (PP) or a linear regression to find the coefficients for this polynomial basis, which is constructed to be orthogonal w.r.t. the underlying probability density. Stochastic Collocation (SC) methods on the other hand, aim at constructing interpolation functions for given coefficients/samples \cite{babuska2007stochastic,xiu2005high}. 
All of these methods are highly efficient when the QoI response is sufficiently smooth.
However, in the case where the stochastic response is of low regularity or even exhibit discontinuities, these global polynomial approximations are not amenable and may  suffer from Gibbs oscillations.
Therefore, {\em adaptive} solution techniques represent alternative approaches that may guarantee and accelerate convergence of the approximations, thanks to local relevant refinements of the response surface. These approaches include adaptive sparse grid collocation approaches \cite{nobile2008anisotropic,Resmini_IJNME2016} or piecewise polynomial approximations such as the Multi-Element gPC (ME-gPC) method \cite{wan2005adaptive, wan2006multi}, and Multi-Element Probabilistic Collocation Method (ME-PCM) \cite{foo2010multi,foo2008multi}. These latter methods split the parameter space into hyperrectangular subdomains on which a local polynomial approximation is constructed. While some anisotropy can be achieved, unless the discontinuities run along one of the principal axes, this partitioning  will not be able to perfectly capture these discontinuities. The degrading effects will still be present in the elements traversed by the discontinuities.
Furthermore, these methods often place the samples according to a tensor structure resulting in a fast increase in the number of samples as the number of elements increases.
Most of these problems are alleviated in the Minimal Multi-Element Stochastic Collocation (MME-SC) method \cite{jakeman2013minimal}, which relies on elements of irregular shapes where a discontinuity detector is used to split the parametric space into a minimal number of elements defined by the discontinuities.
The use of multi-wavelet representations have also been proposed \cite{le2004multi, le2004uncertainty} along with an $h-$refinement technique using hypercubes.\\
A different decomposition is proposed with the Simplex-Stochastic Collocation (SSC) method  \cite{witteveen2009adaptive}, a stochastic finite element-type method which uses a Newton-Cotes quadrature in (non hypercube) simplex elements. Gibbs oscillations are avoided using Local-Extremum Conserving limiter while later developments include higher order interpolation by use of Essentially Non-Oscillatory (ENO) stencils \cite{witteveen2013simplex}, subcell resolution \cite{witteveen2013subcell}, applications to non-hypercube domains \cite{witteveen2012simplex} and higher dimensional problems \cite{edeling2016simplex}.\\

Despite significant improvements obtained with the methods mentioned above, the indicators used to drive adaptive refinement are predominantely numerical {\em heuristics} computed at a local level in the parametric space and often tailored to specific problems. Moreover, these numerical techniques are often applied to problems with complex functionals but simple unimodal and smooth underlying probability measures that do not account for multivariate correlated inputs.
While for some it has been shown that adaptation following the given criterion will lead to convergence, these criteria do not permit a straightforward comparison between the error contribution introduced in the continuous approximation of the stochastic response and the error contribution introduced in each sample by the deterministic solver approximation. As a consequence, one is unable to determine in which approximation space the error is dominant: the deterministic (physical) space or the stochastic (parametric) space?
Potentially large computational gains in efficiency could in fact be realized if one was able to predict what level of mesh refinement is needed in each deterministic computation and what level of sampling refinement is needed in the representation of the stochastic response in order to remain within a given total error level.
By an extension of a goal-based \emph{a posteriori} error estimate, the authors in \cite{bryant2015error} achieved this splitting of the error into deterministic and stochastic contributions in order to drive adaptivity in both spaces. In the parametric space, along with global (anisotropic) $p-$refinement, ME style $h-$ and $h/p-$adaptivity methods are proposed which have the same deficiencies as the standard ME-gPC method: the refinement parameters driving $h-$ and $p-$adaptivity are chosen by the user and need to be tailored to the problem for optimal performance; secondly the use of hypercube elements does not permit efficient representation of discontinuities.
The extension of deterministic {\em a posteriori} goal-oriented error estimates has previously been applied to control numerical error in problems with uncertain input parameters \cite{almeida2010solution, butler2012posteriori,butler2011posteriori,mathelin2007dual}.
Nevertheless, it is known that extracting {\em anisotropic} information from {\em a posteriori} estimates is a difficult task \cite{alauzet2016decade}. Moreover, only very few works demonstrate error splitting and automated adaptivity in both spaces, e.g. \cite{bryant2015error}.
 Other attempts propose estimates of the total error (the deterministic and stochastic contributions combined) \cite{mathelin2007dual}, computable estimates of the error in each deterministic sample \cite{almeida2010solution}, or estimates of the error of the stochastic response approximation \cite{butler2012posteriori,butler2011posteriori} (although in \cite{butler2011posteriori} adaptivity in the parametric space is not addressed).\\

In this work the authors propose to extend {\em metric}-based mesh adaptation \cite{Gruau-2005,loseille2011continuous, loseille2011continuous2, loseille2010fully} to the partitioning of the parametric space.
The idea behind metric-based mesh adaptation is to make use of a powerful mathematical framework for efficient {\em anisotropic} multivariate parametric refinements. 
It relies on the use of a Riemannian metric space as a continuous model for the mesh. From this continuous representation, a unit mesh is computed for which the lengths of the element edges are (approximately) unity  in the {Riemannian metric}. Upon transformation to an Euclidean space, a stretched anisotropic mesh is  obtained.
The use of the metric is then coupled together with an {\em a priori} error estimation through the resolution of an optimisation problem, i.e.: find the {Riemannian metric} minimizing the continuous \emph{a priori} error estimate subject to a constraint on the complexity of the mesh.
The stochastic response is then approximated on the optimized mesh thanks to a tessellation of linear simplex/tetrahedron elements, well adapted for discontinuous responses \cite{witteveen2009adaptive,witteveen2012simplex}. In this respect it resembles the SSC method which uses a similar discretization. In contrast, the method proposed in this work makes use of an error estimate which will not only drive $h-$adaptivity in the parametric space, but also allows for a comparison with the error committed in each deterministic CFD simulations. As a result, one will be able to decide which solution space to refine in order to achieve a prescribed overall error level at minimal computational cost.\\

In summary, the aim of this paper is first: the extension of existing metric-based mesh adaptation to the stochastic space, enabling not only a highly effective adaptive approximation of discontinuous stochastic response surfaces, irrespective of the particular underlying probability density. And secondly we demonstrate computable estimates of both the deterministic and stochastic error contributions which allows for adaptive refinement using a goal-based \emph{a priori} error estimate in the space where the error dominates; in this way a global error level is reached at a reduced overall computational cost. \\

This paper is organized as follows: a short formal presentation of the error estimators are given in Section 2. The Riemannian metric-based approach and its interplay with the aforementioned material is detailed in Section 3. Some adaptive strategies and their algorithms are presented in Section 4. The proposed approach is then applied to several numerical examples and the description of the problems and their results can be found in Section 5. The paper ends with some conclusions and perspectives for future work.

\section{Formal error estimation}
An abstract model formulation is first considered in this section. It consists of a boundary-value problem defined on an open bounded domain $\Omega_{\bx} \subset \mathbb{R}^{d_{\bx}}$ with $d_{\bx}$ the dimension of the physical space. Furthermore, we suppose our problem depends on several uncertain parameters that will be defined as random variables. We therefore introduce a probability space  $(\Omega_{\bxi}, \mathcal{B}, \mathcal{P})$ where $\Omega_{\bxi}$ is the sample space, 
$\mathcal{B}$  is a $\sigma$-algebra and $\mathcal{P}$ the probability measure. We denote $\bxi(\omega)=(\xi_{1}(\omega), \xi_{2}(\omega), \ldots, \xi_{d_{\bxi}}(\omega))$ the vector of random variables that is sufficient to quantify our set of uncertain parameters, where $\omega \in \Omega_{\bxi}$ and $d_{\bxi}$ represents the dimension of the random input space. 
Throughout this paper, it will be assumed that the random input can be represented by a finite-dimensional probability space (the so called \emph{finite-dimensional noise assumption}). This ensures that the input can always be represented by a finite-dimensional set of random variables $\bxi(\omega)$.
Furthermore let $\rho_{\bxi}$ denote the joint probability density function of $\bxi$ and let $\Xi$ be the parameter space to which the $\bxi$'s belong. Each realization in the probability space corresponds, by a mapping defined by the probability density function of the random variables, to a parameter value. In short: $\bxi \in \Xi \equiv \prod_{i=1}^{d_{\bxi}} \Xi_{i}$ where $\Xi_{i}$ is the image of $\xi_{i}(\Omega_{{\bxi}})$.\\
Thus, {for a particular set of parameters $\bxi_{(i)}$},  the model problem can be cast in the following abstract form:
\begin{eqnarray}
    \Psi(\bxi_{(i)},w(\bxi_{(i)},\mathbf{x})) = 0
\label{PDEcont}
\end{eqnarray}
where $\Psi$ is a state equation (for example a steady Navier-Stokes or Euler system) having $w(\bxi,\mathbf{x})$ for exact solution. As we will see in the following, the ``exact" wording does not refer to a finite and fixed solution -- as it is a random quantity -- but designates the reference solution that is free of numerical errors. 
Indeed, the exact solution of a model is often out of reach for real-life engineering applications, and one has to rely on approximate numerical methods in order to approach it. For instance, {for the previously introduced set of parameters $\bxi_{(i)}$}, we may have a numerical tool that is capable of solving the corresponding deterministic discrete problem on a given spatial discretization $\mathcal{H}_{h_{\bx}}$, and produce an approximate solution of (\ref{PDEcont}) at the expense of a certain computational cost, i.e.:
\begin{eqnarray}
    \Psi_{h_{\bx}}(\bxi_{(i)},w_{h_{\bx}}(\bxi_{(i)},\mathbf{x})) = 0,
\label{PDEdisc}
\end{eqnarray}
with $w_{h_{\bx}}$ the deterministic discrete solution associated to sample $\bxi_{(i)}$ and solved on a mesh $\mathcal{H}_{h_{\bx}}$. As we will see later, there obviously exists an implicit dependence or ``coupling" between the choice of an adequate spatial discretization and the parameter value. For instance, certain parametric values might drastically affect the flow regime which will require an adequate mesh adaptation in order to produce a valid and reliable numerical approximation of the solution. 
~\\
In this work, we are more interested in an accurate approximation of a {\em scalar} quantity of interest (QoI) $j$ that 
is computed from the solution $w(\bxi,\mathbf{x})$ (and therefore depends on the uncertain random vector $\bxi$),
 than in the solution of the problem itself. We define our QoI sample obtained from the solution of the deterministic model built for the set of parameters {$\bxi_{(i)}$} and a given spatial discretization {$h_{\bx}$} as:
\begin{eqnarray}
    j(\bxi_{(i)}) = J\big (\bxi_{(i)}, w_{h_{\bx}}(\bxi_{(i)},\mathbf{x}) \big ),
\end{eqnarray}
where $J$ is the observation operator.
Point-wise evaluations of the QoI are limited in practice as each evaluation involves a costly simulation. Predictions of $j$ for new parametric values or 
statistical information (e.g. moments) of interest related to $j$ may be more efficiently computed from continuous approximations, also called response surface models or surrogate models, which are built across the span of the $\bxi$ uncertain parametric ranges. One way to construct these approximations is via a discretization or partition that we call $\mathcal{H}_{h_{\bxi}}$ of the multivariate stochastic (uncertain parameter) space. \\
 Indeed, as it will be explained later in this paper, the random input space is discretized such that a finite number of samples (or design of experiments (DoE)) are chosen to be mapped throughout the deterministic model and a surrogate stochastic model for $j$ is built from the computations on these samples:
\begin{eqnarray}
    j_{h_{\bxi}}(\bxi) = J_{h_{\bxi}} \big (\bxi, w_{h_{\bx}}(\bxi,\mathbf{x}) \big ).
\end{eqnarray}
In this work, we aim to control errors committed on $j$.
We can thus define the total error committed on our QoI for a certain set of parameters:
\begin{eqnarray}
    \delta j(\bxi) \equiv j-j_{h_{\bxi}} = J(\bxi,w) - J_{h_{\bxi}} (\bxi, w_{h_{\bx}} ).
\label{deltaJtotal}
\end{eqnarray}
where $J_{h_{\bxi}}(\bxi,w_{h_{\bx}})$ denotes the approximate QoI.\\

\noindent The \textit{exactness} of the QoI depends thus on two components:
\begin{enumerate}
\item on the {\em deterministic} discrete solution $w_{h_{\bx}}({\bxi})$ error and its choice of spatial discretization $\mathcal{H}_{h_{\bx}}$,
\item on the {\em stochastic} error committed by discretizing the stochastic space $\mathcal{H}_{h_{\bxi}}$ and building the surrogate model $J_{h_{\bxi}}$, {i.e. the representation  of $J$ on $\mathcal{H}_{h_{h_{\bxi}}}$ which is equivalent to an \textit{interpolation error} $J-\Pi_{h_{\bxi}} J$}.  
\end{enumerate} 
Due to the segregation in the approximation space of the problem solution and its inferred QoI, the two sources of errors owing to deterministic and stochastic discretizations/approximations, are treated separately, and an error control strategy is applied either in the deterministic approximation space (e.g. mesh control/adaptation) {\cite{Venditti-2003,Fidkowski2011,palacios2012robust,loseille2008these}} or into the stochastic approximation space \cite{mathelin2007dual, witteveen2009adaptive, butler2011posteriori}. However, 
the interplay between the errors as well as which one dominates the computation of $j$ remain a very important question.\\
In order to give some elements of answer, we split the total error in two contributions:
{
\begin{eqnarray}
    \delta j(\bxi) = \underbrace{J(\bxi,w) - J_{h_{\bxi}}(\bxi,w)}_{\eta(h_{\bxi})} + \underbrace{ J_{h_{\bxi}}(\bxi,w) - J_{h_{\bxi}} (\bxi,w_{h_{\bx}})}_{\varepsilon(h_{\bx})}.
\label{deltaJtotalsplit}
\end{eqnarray}
where $\eta(h_{\bxi})$ only accounts for errors due to the discretization of stochastic space, i.e. the representation of $J$ on $\mathcal{H}_{h_{h_{\bxi}}}$. Therefore, it corresponds to an interpolation error. And, where $\varepsilon(h_{\bx})$ only accounts for point-wise errors due to the evaluation of $J$ in the deterministic space, and therefore amounts to an implicit error.}
While both error contributions, $(\varepsilon,\eta)$ are random quantities depending on $\bxi$, we assume that they will be controlled (in a complementary fashion) in the different approximation spaces: \\
-- $\varepsilon$ will be controlled via {\em deterministic} refinement
and\\
 -- $\eta$ will be controlled via {\em stochastic} refinement.

In practice, $\delta j$ being a random quantity, we will be interested in lowering the {\em average} QoI total error, that we express as:
\begin{equation}
    \overline{\delta j }\equiv \BigE \left [ \delta j(\bxi) \right ] = \underbrace{\BigE \left [ {\eta} \right ]}_{\bar{{\eta}}} + \underbrace{\BigE \left [ {\varepsilon} \right ]}_{\bar{{\varepsilon}}}
\end{equation}
A simple adaptive approach might be chosen in order to lower both contributions of the right-hand side of the previous equation, i.e.: \\
1. optimizing the set of spatial discretizations $\mathcal{H}_{h_{\bx}}$ in order to diminish (or balance) the $\varepsilon$ error contribution of the samples of a given DoE and concurrently \\
2. update the DoE to lower $\BigE \left [ \eta \right ]$ for a given set of $\mathcal{H}_{h_{\bx}}$.

{In the next subsections we propose an estimation for both deterministic $\varepsilon$ and stochastic $\eta$ errors which will serve as indicators for the adaptive process.}

\subsection{Goal-oriented deterministic error estimate}
\label{detErr}
The deterministic error contribution {$\varepsilon(h_{\bx})$ in relation} (\ref{deltaJtotalsplit}), which appears each time the parameter set is given 
 is discussed in this section. Note that we have replaced subscript $\cdot\,_{h_{\bx}}$ by $\cdot\,_h$ for ease of notation. \\
 An {\em a priori} error estimation of the deterministic error $\varepsilon$ is proposed, involving the computation of an adjoint problem. Only the main results are recalled and we refer to \cite{loseille2008these,loseille2010fully} for a more detailed analysis. \\

First, the variational problem associated to state equations (\ref{PDEcont}-\ref{PDEdisc}) on an appropriate Hilbert space of solutions $\mathcal{V}$ and respectively subspace $\mathcal{V}_h$ are introduced hereafter:
\begin{eqnarray}
    \mbox{Find}~w(\bxi_{(i)},\cdot) \in \mathcal{V} ~\mbox{such that}~ \forall \varphi \in \mathcal{V} , 
\quad \left( \Psi(w)\,,\,\varphi \right) = 0,
\label{eqstate}
\end{eqnarray}
with the associated discrete variational formulation:
\begin{eqnarray}
    \mbox{Find}~w_h(\bxi_{(i)},\cdot) \in \mathcal{V}_h ~\mbox{such that}~ \forall \varphi_h \in \mathcal{V}_h , 
\quad \left( \Psi_h(w_h)\,,\,\varphi_h \right) = 0.
\label{eqstatediscret}
\end{eqnarray}
Furthermore, we assume some level of regularity for our QoI such that :
\begin{equation}
    j(\bxi_{(i)}) \in \mathbb{R}~;~ j=J(\bxi_{(i)},w)=(g,w)=\int_{\Omega_{\bx}} g\, w(\bx,\bxi_{(i)}) \mbox{d}\bx
\label{j-a}
\end{equation}
where $g$ represents the deterministic jacobian of $J$.
We introduce the {\it continuous adjoint} solution $w^*$ of the following system:
\begin{eqnarray}
    w^*(\bxi_{(i)},\cdot) \in \mathcal{V}~,~ \forall \psi  \in \mathcal{V}~,~ \left(\frac {\partial \Psi}{\partial w }(w )\psi,w^*\right)=(g,\psi)~.
\label{psi-abis}
\end{eqnarray}
We assume that both state solution $w$ and adjoint solution $w^*$ are smooth enough, such that $w,w^* \in \mathcal{V} \cap \mathcal{C}^0$.\\
 Using the fact that $\mathcal{V}_h\subset \mathcal{V}$,
the following error estimates for the unknown can be written:
\begin{equation}\label{psi-approx}
(\Psi_h(w ),\varphi_h)-(\Psi_h(w _h),\varphi_h) = 
(\Psi_h(w ),\varphi_h)-(\Psi(w),\varphi_h) =
((\Psi_h-\Psi)(w ),\varphi_h).
\end{equation}
The objective here is to estimate the following approximation error committed on the functional:
$$
(g,w)-(g,w_h) = (g,w-w_h)
$$
where $w$ and $w_h$ are respectively solutions of
(\ref{PDEcont}-\ref{PDEdisc}). \\
{The idea is then to define local error estimation} to be used as a guide for anisotropic mesh refinement. {Interpolation errors are known to provide useful local information and in our case, as it will be stated later, directions and sizes for anisotropic mesh refinement. We introduce thus} an interpolation operator:
 $
 \pi_h : \mathcal{V} \cap \mathcal{C}^0 \rightarrow \mathcal{V}_h
 $ 
 which allows for a simple decomposition of the {approximation} error :
 \begin{eqnarray}
     J(\bxi_{(i)},w) - J(\bxi_{(i)},w_h) = (g,w-\pi_h w) + (g,\pi_h w -w_h)
  \label{deltaJsplit}
 \end{eqnarray}
 We recognize the first error term as the \textit{interpolation error} and we can show (see for example \cite{loseille2010fully}) that the second error term, called here \textit{implicit error}, can also be expressed in terms of interpolation errors. Finally we get: 
 \begin{equation}
     \varepsilon = J(\bxi_{(i)},w) - J(\bxi_{(i)},w_h) \approx \left( w^*,\left( \Psi_h(\bxi_{(i)},w) - \Psi(\bxi_{(i)},w)\right)\right).
 \label{DetErrorEstFinal}
 \end{equation}

The {\em a priori} error estimate (\ref{DetErrorEstFinal}) is used here {to control the deterministic error, i.e the implicit error in the stochastic space}. More specifically, we  build our mesh adaptation as an optimisation problem where we seek the optimal mesh that minimises (\ref{DetErrorEstFinal}) under the constraint of a given number of mesh nodes; the adjoint state $w^*$ acting as a Lagrange multiplier associated to the equality constraint (\ref{PDEcont}). 
A continuous formulation to this optimisation problem is proposed in Section \ref{sec:metric}, using the concept of Riemannian metric space. 

\subsection{Stochastic error estimate}
The stochastic error contribution {$\eta(h_{\bxi})$ in relation }(\ref{deltaJtotalsplit}) is discussed next.
In the following, we have replaced $\cdot\,_{h_{\bxi}}$ by $\cdot\,_h$ for ease of notation.
We propose an error estimate of the numerical approximation of the solution in the parametric space inspired by the previous developments. Motivated by the need for anisotropic information, we wish to deploy the notion of Riemannian metric field in the stochastic space. For the type of applications we consider, it is common knowledge that the dependence of the QoI on the random variables is anisotropic as we often encounter singularities and sharp response gradients. To this purpose, we follow and adapt the formulation of \cite{loseille2008these}, in order to express the stochastic error through the $L^p$-norm of the interpolation error:
\begin{eqnarray}
\left \|  J(\bxi, w) - \pi_h J(\bxi, w) \right \| _{L^p(\Xi)} 
\label{StochEE}
\end{eqnarray}
where $\pi_h$ is the (linear) interpolation operator in the parameter space. For the purpose of this paper we will focus on the $L^1$-norm of the interpolation error in the stochastic space. This is a purely practical choice, well adapted for approximation of potentially discontinuous solutions, but there is no restriction in using a different $p$-norm.\\
For deterministic problems, this kind of approach named ``feature-based" mesh adaptation has been proposed to capture all the scales/singularities of the system, and has been applied to deterministic CFD problems \cite{alauzet2010high,LoseilleAIAA2007,alauzet2016decade}. In this case, a sensor is defined (for CFD applications a sensor will be a prescribed field: density, mach, ...) and some norm of the
interpolation error associated to the sensor is controlled by anisotropic mesh refinement.
Several differences appear naturally when transposing this approach in the stochastic context, leading to a different interpretation.
First, what plays the role of the sensor in our case is the stochastic scalar QoI $j$. 
Second, the $L^1$-norm of the interpolation error of $j$ on the parameter space, now equipped with a probability measure $\mathcal{P}$, introduces the parameters probability density function in the formulation:

\begin{eqnarray}
    \bar{\eta} = \mathbb{E}\left [ \eta \right ] = \int_{\Xi} | J(\bxi, w) - \pi_h J(\bxi, w)  |  \rho_{\bxi} \mbox{d}\bxi,
\label{StochErrorDiscrete}
\end{eqnarray}
where $\rho_{\bxi}$ is again the joint probability density function (pdf) of $\bxi$. The probability density function acts as a weighting of the interpolation error, but the formulation is more straightforward than relation (\ref{DetErrorEstFinal}) as it does not involve an adjoint solution.  \\
The convergence of the approximation in the $L^1$-norm, i.e. the fact that $\BigE_{h_{\bxi} \rightarrow 0 }\left [ \eta \right ] \rightarrow 0$ will ensure a {\em convergence in the mean} of the (piecewise linear) interpolated surrogate. Thanks to Markov's inequality, this convergence in the mean will ensure that it converges in probability, which in turn implies convergence in distribution \cite{Grimmett_book2001}. Moreover, triangle inequality\footnote{This is sometimes called the reverse triangle inequality:
\begin{equation}
 \big |  \| a  \| - \| b \| \big | \leq \left \|  a - b \right \|.
\end{equation}} will insure that the mean value of the surrogate $\BigE[\pi_h J]$ will converge to the exact QoI mean $\bar{j} \equiv \mathbb{E}[J] $ at least as fast as the expectation of the interpolation error, i.e.:
\begin{equation}
    \big |  \bar{j} - \BigE[\pi_h J] \big | \leq \mathbb{E}\left [ \eta \right ].
\end{equation}
~\\
Other formulations for the error estimation of $j$ in the stochastic space are conceivable, but at this stage it seems burdensome to envision a fully stochastic adjoint-based approach due to lack of any differential operators along the parametric coordinates. \\

The error estimate (\ref{StochErrorDiscrete}) will be used later in this paper as refinement indicator to drive adaptivity in the parameters space. More precisely, following the deterministic approach, we will formulate the stochastic adaptation problem as an optimisation problem where we seek the optimal ``mesh" that minimises  (\ref{StochErrorDiscrete}) under the constraint of a given number of samples. 
A continuous formulation to this optimisation problem is proposed in Section \ref{sec:metric}, using the concept of Riemannian metric space. The optimal stochastic metric, solution to this problem, will be expressed in terms of Hessians of $j$ in the parameters space. 

\subsection{Application to nonlinear conservation laws: case of steady compressible Euler flows}
\label{subsec:Euler}
The previously defined formal error analysis {of Section \ref{detErr}} is {developed} next {for} the compressible Euler system which will be the CFD model solved in the numerical exemple section. Note that we provide the details of the discretization for a fixed $\bxi$, thus for ease of notation, the dependence on both $\bxi$  and $\bx$ will not be explicitly written out. \\

The two-dimensional steady Euler equations are set in the computational 
domain $\Omega \subset \mathbb{R}^2$ of boundary denoted $\Gamma$.
An essential ingredient  of our discretization and of our analysis 
is the elementwise linear interpolation operator. In order to use it easily,
we define our working functional space as
$\mathcal{V}=\left[H^1(\Omega)\cap \mathcal{C}(\bar \Omega)\right]^4$, 
that is the set of measurable functions that are continuous with 
square integrable gradient.
We formulate the Euler model in a compact variational formulation in the
functional space $\mathcal{V}$ as follows:
\begin{eqnarray}
&& \hspace*{-20mm} \mbox{Find}~w \in \mathcal{V} ~\mbox{such that}~ \forall \varphi \in \mathcal{V} , 
\quad \left( \Psi(w)\,,\,\varphi \right) = 0  \nonumber \\
\mbox{with}~~\left( \Psi(w)\,,\,\varphi \right)  &=&
\int_\Omega \varphi  \, \nabla \cdot {\mathcal F}(w) \, \mathrm{d}\Omega \,    
~-~  \int_{\Gamma} \varphi \, {\hat {\mathcal F}(w)}.{\bf n} \, \mathrm{d}\Gamma \,    \,. 
\label{eulercontinuebis}
\end{eqnarray}
In the above definition, $w = {}^t(\rho, \,\rho {\mathbf u}, \,\rho E)$ is the vector of conservative flow variables and 
$\mathcal{F}(w)=(\mathcal{F}_1(w),\mathcal{F}_2(w))$ is the usual Euler flux:
$$\mathcal{F}(w) = {}^t \left( \rho \mathbf u, \,\rho u {\mathbf u} + p \mathbf e_x, \, \rho v {\mathbf u} + p \mathbf e_y,  
\,\rho {\bf u} H \right) \,.$$ 
We have noted $\rho$ the density, $\mathbf u = (u , v)$ the velocity vector, $H = E + p/\rho$ is the total enthalpy,  
$E=T+\frac{\|\mathbf u \|^2}{2}$ the total energy and 
$p=(\gamma-1) \rho T$ the pressure with $\gamma~=~1.4$ the ratio of specific heat capacities, 
$T$ the temperature, and $(\mathbf e_x, \mathbf e_y)$ the canonical basis. \\
Note that   
$\bf n$ is the outward normal to $\Gamma$,  and 
the boundary flux $\hat{\mathcal{F}}$ contains the different 
boundary conditions, which involve standard inflow, outflow and slip 
boundary conditions.

\paragraph{Discrete state system}

As a spatially semi-discrete model, we consider the Mixed-Element-Volume
formulation \cite{alauzet2010high,cournede2006positivity}. As in \cite{loseille2010fully}, 
we reformulate it under the form of a finite element 
variational formulation.
We assume that $\Omega$ is covered by a finite-element partition in 
simplicial elements denoted $K$. The mesh, denoted by $\mathcal{H}$ is the set of the elements.
Let us introduce the following approximation space:
\begin{eqnarray}
\displaystyle
\mathcal{V}_h ~=~ \left\{ \varphi_h \in \mathcal{V} 
~\big|~  {\varphi_h}_{| K} \mbox{ is affine }
\forall K \in \mathcal{H}   \right\}. 
\nonumber
\end{eqnarray}
Let $\Pi_h:~ \mathcal{V} \rightarrow ~\mathcal{V}_h$ be the usual $\mathcal{P}^1$ projector. 
The weak discrete formulation writes:
\begin{eqnarray}
&& \hspace*{-20mm}  \mbox{Find}~w_h \in \mathcal{V}_h ~\mbox{such that}~
\forall \varphi_h \in \mathcal{V}_h , ~
\; \left( \Psi_h(w_h) \,,\, \varphi_h \right) = 0, \nonumber \\
\mbox{with:} ~~\left( \Psi_h(w_h) \,,\, \varphi_h \right) &=&
 \int_{\Omega} \varphi_h  \nabla \cdot {\mathcal F_h}(w_h) \,\mathrm{d}\Omega \,    
~-~  \int_{\Gamma} \varphi_h {\hat {\mathcal{F}_h}(w_h)}.{\bf n} \,\mathrm{d}\Gamma \,     
~+~ \int_{\Omega} \varphi_h \, D_h(w_h) \, \mathrm{d}{\Omega} \,   , 
\label{eulerdiscrete3}
\end{eqnarray}

with ${\mathcal F_h} = \Pi_h{\mathcal F}$ and $\hat {\mathcal{F}_h}=\Pi_h\hat {\mathcal{F}}$.
~The $ D_h$  term accounts for the numerical diffusion. 
In short, it involves the difference between the Galerkin central-differences 
approximation and a second-order Godunov approximation \cite{cournede2006positivity}. 
In the present study,
we only need to know that for smooth fields, the $D_h$ term 
is a third order term with respect to the mesh size. For shocked fields,
monotonicity limiters become first-order terms.

\medskip
Practical experiments are done with the CFD software {\tt Wolf}. The numerical scheme is vertex-centered and uses a particular edge-based formulation. 
This formulation consists in associating with each vertex of the mesh a control volume (or Finite-Volume cell) 
built by the rule of medians. This flow solver uses a HLLC approximate Riemann 
solver to compute numerical fluxes. 
A high-order scheme is derived according to a MUSCL type method using downstream and upstream tetrahedra. 
Appropriate $\beta$-schemes are adopted for the variable extrapolation which gives us a 
very high-order space-accurate scheme for the linear advection on cartesian triangular meshes. 
This approach provides low diffusion second-order space-accurate scheme in the non-linear case. 
The MUSCL type method is combined with a generalization of the Superbee limiter with 
three entries to guarantee the TVD property of the scheme. 
More details can be found in~\cite{alauzet2010high}. 

\paragraph{Deterministic error estimate applied to Euler model}
We replace in Estimation (\ref{DetErrorEstFinal}) operators $\Psi$ and $\Psi_h$ by their expressions given by Relations (\ref{eulercontinuebis}) and (\ref{eulerdiscrete3}).
As in \cite{loseille2010fully}, where it was observed that even for shocked flows, 
it is possible to neglect the numerical viscosity term, we follow the same guideline. We finally get the following 
simplified error model for the deterministic error associated to the Euler model:
\begin{eqnarray}
\label{estimation-euler0}
\varepsilon \approx 
\int_{\Omega} w^*\,  \nabla . \big(  {\mathcal F}(w) 
- \Pi_h{\mathcal F}(w) \big) \,
\mathrm{d}\Omega -
\int_{\Gamma} w^* \, \big({\hat {\mathcal{F}}}(w) 
- \Pi_h{\hat {\mathcal{F}}(w))}\big)\cdot{\bf n} \;
\mathrm{d}\Gamma~ \,.
\end{eqnarray}
Integrating by parts leads to:
\begin{eqnarray}
\label{estimation-euler}
\varepsilon \approx 
  \int_{\Omega} \nabla w^*\, \big(  {\mathcal F}(w) 
- \Pi_h{\mathcal F}(w) \big) \,
\mathrm{d}\Omega -
\int_{\Gamma} w^* \, \big({\bar {\mathcal{F}}}(w) 
- \Pi_h{\bar {\mathcal{F}}(w))}\big)\cdot{\bf n} \,\mathrm{d}\Gamma\,.
\end{eqnarray}
with $\bar {\mathcal{F}}=\hat {\mathcal{F}}-\mathcal{F}$.
We observe that this error estimate of the deterministic error is expressed in 
terms of interpolation errors of the Euler fluxes weighted by continuous functions $w^*$ and $\nabla w^*$. \\

The integrands in error estimation~(\ref{estimation-euler}) contain positive 
and negative contributions which can compensate for some particular 
meshes. We prefer not to rely on  these parasitic 
effects and to slightly over-estimate the error.
To this end, all integrands are bounded by their absolute 
values:
\begin{eqnarray}
\varepsilon  &\le&
  \int_{\Omega} |\nabla w^*|  \, 
|{\mathcal F}(w) -\Pi_h{\mathcal F}(w)| \,  \mathrm{d}\Omega +   \int_{\Gamma} |w^*| \,  |({\bar {\mathcal F}}(w)
-\Pi_h {\bar {\mathcal F}}(w)).{\bf n}| \, 
\mathrm{d}\Gamma \,    \,.
\label{majobis}
\end{eqnarray}
\section{Total error estimate}
The total error committed on our stochastic QoI $j(\bxi)$ for the Euler model is the summation of the \textit{stochastic error estimate}, the integrand of (\ref{StochErrorDiscrete}) and the deterministic error estimate (\ref{majobis}). These error terms are approximated as weighted interpolation errors. As previously stated, the purpose of this paper is to control both errors through mesh{\footnote{Here the "mesh" nomenclature is used as a generic term and also applies to the sampling and related discretization of the parametric space.}}  adaptivity, in both spaces. The main idea is to set the mesh adaptation problem as an optimisation problem where we seek the optimal mesh that minimises these interpolation errors under the constraint of a fixed number of nodes (or samples for the stochastic problem).  Due to the segregated  construction of the approximation space, the two sources of error will be treated independently. We will use the continuous framework of Riemannian metric space to formulate and solve these optimisation problems. The next section summarizes the main ingredients of the Riemannian metric space and the resulting optimal discretizations are derived for Euler flows.
\section{Riemannian metric-based approach for anisotropic mesh adaptation}
\label{sec:metric}
As previously mentioned, we wish to be able to accurately capture the anisotropic behaviour of our shock-dominated flows in both physical and parameters spaces. 
An efficient approach to generate anisotropic meshes in the physical space has been introduced in \cite{loseille2011continuous,loseille2011continuous2}. In order to generate anisotropic meshes, one must be able to prescribe at each point of the domain privileged sizes and orientations for the elements. The use of Riemannian metric spaces is an elegant and efficient way to achieve this goal. The main idea of metric-based mesh adaptation, initially introduced in \cite{george1991creation}, is to generate a \textit{unit mesh} in a prescribed Riemannian metric space. Consequently, the generated mesh will be uniform and isotropic in the Riemannian metric space while it will be anisotropic and adapted in the usual, Euclidean space. This approach will be used to generate anisotropic meshes in both physical and stochastic spaces. We will briefly describe the main ingredients of this approach and we refer to the works cited in this section for further details. \\

These differential geometry notions are more than just a simple tool for mesh generation; the Riemannian metric spaces can be seen as continuous models representing meshes. The fundamental consequence is that all kind of mathematical analysis can be performed using such spaces for which powerful mathematical tools are available.
More precisely, it allows us to define proper differentiable 
optimization~\cite{absil2009optimization,arsigny2006log}, {\em i.e.,} to use
a calculus of variations on continuous meshes which cannot be applied to the class of discrete meshes. 
A continuous mesh $\mathbf{M}$ of computational domain $\Omega$ 
is identified with a Riemannian metric field~\cite{berger2003panoramic}
$\mathbf{M} = (\mathcal{M}({\bf x}))_{{\bf x}\in \Omega}$.
For all ${\bf x}$ 
of a computational domain $\Omega \subset \mathbb{R}^3$, a Riemannian metric $\mathcal{M}({\bf x})$ 
is a symmetric $3\times 3$ matrix having $(\lambda_i({\bf x}))_{i=1,3}$ as eigenvalues 
along the principal directions 
$\mathcal{R}({\bf x}) = ({\bf v}_i({\bf x}))_{i=1,3}$.  
Sizes along these directions are denoted 
$(h_i({\bf x}))_{i=1,3} = (\lambda_i^{- \frac 1 2}({\bf x}))_{i=1,3}$ 
and the three {\it anisotropy quotients} $r_i$ are defined by: 
$ r_i = h_i^3 \left( \displaystyle h_1h_2h_3  \right)^{-1}$.\\

The diagonalisation of $\mathcal{M}({\bf x})$ writes:
\begin{equation}
\mathcal{M}({\bf x}) = d^\frac{2}{3}({\bf x}) \, \mathcal{R}({\bf x})
\left[ \begin{array}{ccc}
r_1^{-\frac{2}{3}}({\bf x})  \\
&r_2^{-\frac{2}{3}}({\bf x})  \\
&  & r_3^{-\frac{2}{3}}({\bf x})   
\end{array}
\right] \mathcal{R}^{t}({\bf x}),
\label{def-M}
\end{equation}
The {\it node density} $d$ is equal to:
$\displaystyle d= \left(  h_1h_2h_3 \right)^{-1}  = \left( \lambda_1\lambda_2\lambda_3\right)^{\frac{1}{2}} = \sqrt{\det({\mathcal M})}$. By integrating the node density, we define the {\it complexity} $\mathcal{C}$ of a continuous mesh 
which is the continuous counterpart of the total number of vertices: 
$$
\mathcal{C}({\mathbf{M}}) = \int_{\Omega} d({\bf x}) \, \mathrm{d}{\bf x} = 
\int_{\Omega} \sqrt{\det(\mathcal{M}({\bf x}))} \, \mathrm{d}{\bf x} .
$$
A discrete mesh $\mathcal{H}$ is a {\bf unit mesh} with respect to Riemannian metric space $\mathbf{M}$, 
if each tetrahedron $K \in \mathcal{H}$,
defined by its list of edges $({\bf e}_i)_{i=1\dots6}$,
verifies:
$$
\displaystyle \forall i \in [1,6], \quad L_{\mathcal M}({\bf e}_i) \in \left[ \frac{1}{\sqrt{2}},\sqrt{2}\right] 
\quad \mbox{and} \quad 
Q_{\mathcal M}(K) \in [\alpha,1] ~\mbox{ with }~ \alpha > 0 \, ,
$$
in which the length of an edge $L_{\cal M}({\bf e}_i)$ and the quality of an element $\displaystyle Q_{\cal M}(K)$
are defined as follows:
$$
\begin{array}{l}
\displaystyle Q_{\mathcal M}(K) = \frac{36}{3^{\frac{1}{3}}} \frac{|K|_{\mathcal M}^\frac{2}{3}} 
{\sum_{i=1}^6 L^2_{\mathcal M}({\bf e}_i)}  \in [0,1], \, \mbox{ with }
\displaystyle |K|_\mathcal{M} = \int_{K} \sqrt{\det(\mathcal{M}({\bf x}))} \, \mathrm{d}{\bf x}, \\[2.5ex]
\displaystyle  \mbox{and } \,  L_{ \mathcal{M}}({\bf e}_i) =\int_{0}^{1} \sqrt{^t{\bf ab} \; 
{ \mathcal{M}}({\bf a} + t \, {\bf ab}) \; {\bf ab} } \;   , ~~ \mbox{with} \; {\bf e}_i={\bf ab}.
\end{array}
$$
As in \cite{loseille2008these} we choose a tolerance  $\alpha$ equal to  $0.8$.\\

\subsection{Continuous error model}

Let us emphasize first that the set of all the discrete meshes that 
are unit meshes with respect to a unique $\mathbf{M}$ contains 
an infinite number of meshes. \\
Given a smooth function $u$,
to each unit mesh $\mathcal{H}$ with respect to $\mathbf{M}$ corresponds a local interpolation error  
$|u - \Pi_{h} u|$.
In \cite{loseille2011continuous,loseille2011continuous2}, it is shown that all these interpolation errors
are well represented by the so-called continuous interpolation error related
to $\mathcal{M}$, which is expressed locally in terms of 
the Hessian $H_u$ of $u$ as follows: 
\begin{eqnarray}
\label{int-contbis}
(u - \pi_{\cal M}u)({\bf x}) &=&  \frac{1}{10}\mbox{trace}(\mathcal{M}^{-\frac{1}{2}}({\bf x}) \, |H_u({\bf x})| \, 
\mathcal{M}^{-\frac{1}{2}}({\bf x})) 
\nonumber \\
&=&  \frac{1}{10} d({\bf x})^{-\frac{2}{3}} 
\sum_{i=1}^3 r_i({\bf x})^{\frac{2}{3}} {}^t{\bf v_i}({\bf x}) \, |H_u({\bf x})| \, {\bf v_i}({\bf x}),
\end{eqnarray}
where we denoted $\pi_{\cal M}$ the continuous linear interpolate, and $u - \pi_{\cal M}u$ represents the continuous dual of the discrete interpolation error. \\
%
In this paper we aim to control weighted interpolation errors as we have shown in (\ref{DetErrorEstFinal}) and (\ref{StochErrorDiscrete}). We will define next the continuous dual of these error models for the concrete case of the steady Euler flows. We emphasize that the weight of the deterministic error is the adjoint state, while for the stochastic error the probability density function acts as a weight.
\paragraph{Deterministic continuous error model}
Working in the continuous framework  enables us to write Estimate (\ref{majobis}) in the following continuous form:
\begin{eqnarray}
\label{kernel}
 (g,w_h-w) \approx {\mathbf E}(\mathcal{M}) & = & \int_\Omega |\nabla  w^* |  \, 
 |{\mathcal F}(w) -\pi_{\mathcal{M}}{\mathcal F}(w)| \, d\Omega \nonumber \\
 & + & \int_{\Gamma} |w^*| \,  |({\bar {\mathcal F}}(w)-\pi_{\mathcal M} {\bar {\mathcal F}}(w)).{\bf n}| \, d{\Gamma}.
\end{eqnarray}
We observe that the second term introduces a dependency of the error with respect 
to the boundary surface mesh. In the present paper, 
we discart this term and refer to \cite{loseille2010fully} for a discussion of the 
influence of it.
Then, introducing the continuous interpolation error, we can write
the simplified error model as follows:
$$
\mathbf{E}_{\bx} ({\mathcal M}) ~=~ 
  \int_{\Omega}  
\mbox{trace} \left( \mathcal{M}^{-\frac{1}{2}}({\bf x}) \, 
{\mathbf H}{_{\bx}}({\bf x}) \, \mathcal{M}^{-\frac{1}{2}}({\bf x}) \right)  \,
\mathrm{d}\Omega \,   
$$
\begin{equation}\label{Hgo}
\mbox{ with }~~{{\mathbf H}}{_{\bx}}({\bf x}) = 
\sum_{j=1}^4 
\left(  
\left| \frac{\partial w_j^*}{\partial x}({\bf x}) \right| \cdot \big| H({\mathcal F}_{1}(w_j))({\bf x}) \big| +
\left|  \frac{\partial w_j^*}{\partial y}({\bf x}) \right| \cdot \big| H({\mathcal F}_{2}(w_j))({\bf x}) \big|
\right), 
\end{equation}
Here, $H({\mathcal F}_{i}(w_j))$ denotes the Hessian of the $j^{th}$ component of the vector $\mathcal{F}_i(w)$. \\
\noindent
The \textit{deterministic mesh optimization problem} is formulated as: 
\begin{eqnarray}\label{OptPbDet}
    \mbox{Find } {\mathcal M}^{opt}_{\bx} = 
    \argmin_{{\mathcal M}_{\bx}} \, \mathbf{E}_{\bx} ({\mathbf M}),
\end{eqnarray}
under the constraint of bounded mesh fineness:
\begin{eqnarray}\label{adjmet-C}
    {\mathcal C}({\mathbf{M}})  = \mathcal{C}_\bx,
\end{eqnarray} 
where $\mathcal{C}_\bx$ is a specified complexity (i.e. continuous counterpart of the number of nodes $N_{\bx}$). A calculus of variations (see \cite{loseille2008these}) gives the following solution to the deterministic optimisation problem:
\begin{align}
\label{OptMdet}
\mathcal{M}_{\bx}^{opt} = {\mathcal{C}_{\bx}}^{\frac{2}{d_{\mathbf{x}}}} \left(\int_{\Omega_{\mathbf{x}}} \det(\mathbf{H}_{\mathbf{x}})^{\frac{1}{2+d_{\mathbf{x}}}} d\mathbf{x}\right)^{-\frac{2}{d_{\mathbf{x}}}} \det(\mathbf{H}_{\mathbf{x}})^{-\frac{1}{2+d_{\mathbf{x}}}} \mathbf{H}_{\mathbf{x}}.
\end{align}
The error estimate on this optimal metric will be \cite{loseille2010fully}:
\begin{align}
\mathbf{E}_{\bx}^{opt}(\mathcal{M}_{\bx}^{opt}) = d_{\bx} \mathcal{C}_\bx^{-\frac{2}{d_{\bx}}} \underbrace{\left(\int_{\Omega} \det(\mathbf{H}_{\bx})^{\frac{1}{2+d_{\bx}}} \mathrm{d}\mathbf{x} \right)^{\frac{2+d_{\bx}}{d_{\bx}}}}_{\mathcal{K}_{\bx}},
    \label{opt_det_err}
\end{align}

\paragraph{Stochastic continuous error model}
Similarly, Estimate (\ref{StochErrorDiscrete}) in the continuous framework of Riemannian metric spaces, is written as:
\begin{align}
\mathbf{E}_{\bxi} ({\mathcal M}) ~=~ 
\int_{\Xi} 
  \mbox{trace} \left( \mathcal{M}^{-\frac{1}{2}}(\bxi) \, 
 {\rho_{\bxi} \cdot H(j(\bxi))} \, \mathcal{M}^{-\frac{1}{2}}(\bxi) \right)  \,
  \mathrm{d}\bxi \,   
  \label{stoch_err_est}
\end{align}
with $H(j(\bxi))$ the Hessian matrix of $j(\bxi)$ of size $d_{\bxi} \times d_{\bxi}$ .\\
\noindent
The \textit{stochastic optimisation problem} is then formulated as follows:
\begin{eqnarray}
    \mbox{Find } {\mathcal M}^{opt}_{\bxi} = 
    \argmin_{{\mathcal M}_{\bxi}} \, \mathbf{E}_{\bxi} ({\mathcal M}), \  \  \text{ subject to }\  \  \mathcal{C}(\mathbf{M}) = \mathcal{C}_{\bxi}
\label{OptPbStoch}
\end{eqnarray}
where $\mathcal{C}_{\bxi}$ denotes a specified complexity in the parameters space. This is the continuous counterpart of the number of samples $N_{\bxi}$ (or CFD computations) and corresponds to a targeted computational effort constraint.
The notation $\mathcal{M}^{opt}_{\bxi}$ holds for the optimal metric that minimises the expectation of the continuous interpolation error in the parameter space. We will use this metric to build a simplex tessellation of the parameter space, which is, roughly the mesh associated with  $\Xi$. \\
The optimal stochastic metric, solution to the optimisation problem (\ref{OptPbStoch}) is :
\begin{align}
\label{OptMstoch}
 \mathcal{M}_{\bxi}^{opt} = \mathcal{C}_{\boldsymbol{\xi}}^{\frac{2}{d_{\boldsymbol{\xi}}}} \left(\int_{\Omega_{\boldsymbol{\xi}}} \det(\rho_{\boldsymbol{\xi}}|H(j(\bxi))|)^{\frac{1}{2+d_{\boldsymbol{\xi}}}} d\boldsymbol{\xi}  \right)^{-\frac{2}{d_{\boldsymbol{\xi}}}} \det(\rho_{\boldsymbol{\xi}}|H(j(\bxi)|)^{-\frac{1}{2+d_{\boldsymbol{\xi}}}} |\rho_{\boldsymbol{\xi}} H(j(\bxi)|
\end{align}
and the error estimate on this optimal metric is given by:
\begin{align}
    \mathbf{E}_{\bxi}^{opt}(\mathcal{M}_{\bxi}^{opt}) = d_{\bxi} \mathcal{C}_{\bxi}^{-\frac{2}{d_{\bxi}}} \underbrace{\left(\int_{\Xi} \det(\rho_{\bxi} |H(j(\bxi))|)^{\frac{1}{2+d_{\bxi}}} \mathrm{d}\bxi \right)^{\frac{2+d_{\bxi}}{d_{\bxi}}}}_{\mathcal{K}_\xi}
    \label{opt_stoch_err}
\end{align}
The proofs for formulations (\ref{OptMstoch}) and (\ref{opt_stoch_err}) are included in \cite{langenhove2017these}.\\

{\it Practical remark}: In practice, the continuous (exact) states in formulation (\ref{OptMdet}) {and (\ref{OptMstoch})} are approximated with discrete state. Hence, the adjoint state $w^*$ is replaced by a discrete adjoint state, and for the gradients and hessians computations we use derivative recovery methods such as $L^2-$projection or Green formula (see \cite{alauzet2010high}).
\section{Adaptive Strategies for error control}
We have shown in the previous sections how we formulate the adaptation problem in the continuous framework of the Riemannian metric space, and how we solve the optimization problem in this framework, for both stochastic and deterministic systems. A bijection between the continuous and the discrete framework of the discretized mesh exists as it has been detailed in \cite{loseille2008these,belme2011these}. The optimal metric, solution of the optimisation problem, is computed and used by the discrete mesh generator, named {\tt Feflo.a} \cite{loseille2010anisotropic}, to build the new anisotropic adapted mesh where the errors are controlled. \\

In this section, we first briefly recall the mesh adaptation algorithm for the deterministic space. We then introduce the retained adaptive strategy for the approximation error control in the parameters space. Our approach is here facilitated by the choice of stochastic approximation space.
Finally, we propose a coupled algorithmic approach for a more optimal control of {\em both} deterministic and stochastic error contributions.

\subsection{Deterministic adaptive strategy}
\begin{figure}[!h]
\centering
\includegraphics[scale=0.5]{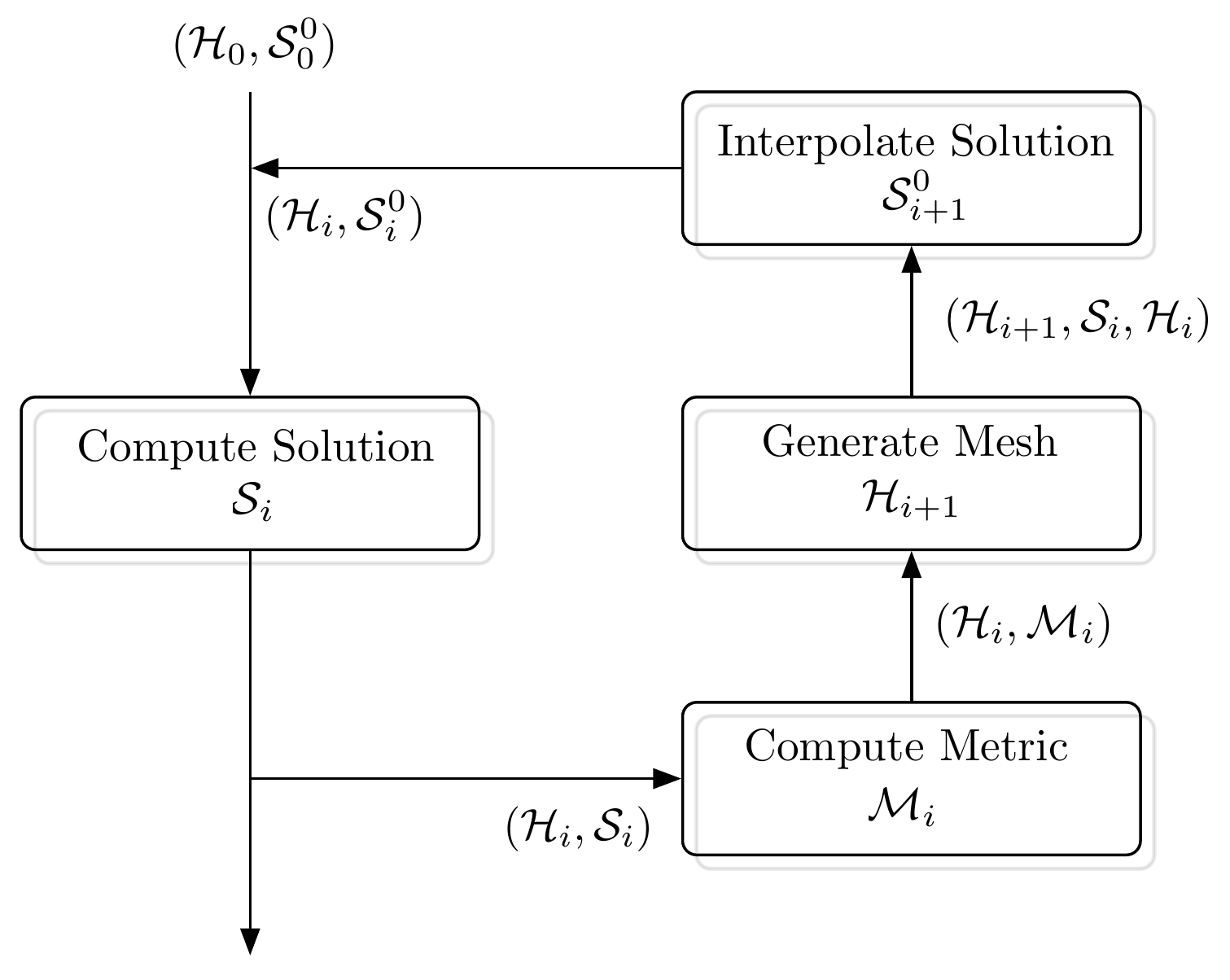}
\caption{Schematic illustrating the mesh adaptation process. The couple (mesh, solution) is denoted $(\mathcal{H}_i, \mathcal{S}_i)$ where the subscript $i$ denotes the fixed-point iteration number; $\mathcal{S}^0_0$ is the initial solution on mesh $\mathcal{H}_0$ whereas $\mathcal{S}^0_{i}$ is the initial solution interpolated on the mesh $\mathcal{H}_{i}$. Each discrete mesh $\mathcal{H}_i$ is based on the metric $\mathcal{M}_i$.}
\label{algoPtfx}
\end{figure}

Regarding deterministic mesh adaptivity, we follow the work of \cite{loseille2008these,belme2011these} and employ a fixed point algorithm as illustrated in Figure (\ref{algoPtfx}). Indeed, the mesh adaptation problem is a non-linear problem, and an iterative algorithm is well suited to converge the couple mesh-solution. The stopping criteria can be a targeted error lever or, as is usually proposed in practice, a maximum number of iterative loops. In general, $5$ fixed-point iterations are enough to reach a satisfying level of convergence. 

\subsection{Stochastic adaptive strategy}
We underline here that we rely on non-intrusive methods for the approximation of the stochastic problem, meaning that the deterministic solver counterpart (e.g. the CFD code) is seen as a black box. Many different choices of non-intrusive stochastic surrogate models are available prior to the sampling adaptation. However, for our type of applications, it is wise to rely on a -- {\em robust}, i.e. relatively insensitive to the response surface smoothness, and -- {\em flexible} and efficient approach, i.e. with an easy handling of local parametric sampling avoiding tensor-based refinements. The Simplex-Stochastic Collocation (SSC) method is a good candidate that meets these requirements by discretizing the domain into simplices, and naturally enforces local extremum conservation to suppress unphysical oscillations.\\
A simple piecewise-linear approximation is chosen over higher degree polynomial approximation in accordance with the error estimate chosen to drive the mesh adaptation in the parameter space. Indeed, the error estimate is based on the interpolation error committed when linearly interpolating a quadratic function.
Relevant information about the stochastic approximation based on the simplex elements we use can be found in \ref{ssc_elem}.

~\\
The adaptive response surface in the parameter space is generated using similar numerical tools to the deterministic ones. Barring a few differences, a fixed-point algorithm is also used here. The optimisation problem (\ref{OptPbStoch}) is solved in the Riemannian metric space and a stochastic metric is proposed. One difference resides in the fact that we do not erase the topology of the parameters space discretization at each fixed point iteration. Instead, we make use of it for computational cost purposes.
In the stochastic case, since each mesh vertex truly represents a sample requiring a (potentially) costly deterministic numerical simulation, one would like to minimize the total number of samples used over the whole adaptation procedure. Therefore samples acquired along the adaptation steps are always kept in the final discretization. In this case, the cumulative sum of mesh vertices from the initial design up to the final mesh is relevant to the numerical cost of the approximation method. 
The adaptation process used to control the stochastic error $\bar{\eta}$ is detailed in Algorithm \ref{it_stoch_algo}. 
The following notations have been used: a total number of adaptive loops $n_{adap}$ is fixed and the iteration count index is the subscript $l$. 
$\mathcal{C}_{\bxi}$ is the target complexity of the optimisation problems set on the parametric domain. This target complexity is passed on to the meshing code which, based on the computed optimal metric $\mathcal{M}_{\bxi}^{opt}$ builds the discrete optimal mesh.

\begin{algorithm}
\caption{Adaptation of the stochastic discretization}
\label{it_stoch_algo}
\begin{algorithmic}
    \State Generate (or read) initial samples $\{\bxi\}_0$ and form initial mesh $\mathcal{H}_{\xi,0}$ by a Delaunay triangulation.
    \State Compute $j(\{\bxi\}_0)$ for this DoE. 
    \For{ $l=1$ to $n_{adap}$}
    \State Compute optimal metric $\mathcal{M}^{opt}_{\bxi,l}$ for complexity value $\mathcal{C}_{\bxi}$ based on {numerical solution approximation constructed on} previous mesh $\mathcal{H}_{\bxi,l-1}$. 
    \State Generate new mesh $\mathcal{H}_{\bxi,l}$ from $\mathcal{M}_{\bxi,l}^{opt}$ containing $N_{\bxi,l} = N_{\bxi,l-1}+N_{\bxi,l}^{new}$ samples.
    \State {Run deterministic solver to} compute QoI at the $N_{\bxi,l}^{new}$ new samples and update $j(\{\bxi\}_l)$. 
    \State Compute the statistical moments of the surrogate model.
    \EndFor
\end{algorithmic}
\end{algorithm}
As outlined in Algorithm \ref{it_stoch_algo}, we start out by constructing an initial mesh. This initial mesh is a Delaunay triangulation of the random samples drawn from $\rho_{\bxi}$. Once the QoI is computed for each sample, a surrogate model of the response is constructed. Next, we compute the optimal metric and the associated error estimate (\ref{stoch_err_est}) using {\tt Metrix} \cite{alauzet2009metrix}. The new, adapted mesh is generated from $\mathcal{M}^{opt}_{\bxi}$  using the adaptive local mesher {\tt Feflo.a} \cite{loseille2010anisotropic}.\\

\subsection{Optimal adaptive strategy}
\label{optStrategy}
We have seen up until now the adaptive strategies and algorithms when solving the optimisation problems (\ref{OptPbDet}) or (\ref{OptPbStoch}). However, as part of the motivation and novelty of this paper we wish to: (a) quantify the deterministic error at each sample of the parameters space and compute a local (sample-wise) optimization deterministic problem with a required complexity $\mathcal{C}_{\bx}$ and to (b) be able to choose which error, stochastic or deterministic, dominates a computation and thus to solve the corresponding problem.\\
There are several options for the first issue (a). Indeed, the two sources of errors: stochastic and deterministic, are strongly coupled. Suppose we compute our QoI for a sample $\bxi_{(i)}$ on a uniform mesh. This mesh is usually not fit to our QoI but to the studied problem in general. Thus, very often, large deterministic errors can propagate to the parameter space. This is even more pronounced when dealing with problems involving shocks. However, we are now able to build an adapted mesh to best observe our QoI and thus reduce the propagation of the deterministic error to the stochastic space. Moreover, the deterministic error on each sample does not necessarily contribute to the overall problem with the same level of error.
Ideally, one would want the mean deterministic error to be lower than or equal to some target error $\bar{\varepsilon}^t$, best chosen to be comparable to the surrogate model error $\bar{\eta}$. Secondly, one would want the variance of the error contained in all the sample to be as close as possible to zero. This ensures that all the samples used to construct the surrogate model of the QoI are of equal accuracy. Since an interpolation method is used for the surrogate model, having a large error in a few samples can be detrimental to the quality of the surrogate model.
A trivial solution is to require that the deterministic error in each sample $\varepsilon (\bxi_{(i)})$ be equal to $\bar{\varepsilon}^t$. This results in the expected deterministic error to be equal to $\bar{\varepsilon}^t$ and its variance to be (as close as possible to) zero.\\
Following the error model (\ref{opt_det_err}) the required complexity for this case can be computed as:
\begin{align}
    \mathcal{C}_{\bx}^{t} = \left(\frac{d_{\bx} \mathcal{K}_{\bx}}{\bar{\varepsilon}^{t}}\right)^{\frac{d_{\bx}}{2}}.
    \label{det_goal_compl}
\end{align}
using the notation $\mathcal{K}_{\bx}$ introduced in (\ref{opt_det_err}) and we recall $d_{\bx}$ is the deterministic space dimension.\\
The detailed optimisation algorithm is presented in Algorithm \ref{err_control_algo} hereafter.

\begin{algorithm}
\caption{Sample-wise control of the deterministic error contribution over the parametric domain}
\label{err_control_algo}
\begin{algorithmic}
    \State Compute mean deterministic error over the $N_{\bxi}$ samples of the parameter space $\bar{\varepsilon} = \mathbb{E}[\varepsilon (\bxi)] = \int_{\Xi} \varepsilon(\bxi) \rho_\xi \mathrm{d} \bxi$.
    \For{$i=1$ to $N_{\bxi}$}
    \If{$\varepsilon(\bxi_{(i)}) > \bar{\varepsilon}^t$} 
    \State Compute required complexity $\mathcal{C}^{t}_{\bx_{(i)}}$ from (\ref{det_goal_compl}) to reduce deterministic error.
   \State Solve optimisation problem (\ref{OptPbDet}) associated to this computed complexity.
   \State Update $j(\bxi_{(i)})$ value to sample point $\bxi_{(i)}$ in the parametric domain.
    \Else
    \State Keep $j(\bxi_{(i)})$ value.
    \EndIf
    \EndFor
\end{algorithmic}
\end{algorithm}

 This optimisation strategy can be coupled with Algorithm \ref{it_stoch_algo} in order to adress issue (b) and control the total error by performing adaptations in both physical and parameter spaces. The required complexity $\mathcal{C}_{\bxi}^t$ (equivalent to a number of $N_{\bxi}^{t}$ samples in the discrete parameter space) for the stochastic problem will be computed following a similar approach from (\ref{opt_stoch_err}) :
\begin{align}
    \mathcal{C}_{\bxi}^{t} = \left(\frac{d_{\bxi} \mathcal{K}_{\bxi}}{\bar{\eta}^{t}}\right)^{\frac{d_{\bxi}}{2}},
    \label{stoch_goal_compl}
\end{align}
The coupled adaptation strategy is detailed below in Algorithm \ref{opt_err_control_algo}.
\begin{algorithm}
\caption{Total error control strategy}
\label{opt_err_control_algo}
\begin{algorithmic}
    \State Generate and compute (or read) initial samples $\{\bxi\}_0$ on uniform or adapted initial deterministic mesh $\mathcal{H}_{\bx,0}$.
    \State From initial stochastic mesh $\mathcal{H}_{\bxi,0}$ by a Delaunay triangulation (or load initial existing stochastic mesh). 
    \State Compute the stochastic error $\bar{\eta}_0$ and the mean deterministic error $\bar{\varepsilon}_0$.
    \State Set maximal number of iteration cycles $it_{MAX}$ and set total target error value $\bar{\delta j}^{t}$.
    \While{$l < it_{MAX}$ \bf{or} $\bar{\delta j} \leq \bar{\delta j}^{t}$}
    \If{$\bar{\varepsilon}_l > \bar{\eta}_l$}
    \State Adapt deterministic computations following Algorithm \ref{err_control_algo} with $\bar{\varepsilon}^{t} = \bar{\eta}_l$.
    \Else
    \State Adapt in parametric domain following Algorithm \ref{it_stoch_algo}. 
    \EndIf
    \State $l=l+1$
    \EndWhile
\end{algorithmic}
\end{algorithm}

\section{Numerical applications}

In the previous sections, error estimates for adaptive control of stochastic and deterministic approximation errors have been proposed and a continuous metric-based approach is used to build anisotropic adapted meshes in both domains.  
In this section the effectiveness of the proposed approach will be demonstrated on several test cases. In particular, we will underline first the effectiveness of the adaptive approximation in the stochastic space where we emphasize also the impact of the parameters probability density functions for sensitive nonlinear functionals: first, on a discontinuous analytical function and then on the fluid mechanics stochastic piston problem. The coupled approach where both deterministic and stochastic error are controlled is demonstrated on a supersonic/hypersonic inlet problem.

\subsection{Validation of stochastic test functions}

The test function treated here is one with multiple curved and straight discontinuities proposed by Jakeman \textit{et al.} in \cite{jakeman2013minimal}.  We consider thus a function $y$ defined on $[-1, 1]^{d_{\mathbf{\bxi}}}$ and of analyical form given hereafter:
\begin{align}
    y({\bxi}) = \left\{
    \begin{array}{ll}
        f_1({\bxi}) - 2 \  &\text{ if }\   3\xi_1 + 2\xi_2 \geq 0 \  \text{ and }\  -\xi_1 + 0.3\xi_2 < 0,\\
        2f_2({\bxi}) \  &\text{ if }\   3\xi_1 + 2\xi_2 \geq 0 \  \text{ and }\  -\xi_1 + 0.3\xi_2 \geq 0,\\
        2f_1({\bxi}) + 4 \  &\text{ if }\   (\xi_1+1)^2+(\xi_2+1)^2 < 0.95^2 \  \text{ and }\  d_{{\bxi}} = 2,\\
        f_1({\bxi}) \  &\text{ otherwise}.
    \end{array}
\right.
\label{jakemanTestRef}
\end{align}
where $f_1$ and $f_2$ are given by
\begin{align*}
    \begin{split}
        f_1({\bxi}) = \text{exp}\left(-\sum_{i=1}^2 \xi_i^2\right) - \xi_1^3 - \xi_2^3, \quad \text{and} \quad
        f_2({\bxi}) = 1+f_1(\xi) + \frac{1}{4 d_{\bxi}} \sum_{i=2}^{d_{{\bxi}}} \xi_i^2.
    \end{split}
\end{align*}

For a graphical representation of $y$ we refer to \cite{jakeman2013minimal}.
\subsubsection{Test 1: uniform probability density function}

\begin{figure}
    \centering
    \includegraphics[width=0.75\linewidth]{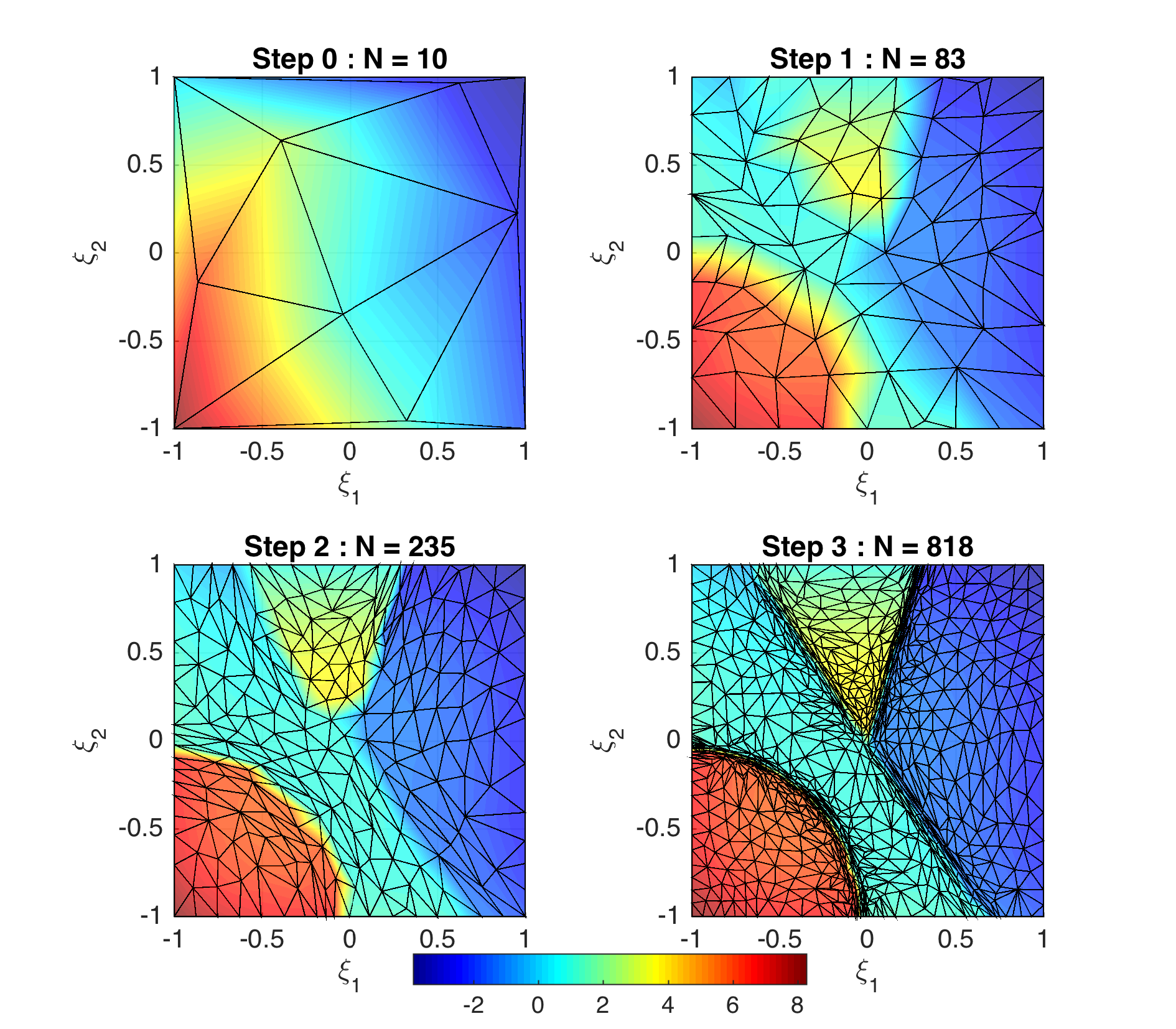}
    \caption{Test 1: isocontours of the approximated solution of (\ref{jakemanTestRef}) on the mesh at successive refinement levels.}
    \label{fig:disco_sol}
\end{figure}

Figure \ref{fig:disco_sol} shows the results of the adaptive approximations of the discontinuous functional (\ref{jakemanTestRef}) in two dimensions,  when both $\xi_1$ and $\xi_2$ follow a uniform distribution $\mathcal{U}_{[-1, 1]}$. The meshes are displayed here at different refinement steps/cycles.
The initial mesh shown in Figure \ref{fig:disco_sol} (top left) contains $10$ samples following a Latin Hypercube Sampling (LHS) DoE, with very little information as to the whereabouts of the discontinuities.
Since for the approximation of the response surface a first degree Newton-Cotes quadrature is used, one only computes samples on the vertices of the simplices. The response surface inside each element is approximated by a linear interpolation on the simplex.
Three subsequent adaptation steps were performed. 
As the mesh becomes more refined, one can observe that the characteristics, notably the discontinuities, of the response surface are better captured. 
The sample point density is increased in the vicinity of the discontinuity and the elements become stretched as to be aligned with these discontinuities.
Note that no discontinuity capturing technique is used here, nor are there any parameters that need to be set; the results are obtained by applying the metric-based adaptive method through Algorithm \ref{it_stoch_algo}.\\

\begin{figure}
    \centering
    \begin{subfigure}{.49\linewidth}
        \centering
        \includegraphics[width=1.0\linewidth]{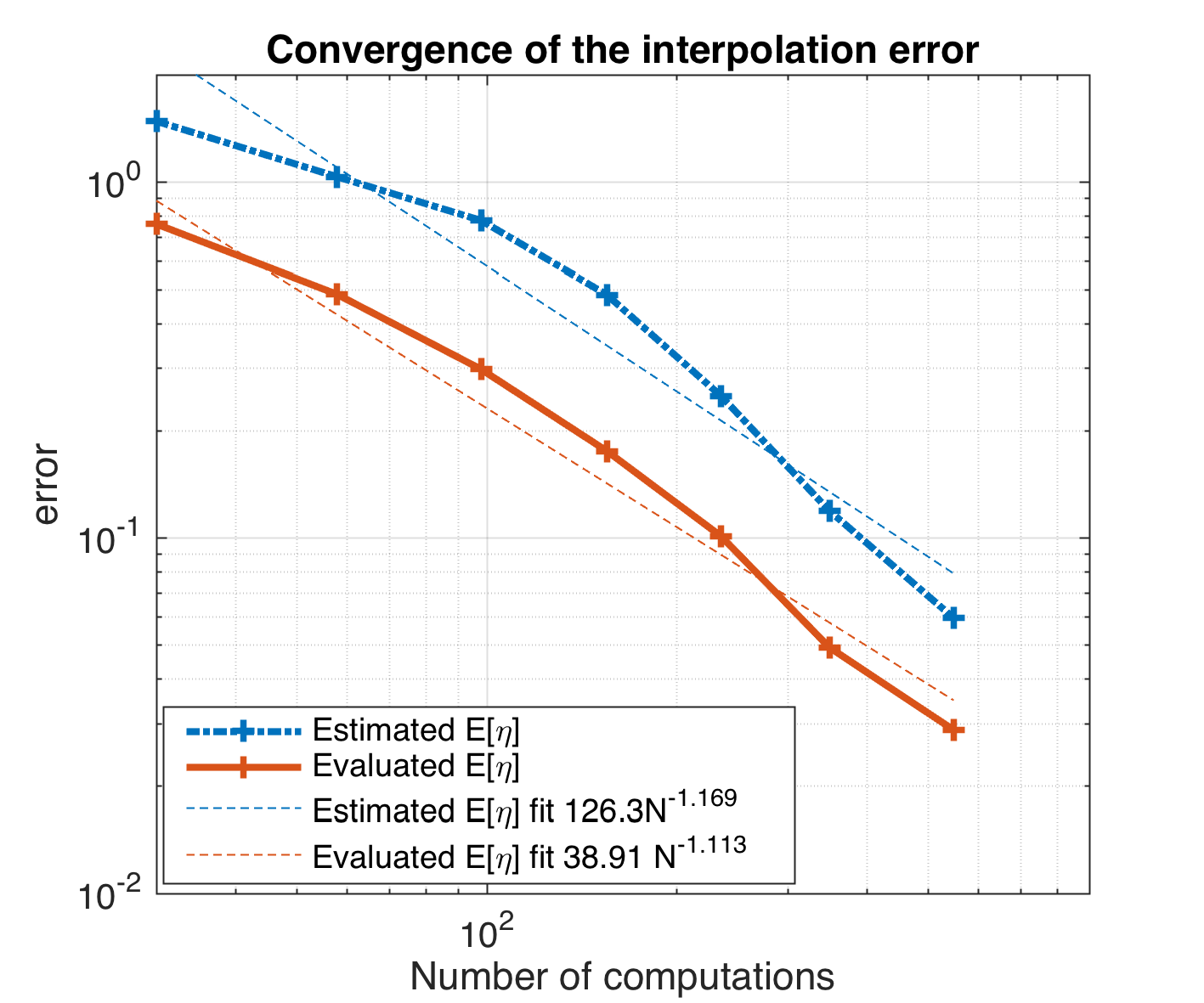}
        \caption{8 refinement steps}
        \label{fig:disco_8_conv}
    \end{subfigure}
    \begin{subfigure}{.49\linewidth}
        \centering
        \includegraphics[width=1.0\linewidth]{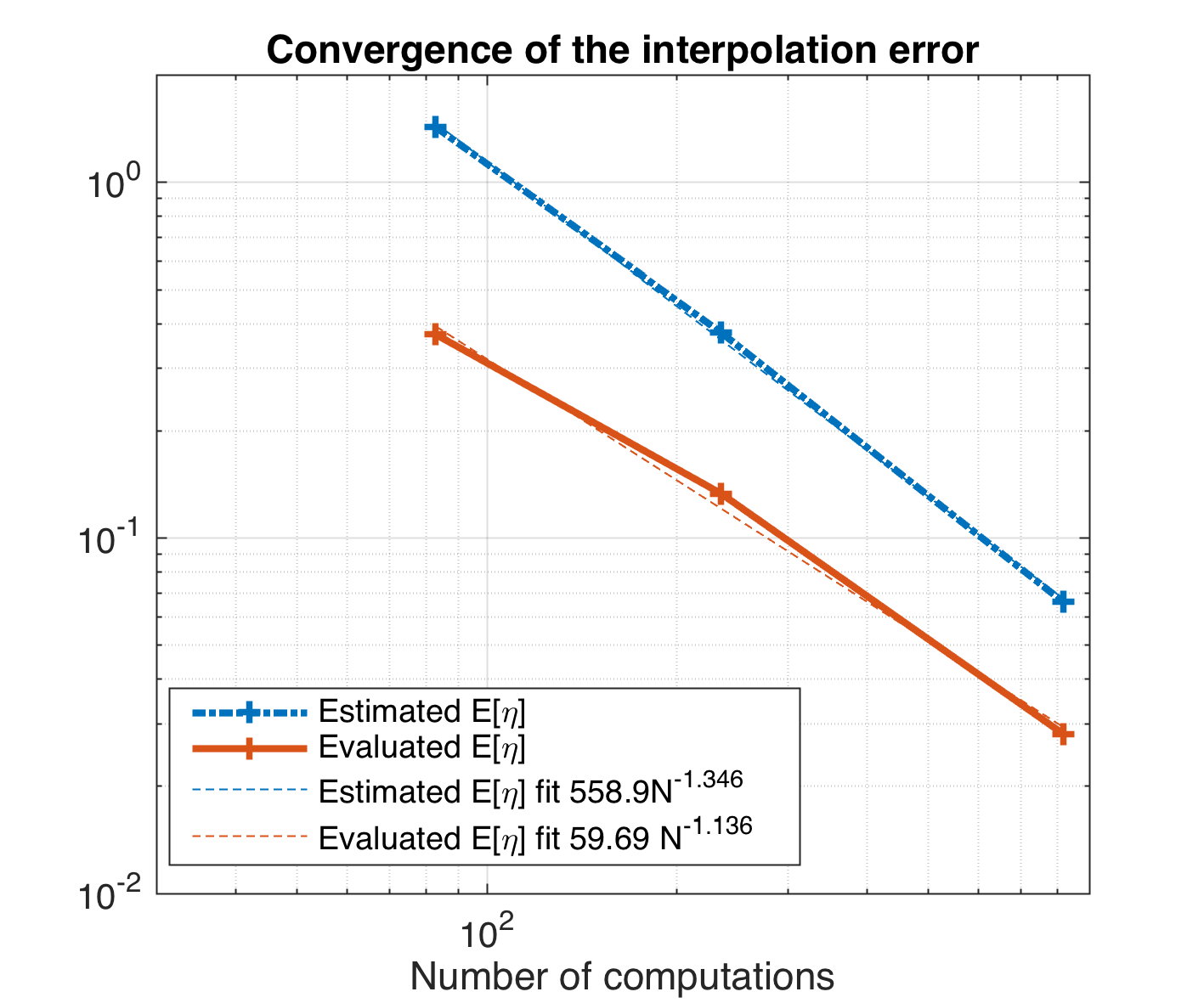}
        \caption{3 refinement steps}
        \label{fig:disco_3_conv}
    \end{subfigure}
    \caption{Test 1: $L^1-$convergence of the interpolation error and the error estimate. The interpolation error is computed on a third degree Newton-Cotes quadrature, a sort of subgrid constructed in order to evaluate the interpolation error. The blue dot-dash line represents the expectation of the estimated optimal continuous interpolation error while the red full line shows the expectation of the interpolation error evaluated on the fine quadrature within each element. For both curves a least-squares fit to an exponential curve was computed, these are shown as the dashed lines.}
    \label{fig:disco_conv}
\end{figure}
Since this is an analytic test function, for which we know the exact solution, once the adaptation is completed one can easily compute the expectation of the error committed by the response surface approximation relative to the exact solution. This is what we call the ``evaluated" error in the plots vs. the ``estimated" error, the latter being provided by the error estimation and the algorithm and not involving the exact reference solution.  
The convergence of the evaluated and estimated errors are shown in Figure \ref{fig:disco_conv} for two separate cases. Both results illustrated in figures (\ref{fig:disco_3_conv}) and (\ref{fig:disco_8_conv}) started with the same initial mesh containing $10$ samples. The difference is that in Figure \ref{fig:disco_8_conv} $8$ small refinement steps were used (at each step the targeted complexity is doubled) while in Figure \ref{fig:disco_3_conv} only $3$ steps were taken (at each step the mesh complexity was increased by a factor of $5.5$). This allows to evaluate the influence of the step size on the convergence of the error.

The theoretical error convergence is given by $\Ocal ( N_{\bxi}^{-{\frac{2}{d_{{\bxi}}}}} )$ (see \cite{LoseilleAIAA2007}). 
It can be seen in Figure \ref{fig:disco_conv} that the fits provide convergence rates that are close but superior to the theoretical one. Furthermore, convergence rates are not significantly altered when the adaptation step size is significantly increased.
Despite the fact that the theoretical convergence results were based on the linear interpolation of a quadratic form, a second order convergence is maintained even for this discontinuous function underlining the effectiveness of the proposed method for discontinuous response surfaces.
The expectation of the \emph{exact interpolation error} shown in Figure \ref{fig:disco_conv} is computed on a third degree Newton-Cotes quadrature within each element. 
The convergence constant will both depend on the function that is being approximated and on the discrepancy between the continuous number of mesh vertices, i.e the complexity $\mathcal{C}$ and the number of vertices in the discrete mesh $N$. In the computation of the continuous interpolation error estimate, the continuous number of vertices is used, while in Figure \ref{fig:disco_conv} it is plotted as a function of the realized number of vertices in the discrete mesh $N$. For a detailed numerical validation between the continuous error estimate and the actual interpolation error, the reader is referred to \cite{loseille2011continuous2}. \\
{Another interesting point here is the evolution of the error with the choice of the number of refinement steps of the adaptive process. At the end of the process, we observe that we reach the same error level with both strategies but the computational cost is not the same: it has required $800$ CFD computations with $3$ refinement steps, compared to $550$ computations with $8$ refinement steps. Increasing moderately the number of additional simulation runs at each step, therefore seems a better strategy. 
This reduces the number of  computations to reach a given threshold. Such a strategy is peculiar to the stochastic approximation process and is not needed for standard deterministic CFD mesh adaptations for which all nodes maybe redistributed and enriched.}

\subsubsection{Test 2: piecewise-uniform probability density function}
The same test function is used, but this time the underlying probability density function is now discontinuous and piecewise-constant, i.e. : 
\begin{align}
    \rho_{\bxi} = \left\{
    \begin{array}{ll}
        \frac{1-(2.6-\frac{1}{2}0.2^2 \pi) 0.005 - (\frac{1}{2} 0.2^2 \pi) 0.9}{1.4}\approx 0.66\  &  \text{ if }\  \xi_2 \geq -0.3\xi_1 + 0.3,\\
        0.9 \  &  \text{ if } \xi_1^2 + (\xi_2+1)^2 \leq 0.2^2,\\
        0.005 \  &  \text{ otherwise.}
    \end{array}
\right.
\label{discoRhoDef}
\end{align}
The method presented in this manuscript does not require an orthogonal basis to be found with respect to the local probability measure. However, the SSC approximation must  take into account the pdf, the effect of which shows in the quadrature weights. Here, it  makes a numerically challenging case as the pdf is discontinuous. 
In order to compute the weight associated to each quadrature point, the linear Lagrange polynomial associated to that sample and weighted by the pdf is numerically integrated on the local simplex.
Once the tessellation is given, that step can be done offline and only involves the interrogation of the pdf on a finer subgrid within each element.
Details about this procedure are presented in \ref{ssc_elem}.
Table \ref{tab:subquad_effect} collects the sum of the quadrature weights on the DoEs at different refinement steps for different subgrid quadratures; disregarding numerical approximations this sum should exactly be equal to $1$.
\begin{table} 
    \footnotesize
    \arraycolsep=0.5pt
    \medmuskip=0.3mu
\centering
\begin{tabular}{| l || c | c | c | c |}
\hline 
Subgrid NC deg.  & $N_{\bxi} = 10$ & $N_{\bxi} = 77$ & $N_{\bxi} = 218$ & $N_{\bxi} = 707$\\ 
 \hline 
\hline 
$1$  & $1.209$ & $1.195$ & $1.024$ & $0.996$\\ 
 \hline 
$2$  & $0.589$ & $0.986$ & $0.984$ & $1.007$\\ 
 \hline 
$3$  & $1.199$ & $0.984$ & $1.003$ & $1.002$\\ 
 \hline 
$4$  & $0.965$ & $0.980$ & $1.004$ & $1.002$\\ 
 \hline 
$5$  & $1.012$ & $0.999$ & $1.001$ & $0.999$\\ 
 \hline 
$6$  & $0.856$ & $1.002$ & $0.993$ & $1.000$\\ 
 \hline 
$7$  & $0.965$ & $1.002$ & $1.006$ & $1.000$\\ 
 \hline 
$8$  & $1.073$ & $1.003$ & $1.000$ & $1.000$\\ 
 \hline 
\end{tabular}
\caption{Overview of the sum of the quadrature weights for the discontinuous pdf (\ref{discoRhoDef}). The subgrid NC quadrature degree used to compute the weights according to (\ref{c_def}) is indicated in the column on the left.}
\label{tab:subquad_effect}
\end{table}
One can see that the error committed is more significant for meshes with very few samples. However, with increasing quadrature degree and increasing $N_{\bxi}$, this error quickly drops to $\mathcal{O}(10^{-3})$ or less. 
We emphasize that for smooth pdf the error is much lower.
Further discussions about the numerical conditioning of higher-order NC quadratures are also given in \cite{langenhove2017these}.\\

\begin{figure}
    \centering
    \includegraphics[width=0.6\textwidth]{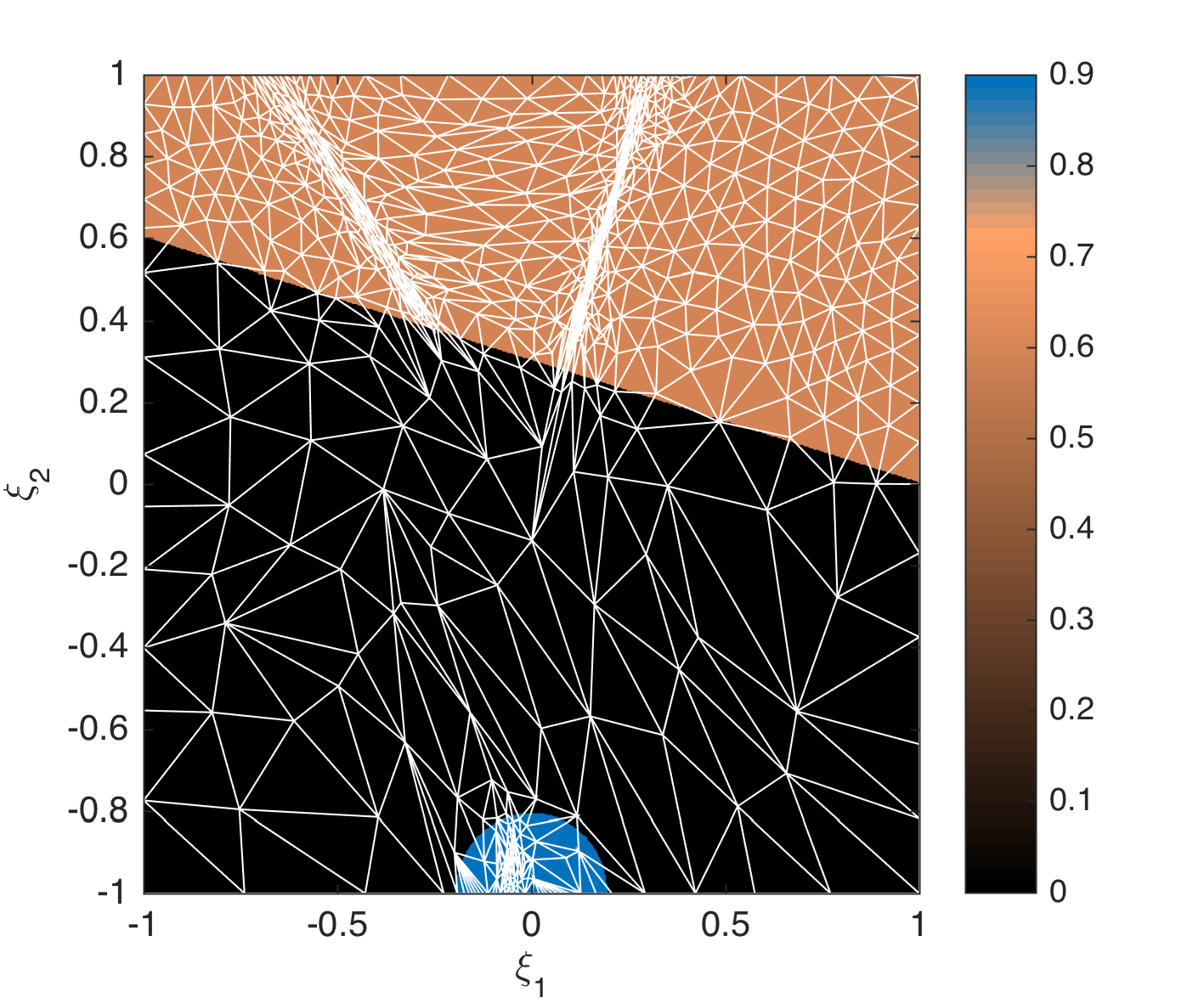}
    \caption{Test 2: the adapted mesh plotted is obtained at the end of three iteration steps and for which the approximation errors are shown in Figure \ref{fig:discoDisco_3_conv}. The color background represents the discontinuous probability density function.}
    \label{disco_rho}
\end{figure}

In Figure \ref{disco_rho} the color contours visualise the discontinuous pdf while the overlaid mesh shown is the result of three subsequent adaptation steps. 

The resulting mesh clearly shows how the probability distribution is taken into account in the mesh refinement. Indeed, as we have shown in Section 2, the probability density function acts as a weight of the interpolation error. Thus, despite the presence of singularities, the regions of low to zero event probability correspond to coarse regions where very few samples will be placed, while in regions with high event probability more sample points will be placed and the sample point density is especially high in the vicinity of the discontinuities.

\begin{figure}
    \centering
    \begin{subfigure}{.49\linewidth}
        \centering
        \includegraphics[width=1\linewidth]{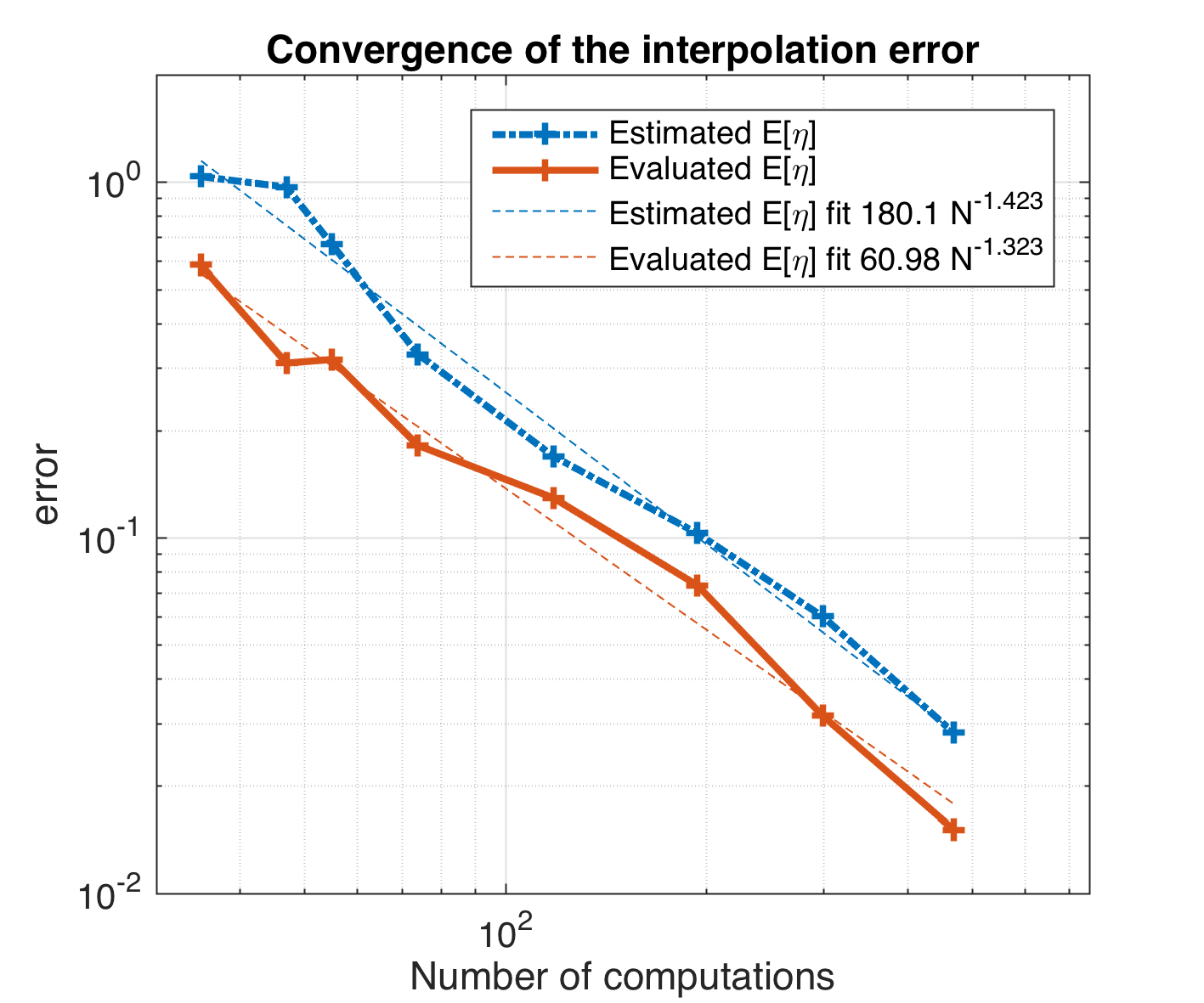}
        \caption{8 refinement steps}
        \label{fig:discoDisco_8_conv}
    \end{subfigure}
    \begin{subfigure}{.49\linewidth}
        \centering
        \includegraphics[width=1\linewidth]{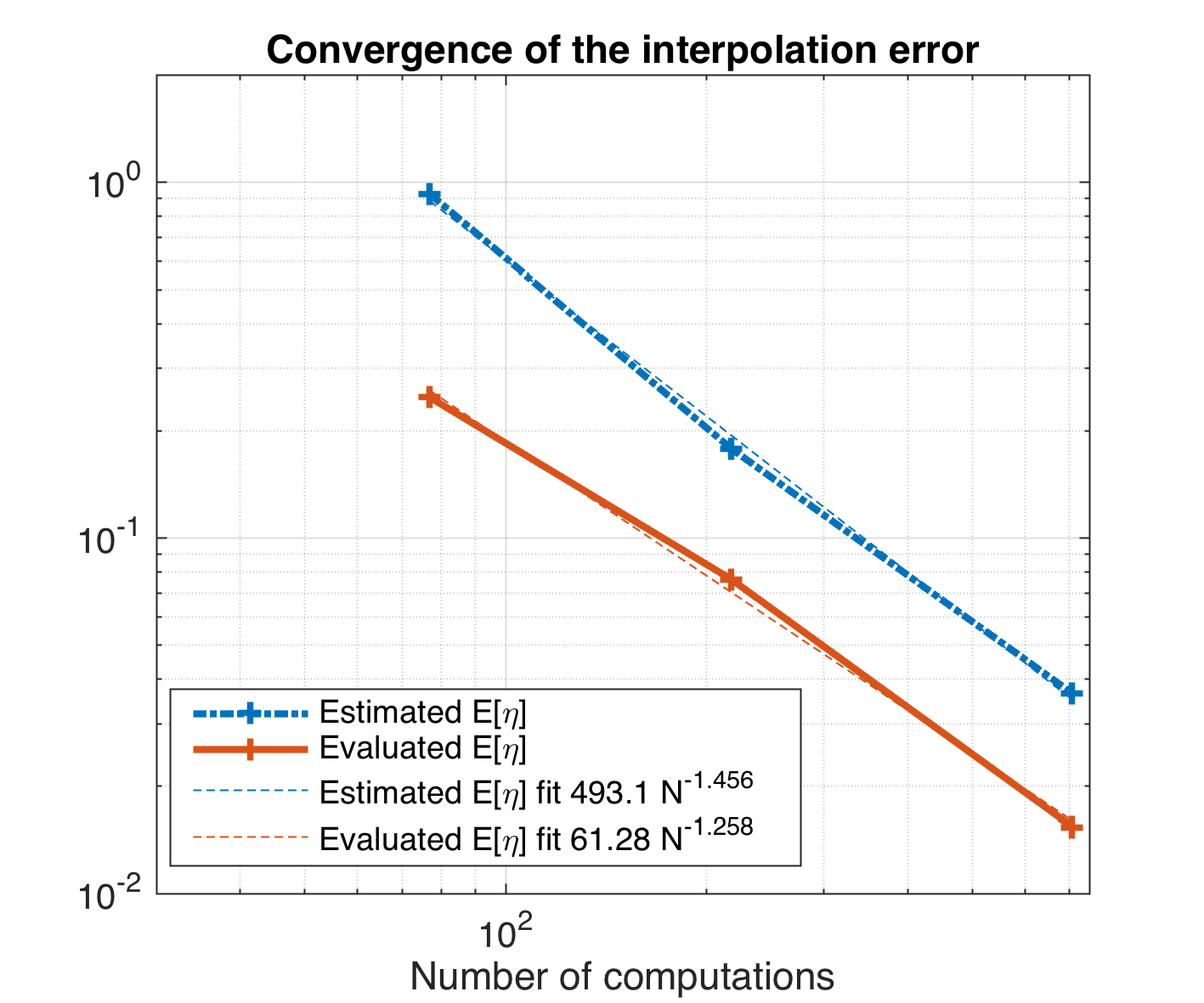}
        \caption{3 refinement steps}
        \label{fig:discoDisco_3_conv}
    \end{subfigure}
    \caption{Test 2: $L^1-$convergence of the interpolation error and the error estimate. The interpolation error is computed on a third degree Newton-Cotes quadrature, a sort of subgrid constructed in order to evaluate the real interpolation error.}
    \label{fig:discoDisco_conv}
\end{figure}
Similarly to the previous case, the adaptation process was executed with different step sizes. The resulting convergence of both the expectation of the actual interpolation error and the expectation of the error estimate, shown in Figure \ref{fig:discoDisco_conv}, exceed $2^{\text{nd}}$-order convergence even though the test function is highly discontinuous. The observations made for the case of uniform distribution are confirmed here: the convergence rate is at least what was predicted by theory and is not significantly affected when the step size is  increased nor does the discontinuous probability density function negatively affect this convergence rate. The convergence constant of the error estimate does increase when the step size is increased.\\
{Moreover, we again observe that less samples (here $375$ vs. $600$) are required to reach the same accuracy when more refinement steps are performed, showing that a slow increase is a better strategy.}
\subsubsection{Test 3: uniform distributions in three dimensions}
\begin{figure}
    \centering
    \begin{subfigure}{.49\linewidth}
        \centering
        \includegraphics[width=1\linewidth]{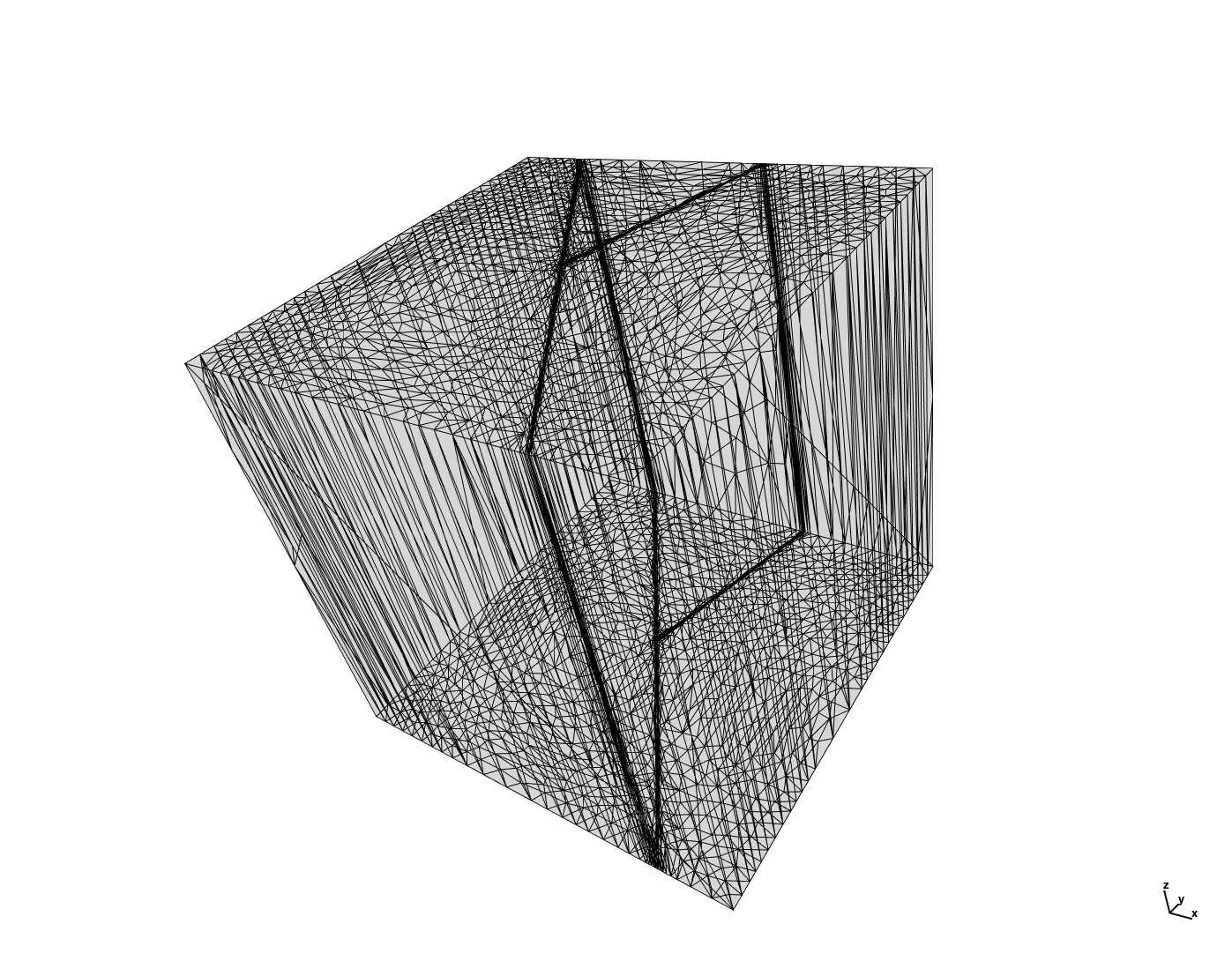}
        \caption{Test 3: Wireframe of the adapted stochastic domain}
        \label{fig:disco3d_cut}
    \end{subfigure}
    \begin{subfigure}{.49\linewidth}
        \centering
        \includegraphics[width=1\linewidth]{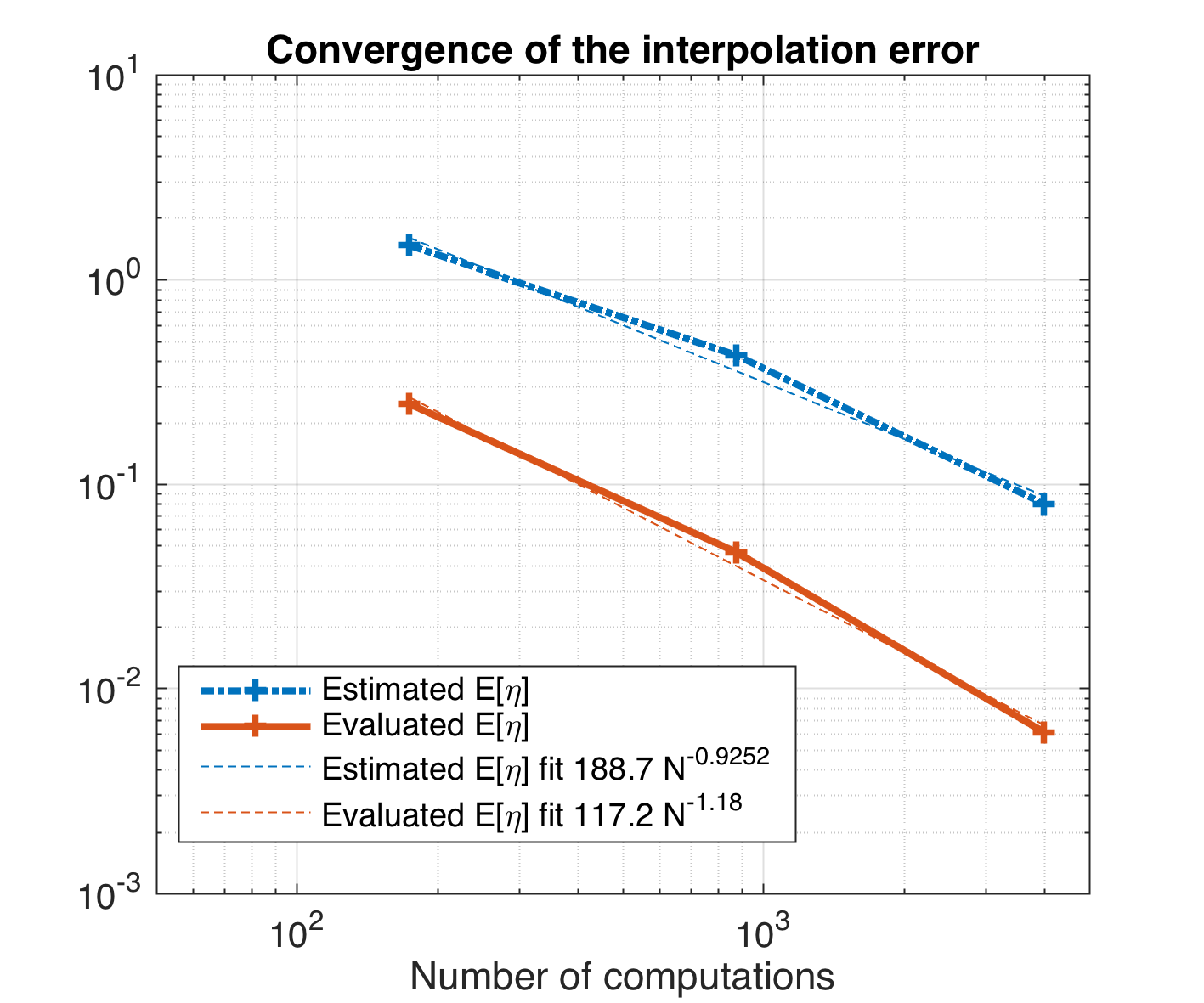}
        \caption{Test 3: convergence of the interpolation error and the error estimate with three uncertain variables and three refinement steps. The interpolation error is computed on a third degree Newton-Cotes quadrature, a sort of subgrid constructed in order to evaluate the real interpolation error.}
        \label{fig:disco3d_3_conv}
    \end{subfigure}
\caption{Test 3: $3D$ problem. Illustration of the resulting adaptive mesh (a) and the error convergence plot (b) after three adaptive steps.}
    \label{fig:disco3d_conv}
\end{figure}

The case of three uncertain variables is also examined. We consider again uniform distribution for the three variables and the resulting wireframe of the $3D$ mesh is shown in Figure \ref{fig:disco3d_cut}. As before, the discontinuities are well captured. The error convergence is examined as well and the convergence plot obtained using three adaptive steps are shown in  Figure \ref{fig:disco3d_3_conv}.

The theoretical convergence should be of $O(N_{\bxi}^{-\frac{2}{3}}) \approx O(N_{\bxi}^{-0.67})$ and it can be observed that the convergence rate obtained surpasses this prediction. In contrast to the $2D$ case, the expectation of the evaluated interpolation error decreases faster than the error estimate. This could be due to the fact that in the $3D$ version of this test function the solution is mostly discontinuous w.r.t. the first and second dimension but not in the third. It can be seen from (\ref{jakemanTestRef}) that the test function depends quadratically on $\xi_3$ whilst the dependence on $\xi_1$ and $\xi_2$ is discontinuous. As this somewhat reduces the effective dimensionality of the problem it seems that the evaluated expectation of the interpolation error diminishes quicker than the expectation of the interpolation error estimate.

\subsection{The fluid mechanics stochastic piston problem}
\label{subsec:piston}
The metric-based stochastic mesh adaptation is tested for a classical application of fluid mechanics: the piston problem. This fundamental problem was revisited several times in the context of uncertainty quantification, e.g. \cite{LinSuKar05PNAS,Zhang_2013}. A particular version was well described by \cite{witteveen2009adaptive} and approximated with simplex elements. In the following, we use the same notations as in that reference. The setup consists of a tube filled with air, assumed to be an ideal gas. A fast piston pushes the fluid from left to right in the tube. The flow domain is one-dimensional and the piston starts moving to the right  with constant velocity $u_{\text{piston}} > 0$ at the initial time $t=0$. The flow initial conditions in the tube are described by the initial density $\rho_{\text{pre}}$, pressure $p_{\text{pre}}$ and velocity $u_{\text{pre}} = 0$. 
As the piston moves to the right a shock forms ahead and moves forward into the gas with speed $u_{\text{shock}}$. The conditions behind the shock are given by $\rho_{\text{post}}$, $p_{\text{post}}$ and $u_{\text{post}}=u_{\text{shock}}$. The effects of viscosity are neglected and so the pressure behind the shock can be obtained from the Rankine-Hugoniot relation:
\begin{align*}
    p_{\text{post}} - p_{\text{pre}} = \rho_{\text{pre}} c_{\text{pre}} (u_{\text{post}}-u_{\text{pre}}) \sqrt{1+\frac{\gamma-1}{2 \gamma} \frac{p_{\text{post}}-p_{\text{pre}}}{p_{\text{pre}}}},
\end{align*}
where the initial sound speed is given by $c_{\text{pre}} = \sqrt{\gamma p_{\text{pre}}/\rho_{\text{pre}}}$ and the specific heat ratio is chosen as $\gamma=1.4$.
The Mach number of the shock is obtained from one-dimensional shock wave relations \cite{anderson2001fundamentals}:
\begin{align*}
    M_{\text{shock}} = \sqrt{1+\frac{\gamma+1}{2\gamma} \left(\frac{p_{\text{post}}}{p_{\text{pre}}} -1 \right)}.
\end{align*}
The QoI is an instantaneous mass flow measured by a sensor positioned at a distance $L$ from the initial position of the piston. The mass flow is discontinuous in time, depending on the position of the shock relative to the sensor.
\begin{align*}
    m(t) = \left\{
    \begin{array}{ll}
        \rho_{\text{pre}} u_{\text{pre}} \  \text{ if }\   t<\frac{L}{u_{\text{shock}}},\\
        \rho_{\text{post}} u_{\text{post}} \  \text{ if }\  t>\frac{L}{u_{\text{shock}}}.
    \end{array}
\right.
\end{align*}

\subsubsection{Two uncertain parameters}
\begin{figure}
    \centering
    \includegraphics[width=1.0\textwidth]{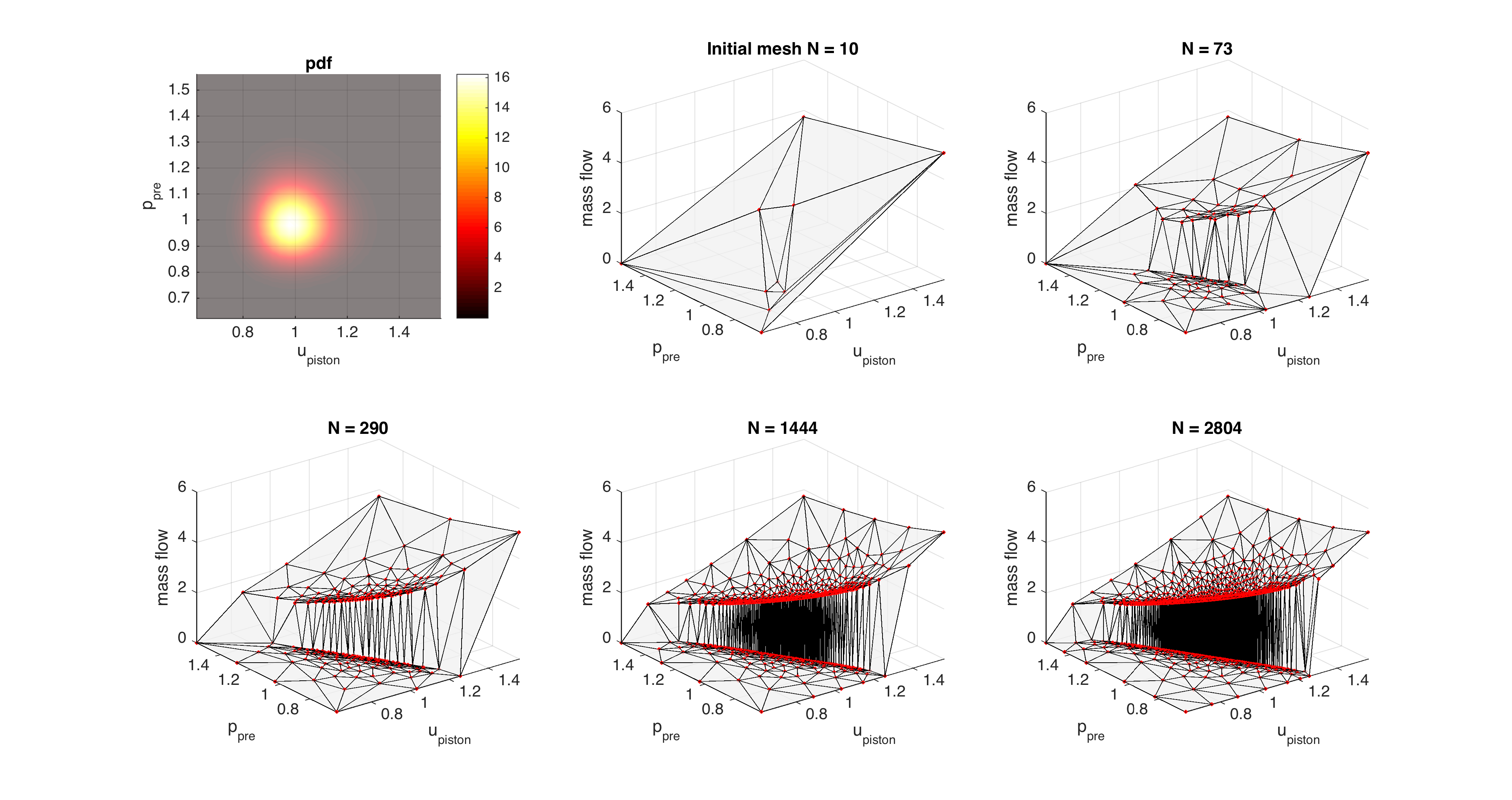}
    \caption{Piston problem: adaptive refinement of the response surface. Isocontours of the joint probability density function $\rho_{u_{\text{piston}},p_{\text{pre}}}$ is shown on the top left plot.}
    \label{piston_resp}
\end{figure}

\begin{figure}
    \centering
    \includegraphics[width=0.9\textwidth]{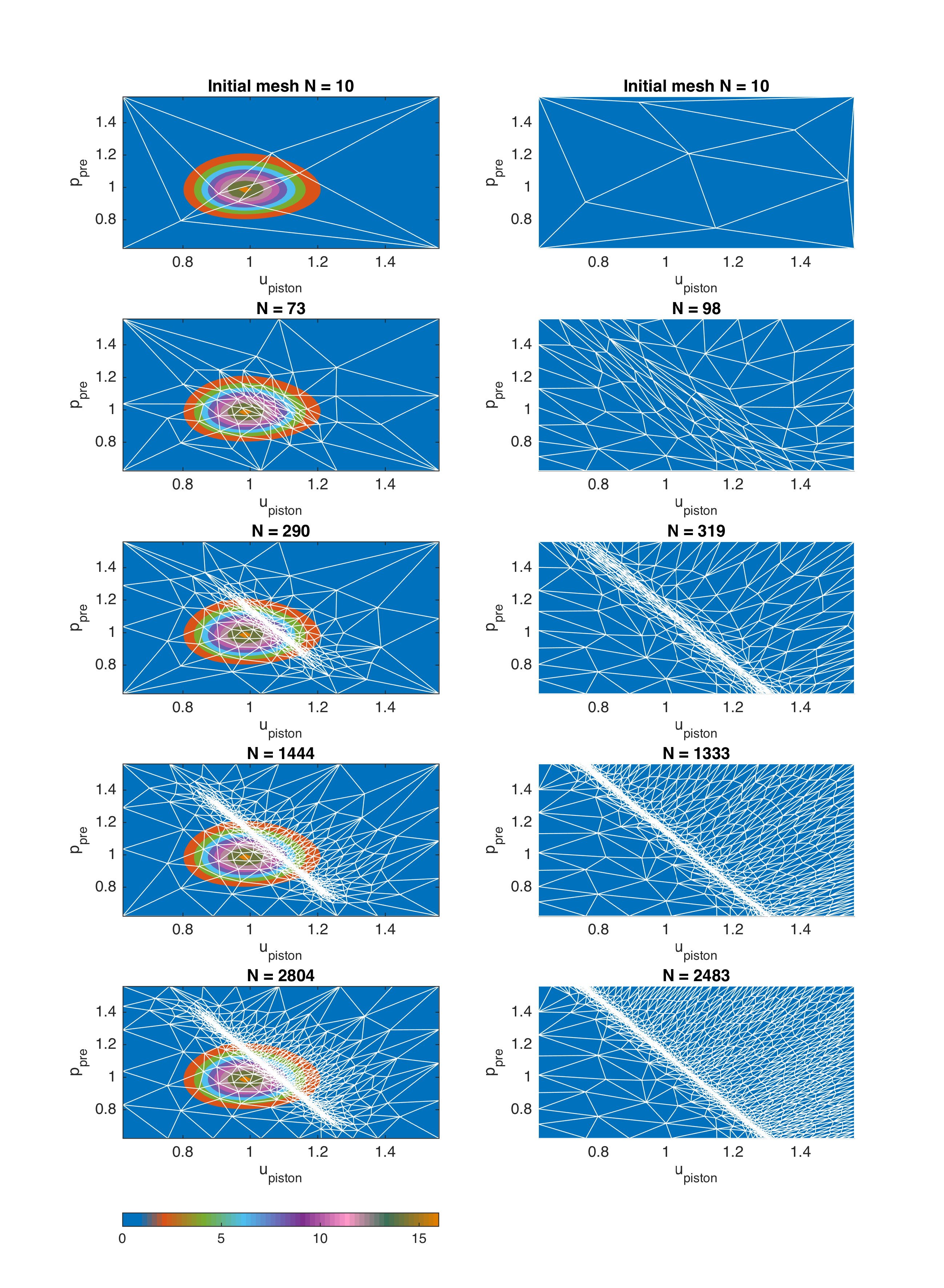}
    \caption{Piston problem: effect of the probability measures (lognormal pdf (left) vs. uniform pdf (right)). Sequences of stochastic mass flow mesh adaptations. 
    Discontinuous and nonlinear (upper) responses are only refined over regions of higher probability density values. Initial LHS DoEs follow underlying pdf.
    }
    \label{piston_lognorm_vs_unif}
\end{figure}

Identically to \cite{witteveen2009adaptive}, we first consider two uncertain parameters: $u_{\text{piston}}$ and $p_{\text{pre}}$, both being modeled as random variables following a lognormal distribution with mean $\mu_{u_\text{piston}} = \mu_{p_{\text{pre}}}=1$ and coefficient of variation (CV) such that $\text{CV}\equiv \sigma/\mu = 10\%$ (the equations are nondimensionalized). The QoI is the mass flow $m(t=0.5)$ at the sensor which is located at $L=1$.
For our application, depending on the random value of $u_{\text{shock}}$, $m(t=0.5)$ may sometimes be equal to $\rho_{\text{pre}}\,u_{\text{pre}}=0$ (because $u_{\text{pre}}=0$), while it may be $m(t=0.5)>0$ for other realizations. \\
The mesh adaptation process with the final response surface are shown in Figure \ref{piston_resp}. One can see how the discontinuous response with a complex curved front is gradually captured whilst taking into account the probability distribution $\rho_{{\bxi}}$.
Mesh adaptations of Figure \ref{piston_lognorm_vs_unif} complete the results by showing the impact of the probability measures on the localization of the discretization refinements. 


\begin{figure}
    \centering
    \begin{subfigure}{.49\linewidth}
        \centering
        \includegraphics[width=1\linewidth]{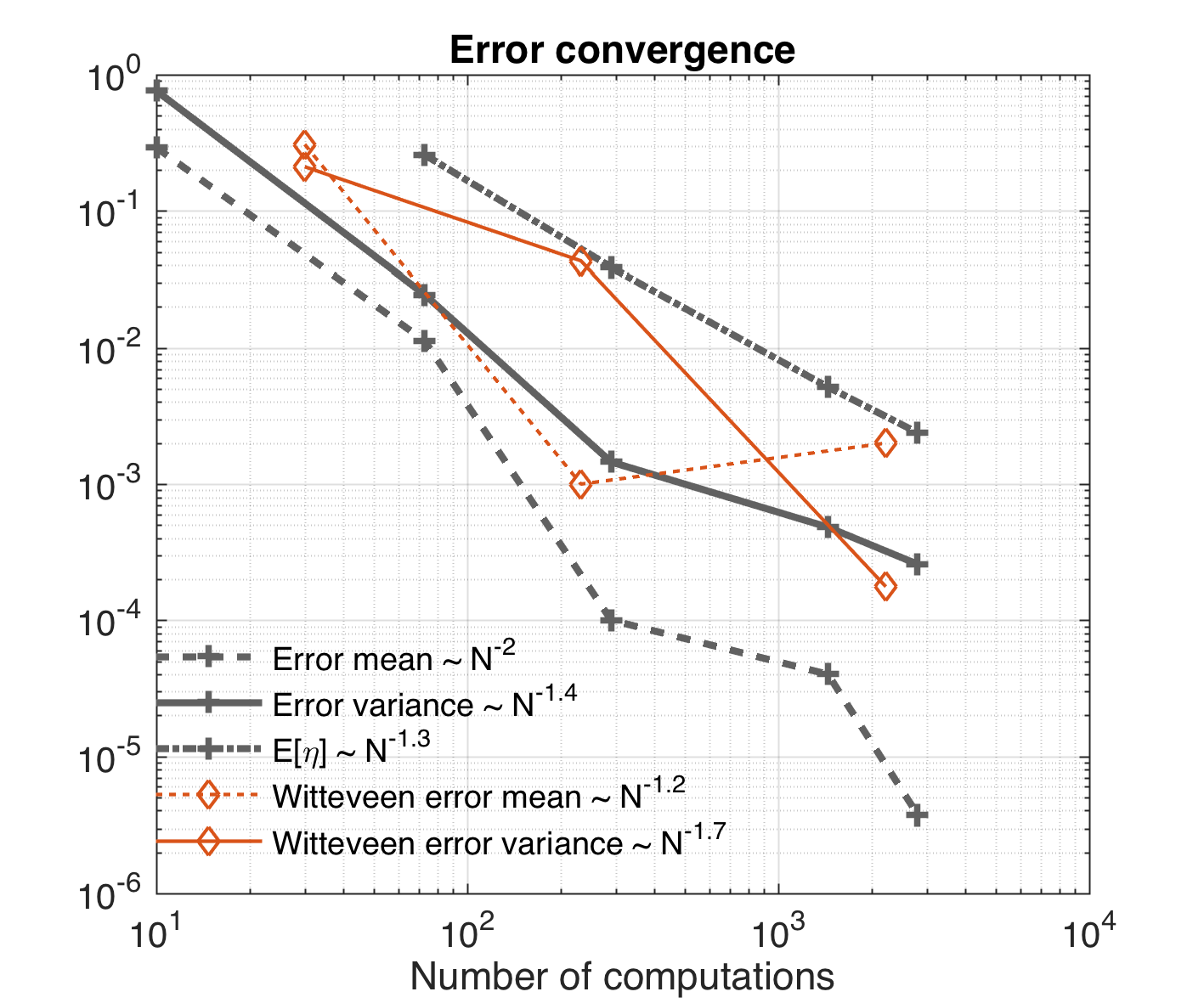}
        \caption{Piston problem with $2$ uncertain parameters}
        \label{piston_convergence}
    \end{subfigure}
    \begin{subfigure}{.49\linewidth}
        \centering
        \includegraphics[width=1\linewidth]{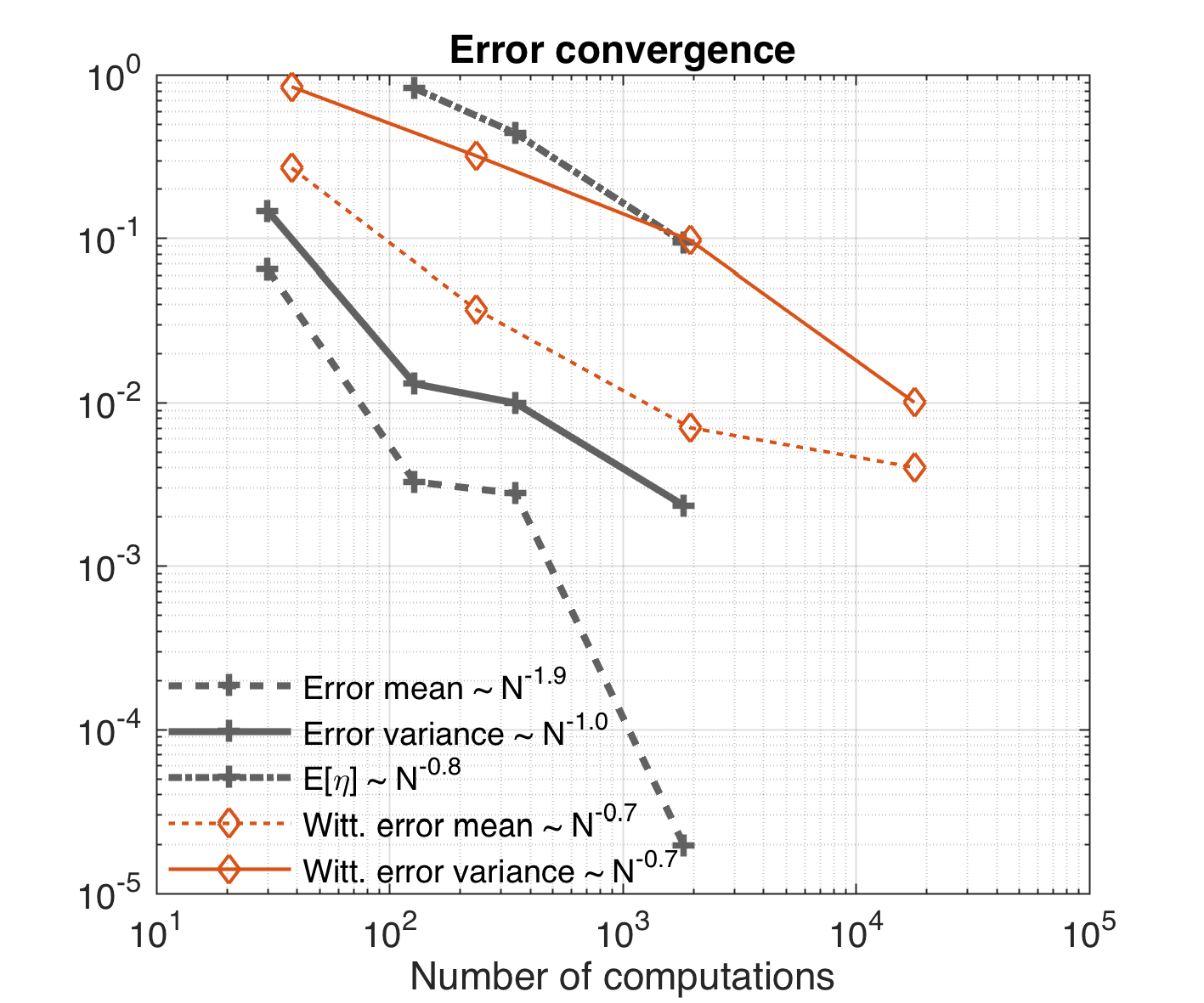}
        \caption{Piston problem with 3 uncertain parameters.}
        \label{piston_3d_convergence}
    \end{subfigure}
    \caption{Piston problem: convergence error of the mean, variance and the expectation of the estimated interpolation error 
    in $L^1-$norm. The results obtained by \cite{witteveen2009adaptive} are included for comparison.}
    \label{fig:piston_conv}
\end{figure}

Mass flow low-order moments (mean and variance) convergence results are given in Figure \ref{piston_convergence}. 
The errors in mean and variance are computed w.r.t.  mean and variance obtained from Monte Carlo simulation (obtained with $10^7$ samples). In the same figure the results obtained by Witteveen et al. \cite{witteveen2009adaptive} are also included which allows for a quantitative comparison.

\subsubsection{Three uncertain parameters}
Similarly to \cite{witteveen2009adaptive} a third lognormal random variable, the sensor position $L$, is added to the previous test case. 
In Figure \ref{piston_3d_convergence}, the convergence of the expectation of the interpolation error, the error in the mean and variance are shown. Plotted in the same figure are the results obtained by Witteveen et al.
~\\
~\\
Before {testing} our approximation method on more complex stochastic systems -- such as compressible CFD problems for which a very accurate Monte-Carlo reference solution is out of reach -- it is worth summarizing our findings for all of the examples treated previously:
\begin{itemize}
    \item we show a {\em convergence in the mean} of the proposed metric-based approximation method for multivariate nonlinear {\em discontinuous} functionals. More specifically, the method is at least $2^{\text{nd}}-$order in the $L^1-$norm, i.e. the mean error is such that: $\mathbb{E}\left [ \eta \right ]  \sim \Ocal \left (\beta N^{-2\kappa /d_{\bxi}} \right )$, with $(\beta>0,\kappa \gtrsim 1)$ positive constants
    \item the convergence is {\em robust}: -- for a given  functional with a fixed number of random parameters, a change in the nature and/or the regularity of the probability measure of the parameters does not affect the convergence rate of this mean error, -- the choice of the refinement step size does not greatly affect the convergence rate of the error, i.e. $\kappa$
    \item {the choice of the refinement step size affects the accuracy: lowering the step size overall requires less samples for a given accuracy}
    \item we have verified that the $L^1-$error of the surrogate solution mean is always lower in magnitude than the mean error of the surrogate; while its convergence rate is significantly higher than the latter one
    \item we have noticed that the $L^1-$error of the surrogate solution variance is larger than the error of the mean solution with a rate of convergence lower than the latter one, but higher than the convergence rate of the mean error
    \end{itemize}
    For the piston problem, we have compared our results to the SSC results proposed by \cite{witteveen2009adaptive}, and we have noticed that:
\begin{itemize}
    \item the proposed metric-based approximation method generates convergence rates -- for the surrogate solution mean that are larger  than the ones obtained with the SSC method and -- for the surrogate solution variance that are of comparable order as the ones obtained with the SSC method
    \item convergence constants $\beta$ are smaller with our method
    \item our convergence is more monotone than the SSC convergence that is based on local refinements.
\end{itemize}

In the following, we move to more complex applications involving CFD models and simulations.

\subsection{CFD applications}
The numerical method proposed in this work has been tested for several CFD configurations, not presented here, such as internal and external Euler flows, e.g. in a scramjet engine and around a NACA0012 \cite{langenhove2017these}.
In the following, we consider the fully coupled case, for which we attempt to control both deterministic and stochastic numerical errors of a supersonic/hypersonic inlet problem. For those problems no exact solutions are available, therefore only estimated errors provided by the method will be reported.
\subsubsection{Supersonic/Hypersonic inlet problem}

\begin{figure}[!h]
\centering
\includegraphics[trim=0cm 1cm 0cm 1cm, clip=true, scale=1.1]{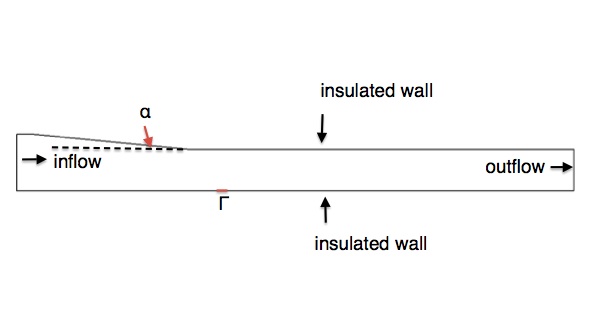}
\caption{Inlet problem: geometry configuration with an illustration of the targeted area $\Gamma$.}
\label{fig:inletGeom}
\end{figure}

Scramjet and ramjet propulsions are developing techniques of high technical and economical interest for high-speed engines. 
While scramjet engines are
fundamentally simple in concept, they are difficult in realization as current simulation
capability of in-flight performance is overwhelmed
by numerous uncertainties (e.g. natural variability of
flight scenario, effect of geometrical variability and manufacturing
tolerances, fuel conditions and combustion kinetics,...) and errors due to the multi-physics nature of
the problem.\\
In this type of dual-mode engines, the inlet (or isolator) plays an important role. For example, in ramjet mode, a strong precombustion shock forms in the inlet, resulting in a subsonic combustor entrance flow. More generally, the pressure typically rises in the system  through a series of shocks known as \textit{shock train}. The structure and shape of the shock train depends on the inlet entrance conditions: for instance, -- normal shock occurs for moderate Mach numbers of about 2-3, whereas for higher  Mach numbers -- oblique shocks occur. \\
In this section, we are interested by the internal flow for an inlet configuration described in Figure \ref{fig:inletGeom}. 
While this configuration with sharp angles induces a solution
with numerous shock waves, the position of the shock-train directly affects the combustion
and engine performance. 
Experimental investigations for this type of configuration can be found in \cite{wagner2007experimental,wagner2009experimental}. 
Another challenge of hypersonic flows is the simulation of the strong interaction
between those shocks and the turbulent boundary
layers. In this preliminary work, we do not address this issue and neglect viscosity and therefore we rely on compressible
Euler equations modeling, cf. Section \ref{subsec:Euler}.

\begin{figure}[!h]
\centering
\begin{tabular}{c}
(a)\includegraphics[scale=0.25]{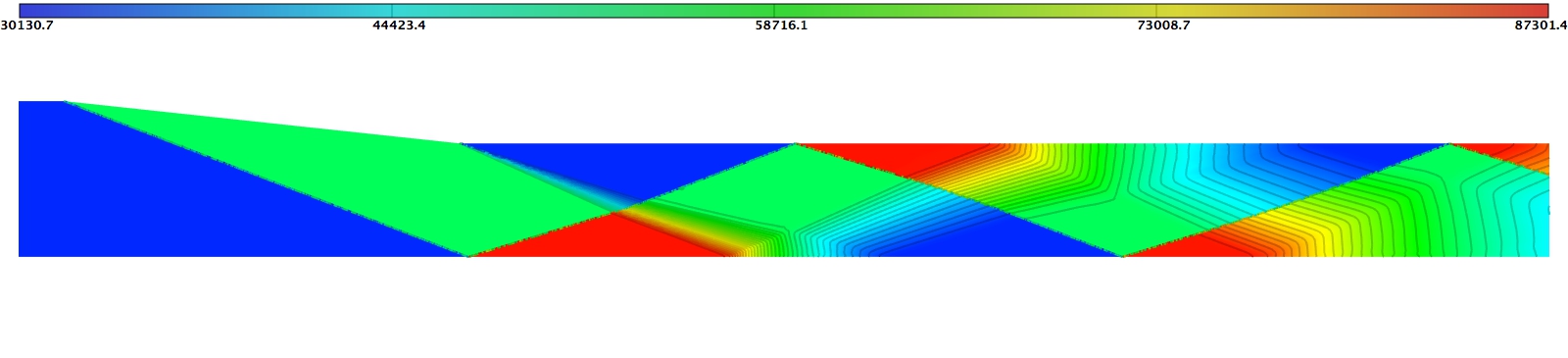} \\
(b)\includegraphics[scale=0.25]{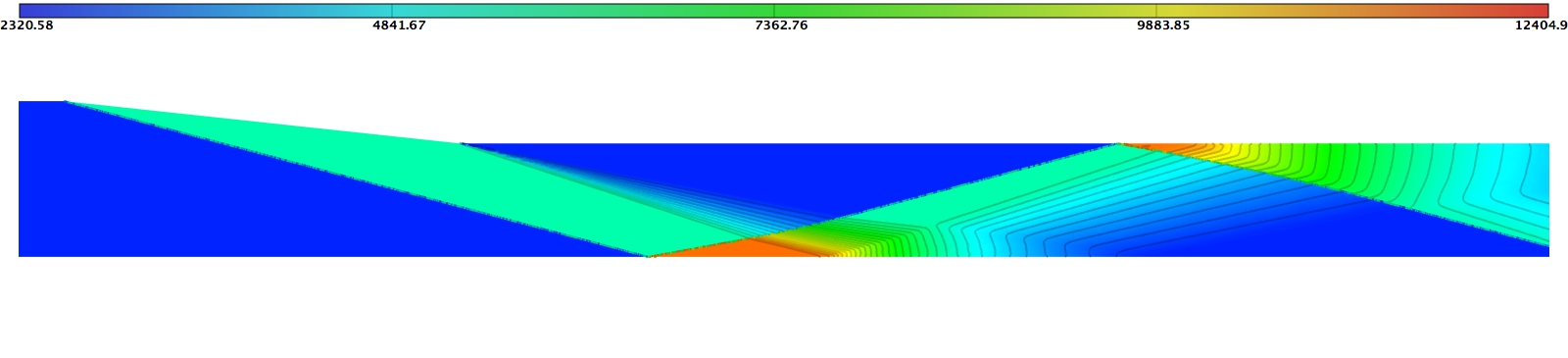} \\
(c)\includegraphics[scale=0.25]{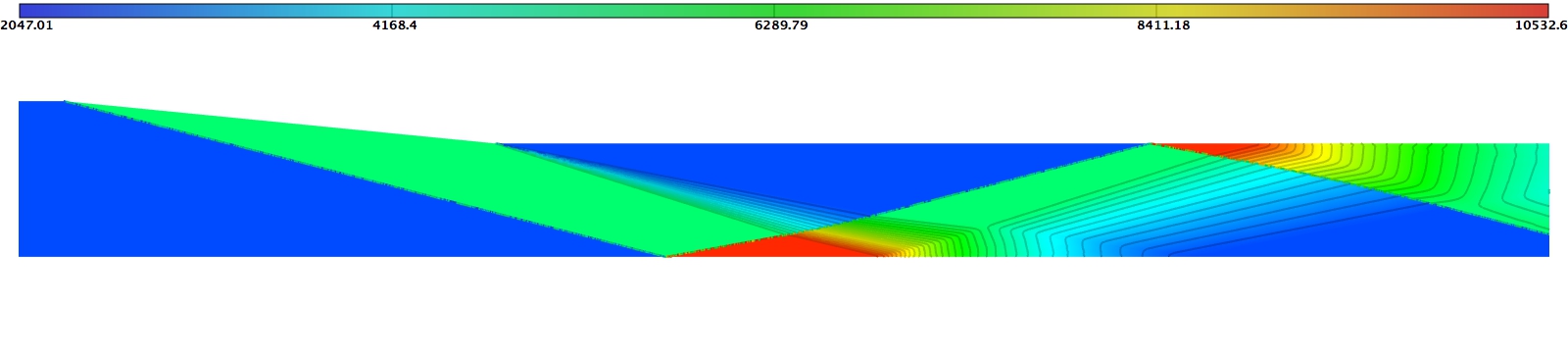} 
 \end{tabular}
\caption{Inlet problem: few pressure field isocontour realizations, for $M_\infty = 3.5$ (a) and $M_\infty=5.5$ (b-c) with ramp angles $\alpha = 6.1 ^\circ$ (a-b) and respectively $\alpha=5.6 ^\circ$ (c). Here, well resolved globally adapted meshes have been used in order to capture and show all flow features. }
\label{fig:inletPres}
\end{figure}

In our case, the QoI is related to the pressure signature on the lower surface wall, namely we define it as the integrated pressure coefficient $C_p$ over a short segment: $j=\int_{\Gamma}\frac{p-p_\infty}{\frac{1}{2}\rho_\infty || \boldsymbol{u}_{\infty} ||^2}\, d\bx$, where $\Gamma$ is a small region of size $5.08 mm$, located on the lower wall at $x=116.08mm$ from the entrance, slightly downstream on the end of the inlet ramp, cf. Figure \ref{fig:inletGeom}. 
It is thus interesting to analyse how the shocked flow and pressure distribution are impacted by operational and geometrical uncertainties: i.e. changes in the free stream Mach number $M_\infty$ and in the ramp angle $\alpha$. A change in the ramp angle will also affect its length, i.e. the boundary of the domain being modified. Palacios {\em et. al} \cite{palacios2012robust} analyzed this stochastic system with an identical configuration and obtained a discontinuous response surface in the QoI. The authors  were interested in developing an adaptive deterministic mesh associated to one nominal condition, in order to obtain a representative reference value of the QoI for the entire variability range.
Here, we have chosen uniformly distributed parameters with large variabilities: i.e. $\alpha \in \mathcal{U}_{[5.6;6.1]}$ and $M_\infty \in \mathcal{U}_{[3.5;5.5]}$. This is a challenging numerical problem since singularities arise in both physical and parameter spaces. Indeed, we can see from the study of  some pressure fields associated to several conditions, as illustrated in Figure \ref{fig:inletPres}, how variations in the $M_\infty$ and (to some extent) in the ramp angle $\alpha$ affect the shock train and pressure values along the computational domain and more specifically in the targeted zone $\Gamma$. One may also guess that depending on the combination between the flow speed and the angle of attack, the small area under interest may experience some pressure discontinuities (or not), strongly impeding on the numerical errors of the approximations.
Changes in the flow Mach number fields are also analyzed in \cite{langenhove2017these}.

\begin{figure}[!h]
\centering
\begin{tabular}{c}
(a)\includegraphics[trim=0cm 2cm 0cm 2cm, clip=true, scale=0.25]{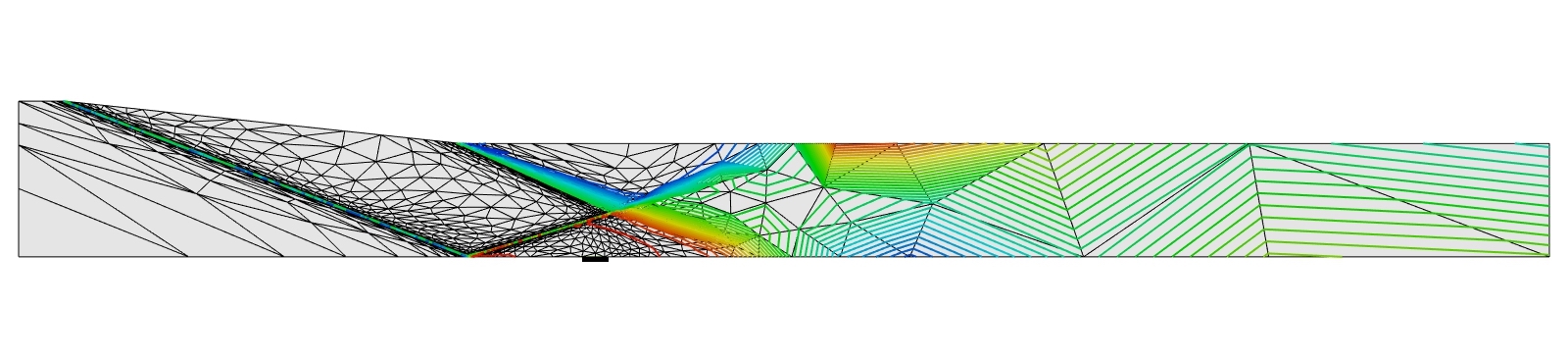} \\
(b)\includegraphics[trim=0cm 2cm 0cm 2cm, clip=true, scale=0.25]{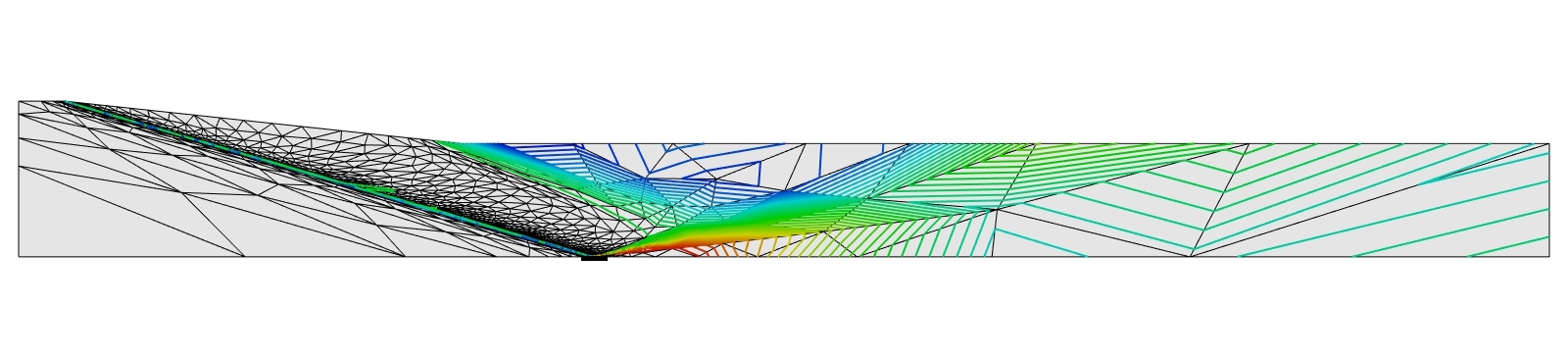} \\
(c)\includegraphics[trim=0cm 2cm 0cm 2cm, clip=true, scale=0.25]{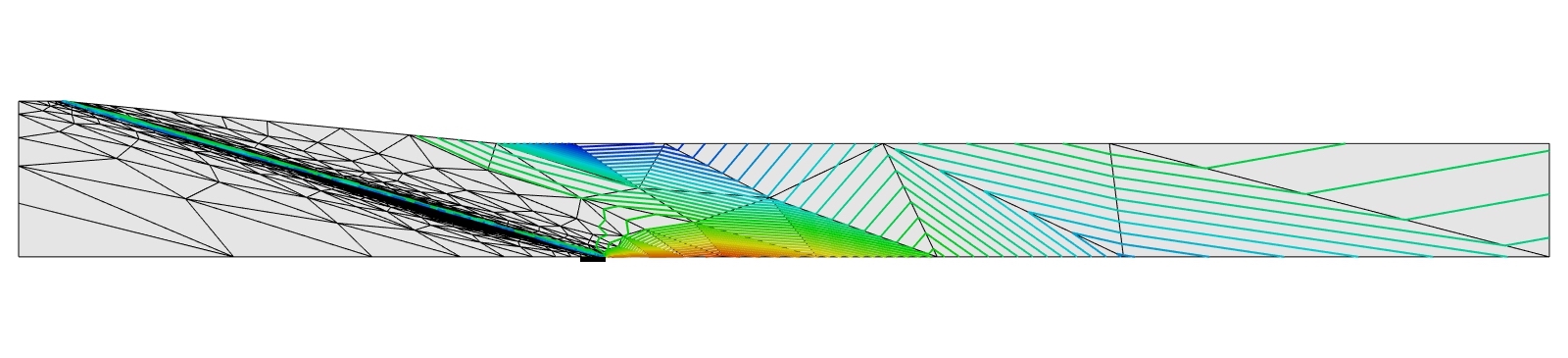}  
 \end{tabular}
\caption{Inlet problem: density iso-contours on rather coarse ($C_{\bx} =1000$) goal-based adapted meshes for $M_\infty = 3.5$ (a) and $M_\infty=5.5$ (b-c) with ramp angles $\alpha = 6.1 ^\circ$ (a-b) and respectively $\alpha=5.6 ^\circ$ (c). The location of $\Gamma$ is represented by a black thick line segment.}
\label{fig:inletGBmeshAll}
\end{figure}

\begin{figure}[!h]
\centering
\begin{tabular}{c}
(a)\includegraphics[trim=9cm 0cm 0cm 0cm, clip=true, scale=0.25]{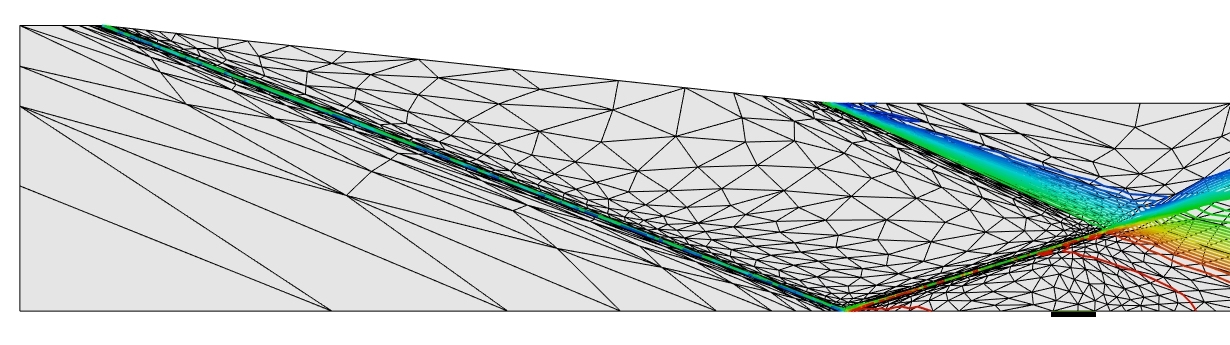} \\
(b)\includegraphics[trim=9cm 0cm 0cm 0cm, clip=true,scale=0.25]{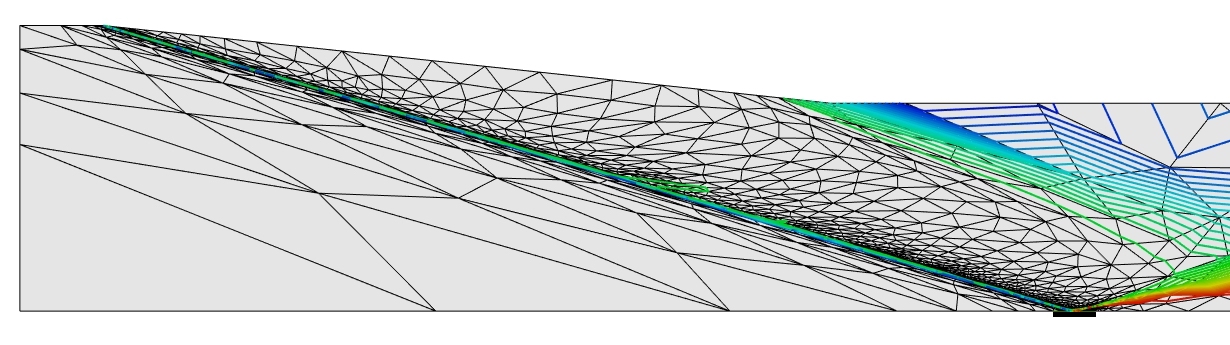} \\
(c)\includegraphics[trim=9cm 0cm 0cm 0cm, clip=true,scale=0.25]{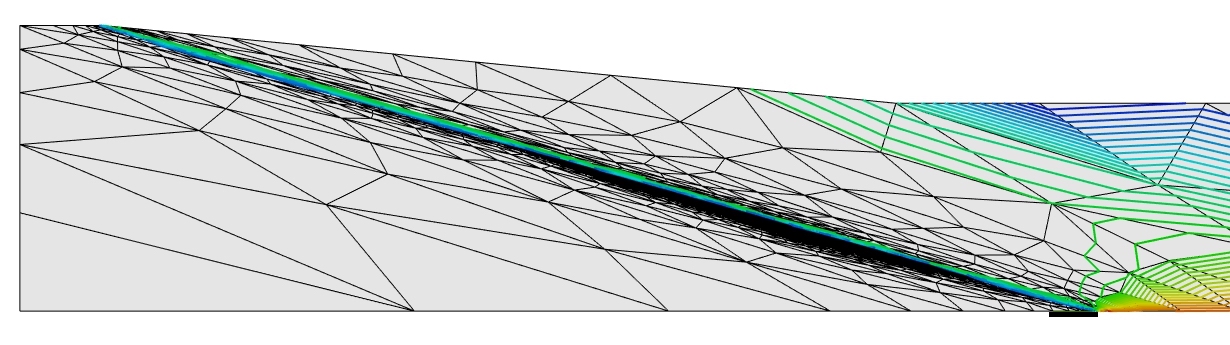}  
 \end{tabular}
\caption{Inlet problem:  close-ups of figures \ref{fig:inletGBmeshAll} (same caption).}
\label{fig:inletGBmeshZoom}
\end{figure}

We first need to make sure we control the discretization error for a given flow speed and geometry. While it is obvious that a good ``shock-capturing" method is needed for this problem, depending on the available computing ressources, it is in general not recommended to refine all shocks present in the domain. Our proposed adaptive method based on optimal control of both stochastic and deterministic errors is a sound approach to make the right selection for refinements. Thanks to the adjoint-based method, the mesh is efficiently adapted only in the regions with large impact for our QoI. 
Representative examples of goal-based adapted meshes for various ramp angles and Mach numbers are displayed in Figures \ref{fig:inletGBmeshAll} and \ref{fig:inletGBmeshZoom}.
With a closer look at the meshes, cf. Figure \ref{fig:inletGBmeshZoom}, we observe that the remeshing effort is solely focused in the regions impacting the QoI. When the first compression shock emitted at the beginning of the ramp hits the target domain or one of its boundaries, cf. (b) and (c), spatial mesh adaptation is needed to better capture this shock and its reflection. However, when the reflection is located upstream of $\Gamma$, cf. (a), then the expansion fan region generated at the lower ramp angle also needs refinement, because of its interaction with the first shock reflected. Further downstream the target, no mesh refinements are particularly needed for all cases. \\

We have applied our optimal adaptive strategy described in Section \ref{optStrategy} to control \textit{both} stochastic and deterministic errors. We arbitrarily chose a total of $n_{cycle}=19$ adaptive cycles where, at each cycle, according to Algorithm \ref{opt_err_control_algo} we choose to perform either deterministic or stochastic adaptation to lower the numerical error. 
Nevertheless, we have here voluntarily pushed the iterations further than necessary in order to monitor some asymptotic convergence trends.\\
For a practical purpose, we set maximal complexity values $\Ccal^{max}_{\bx}$ and $C^{max}_{\bxi}$ in order to avoid a cost exceeding our computational ressources. The initial DoE is set to $N_{\bxi,0}=30$ with a repartition drawn according to LHS plan mapped to the underlying pdf, while the initial spatial discretization for all samples is chosen to be a uniform finite-element mesh with $N_{\bx,0}=4000$ vertices.\\
\begin{figure}[!h]
\centering
\begin{tabular}{c c c}
\hspace{-0.35in} \includegraphics[scale=0.3]{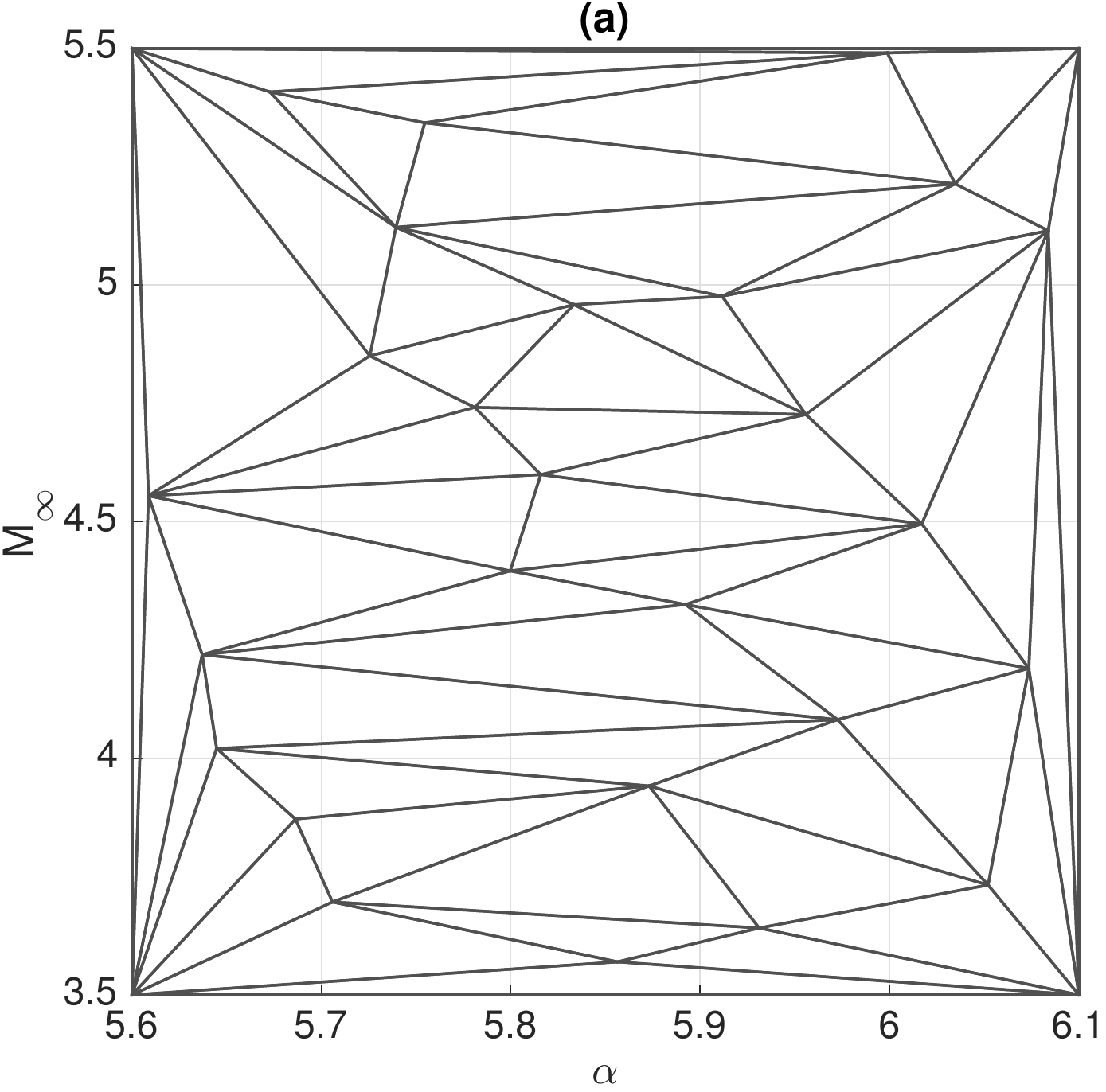} & \hspace{-0.35in} \includegraphics[scale=0.3]{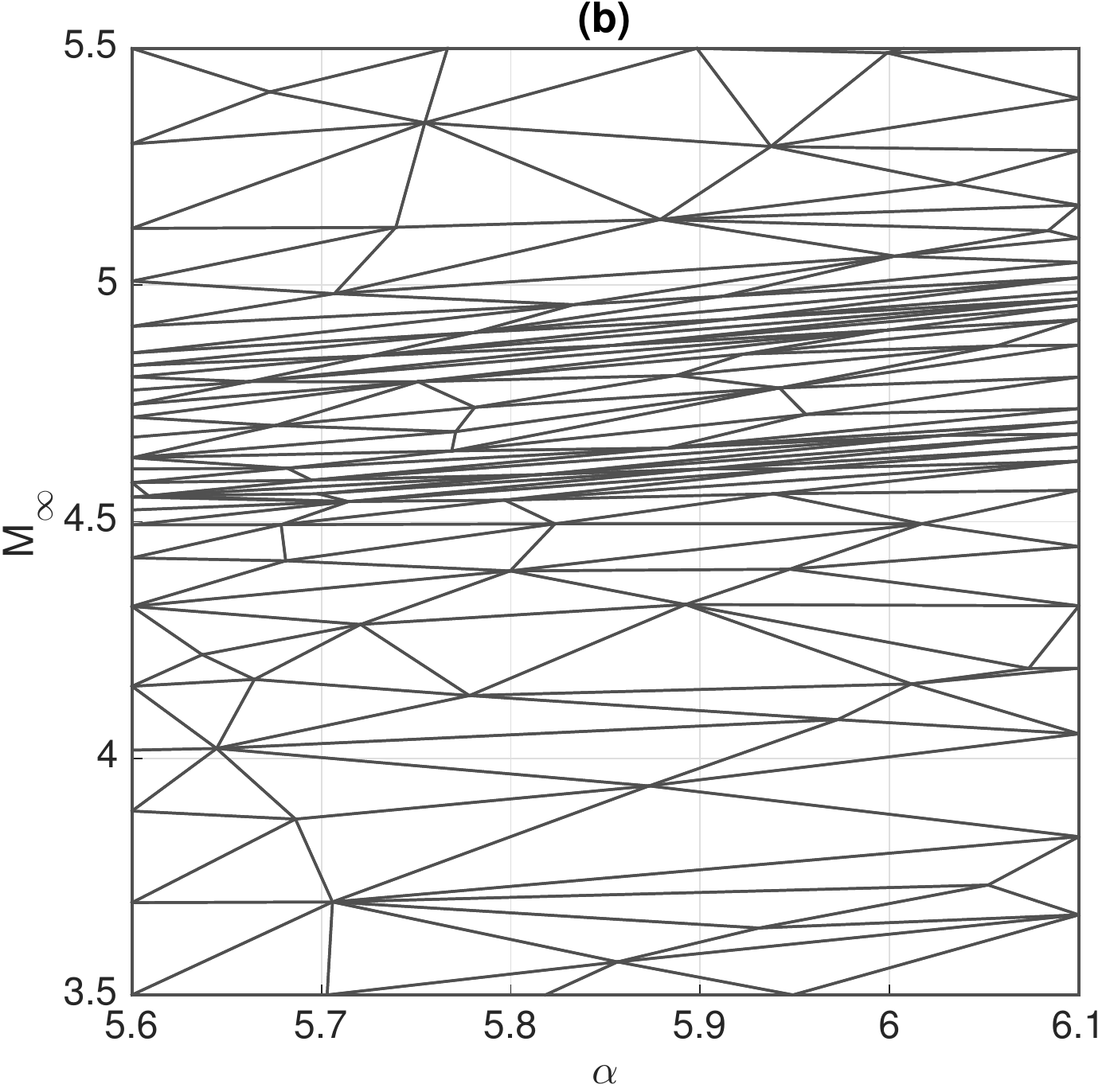} & \hspace{-0.35in} \includegraphics[scale=0.3]{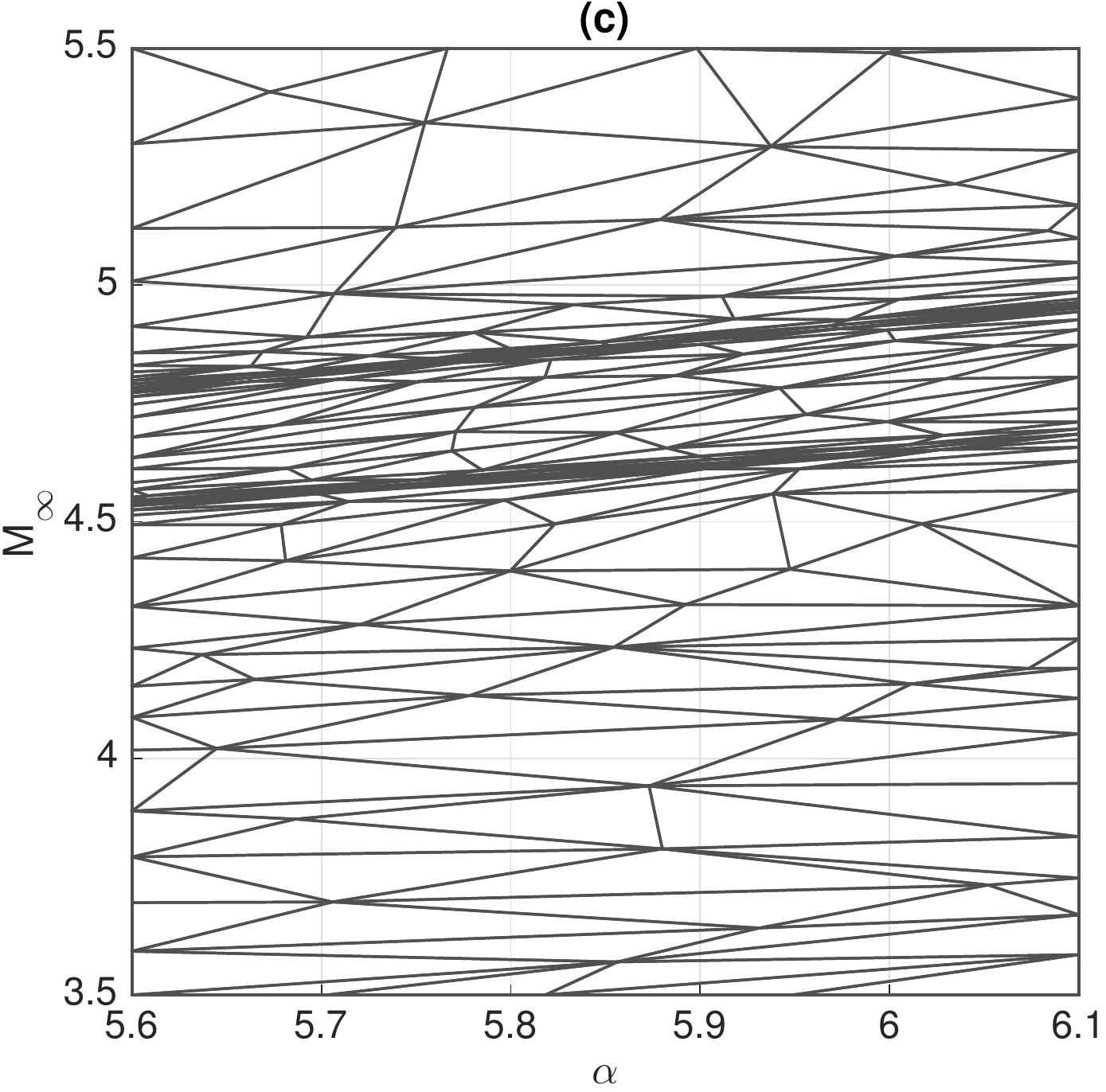}\\
\includegraphics[trim=2cm 7.5cm 2cm 7.5cm, clip=true, scale=0.3]{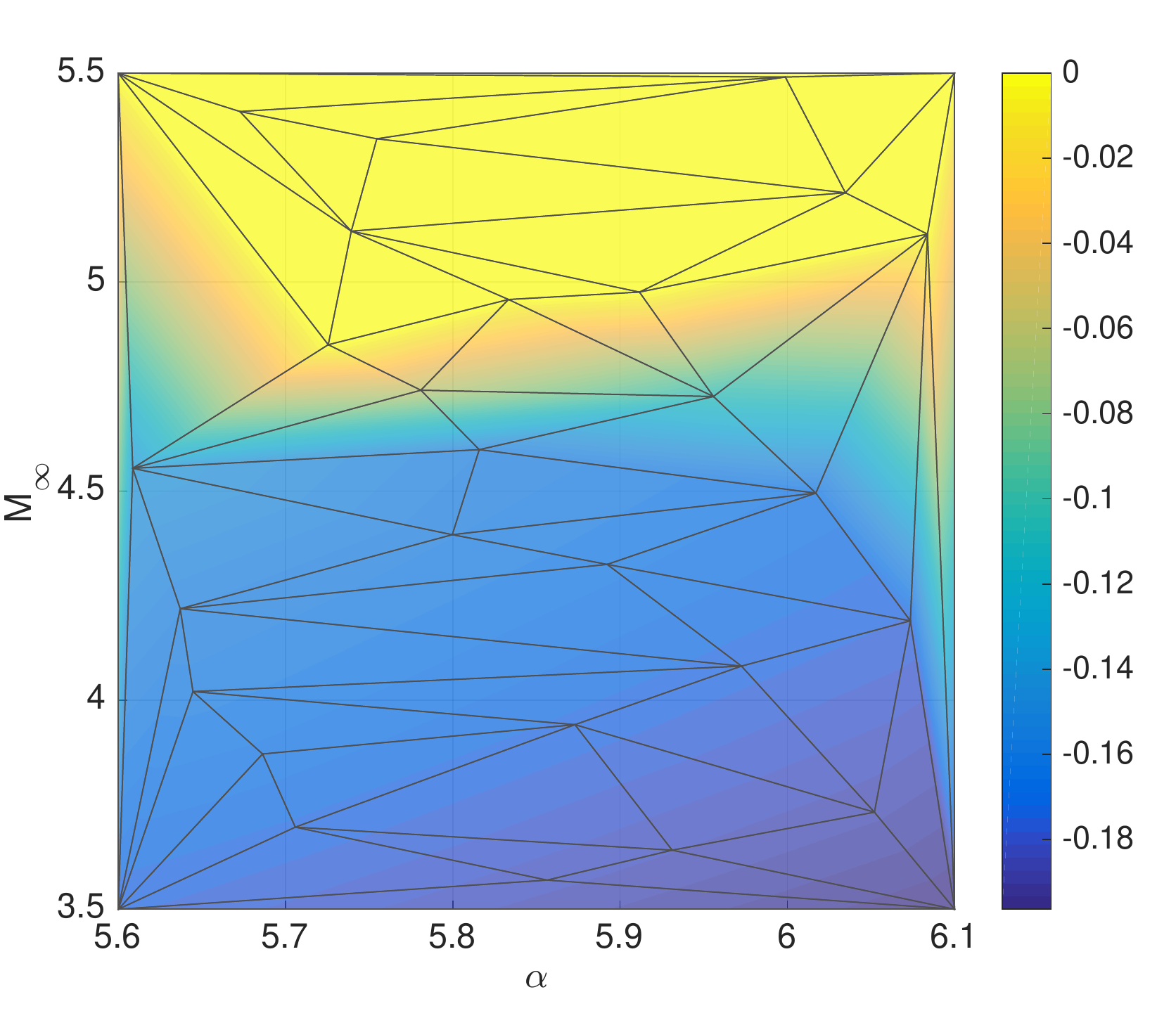} & \includegraphics[trim=2cm 7.5cm 2cm 7.5cm, clip=true, scale=0.3]{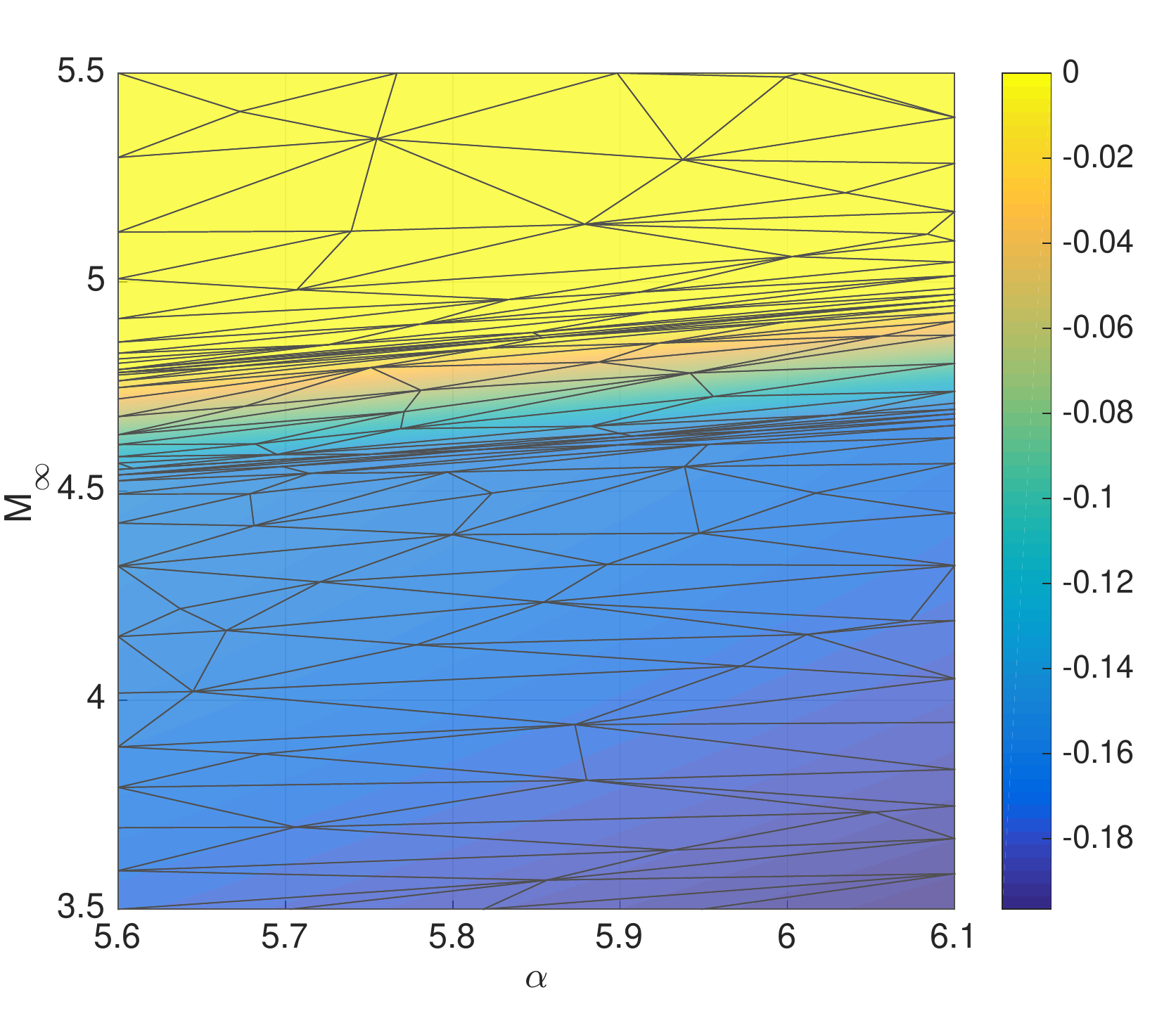} & \includegraphics[trim=2cm 7.5cm 2cm 7.5cm, clip=true, scale=0.3]{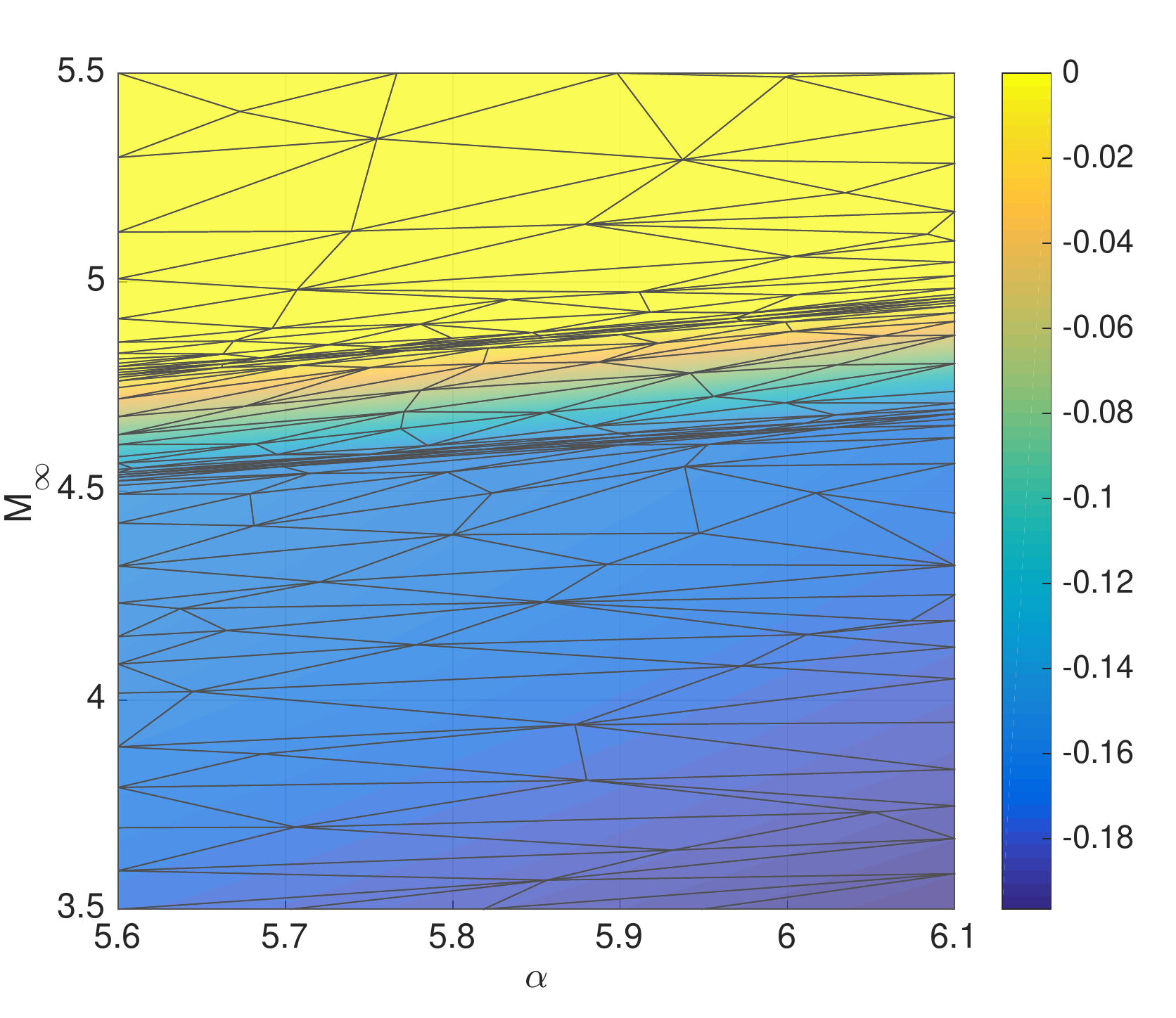}\\
\includegraphics[trim=2cm 7.5cm 2cm 7.5cm, clip=true, scale=0.3]{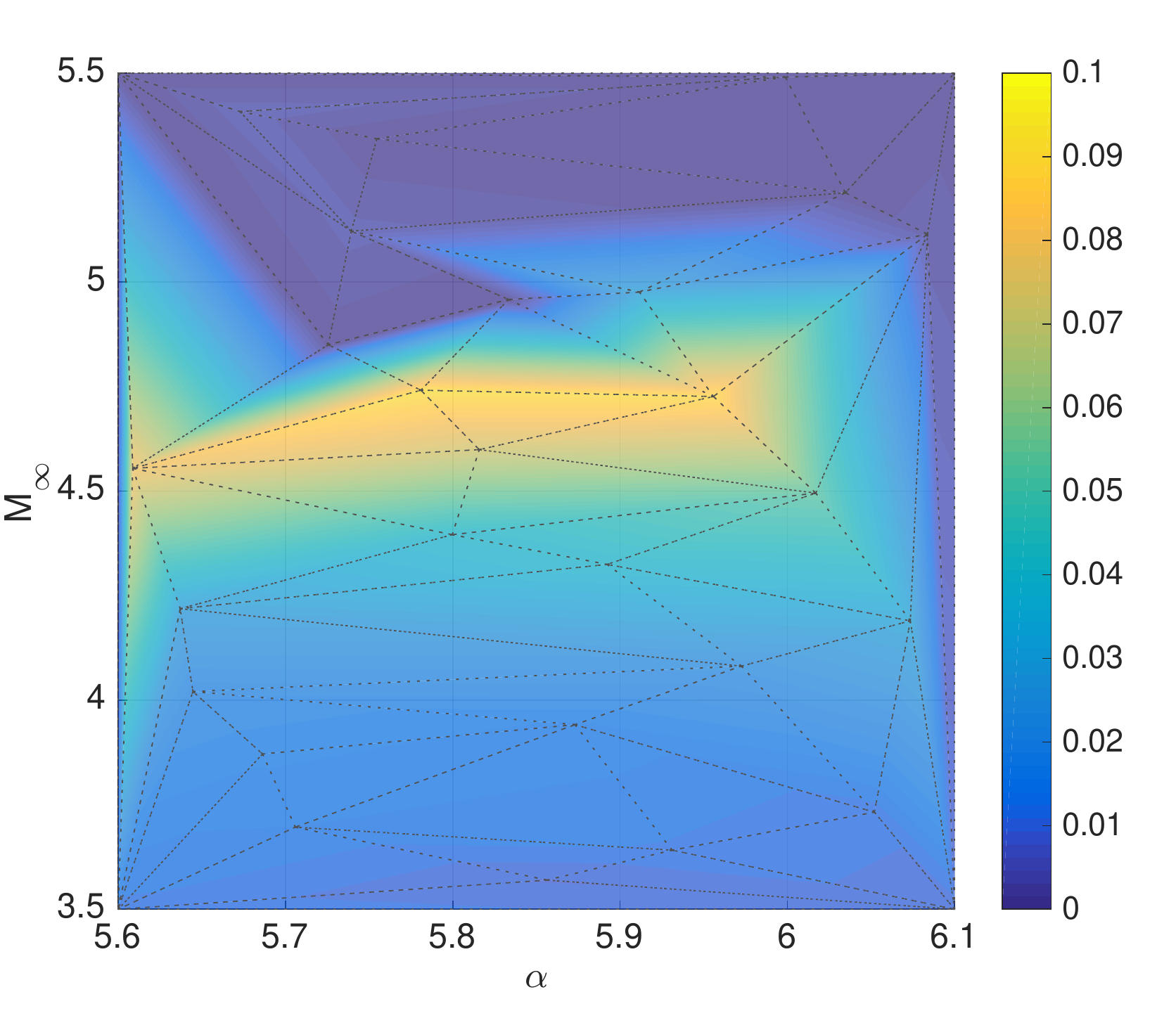} & \includegraphics[trim=2cm 7.5cm 2cm 7.5cm, clip=true, scale=0.3]{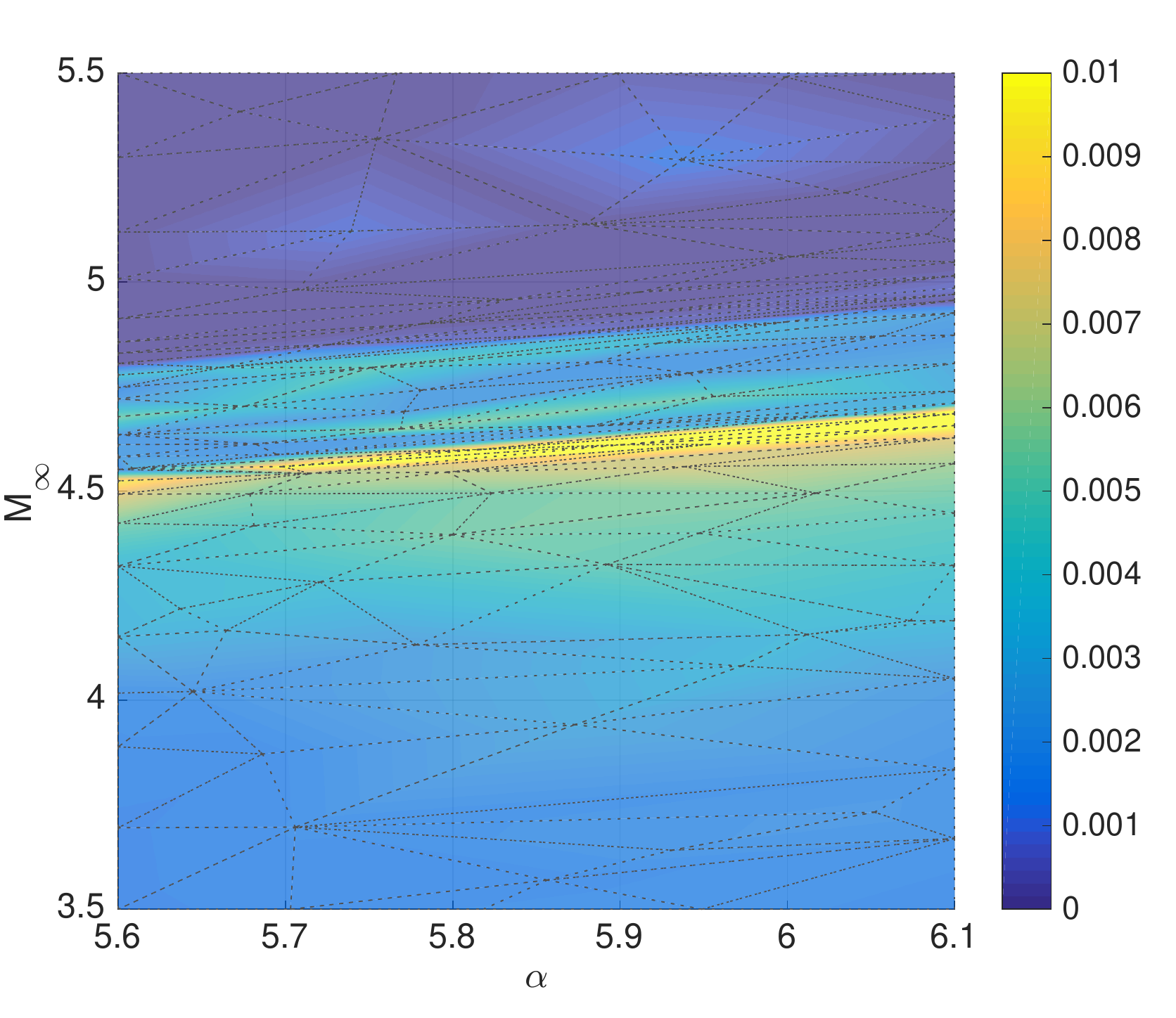} & \includegraphics[trim=2cm 7.5cm 2cm 7.5cm, clip=true, scale=0.3]{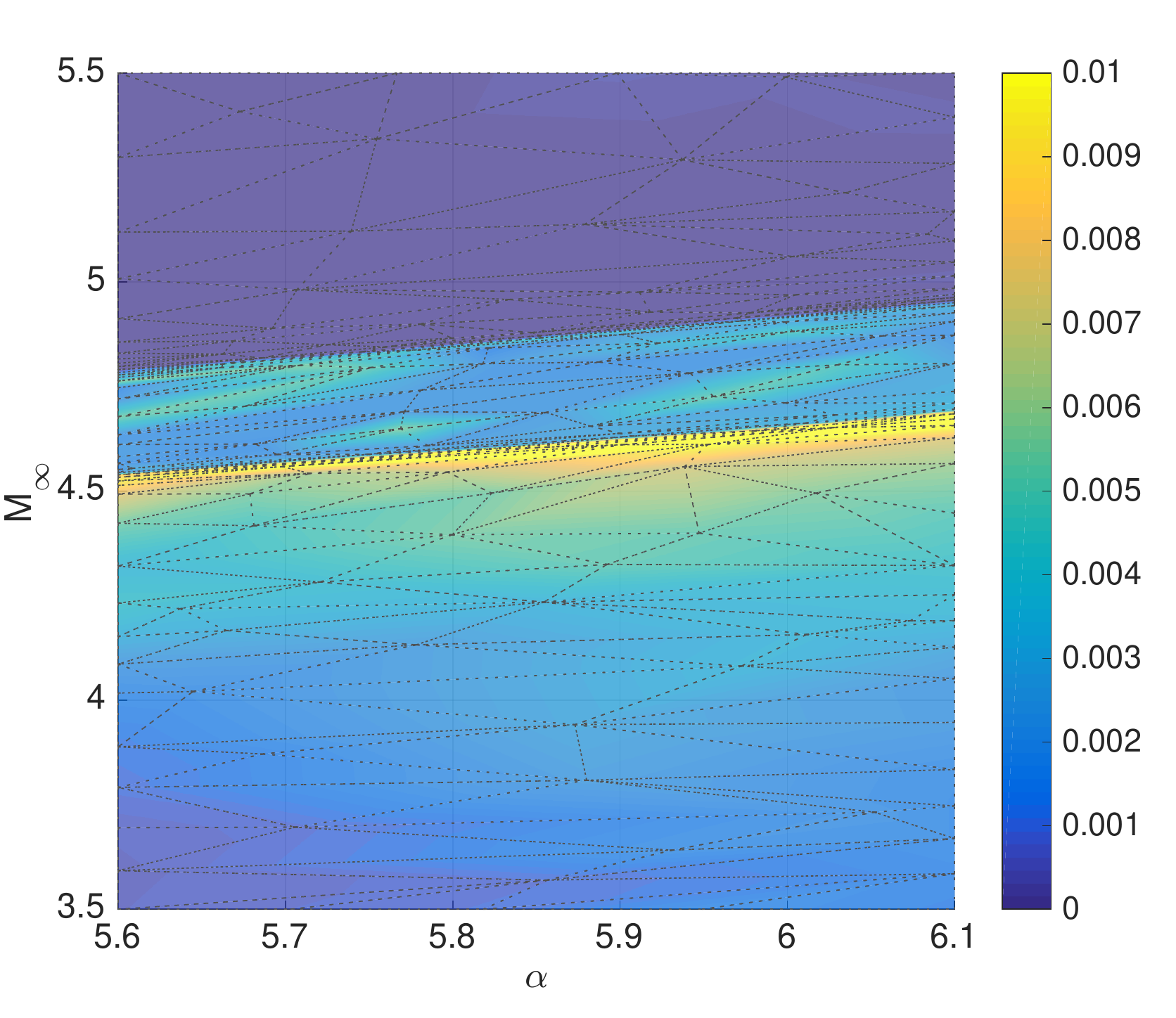}\\
\includegraphics[trim=2cm 7.5cm 2cm 7.5cm, clip=true, scale=0.3]{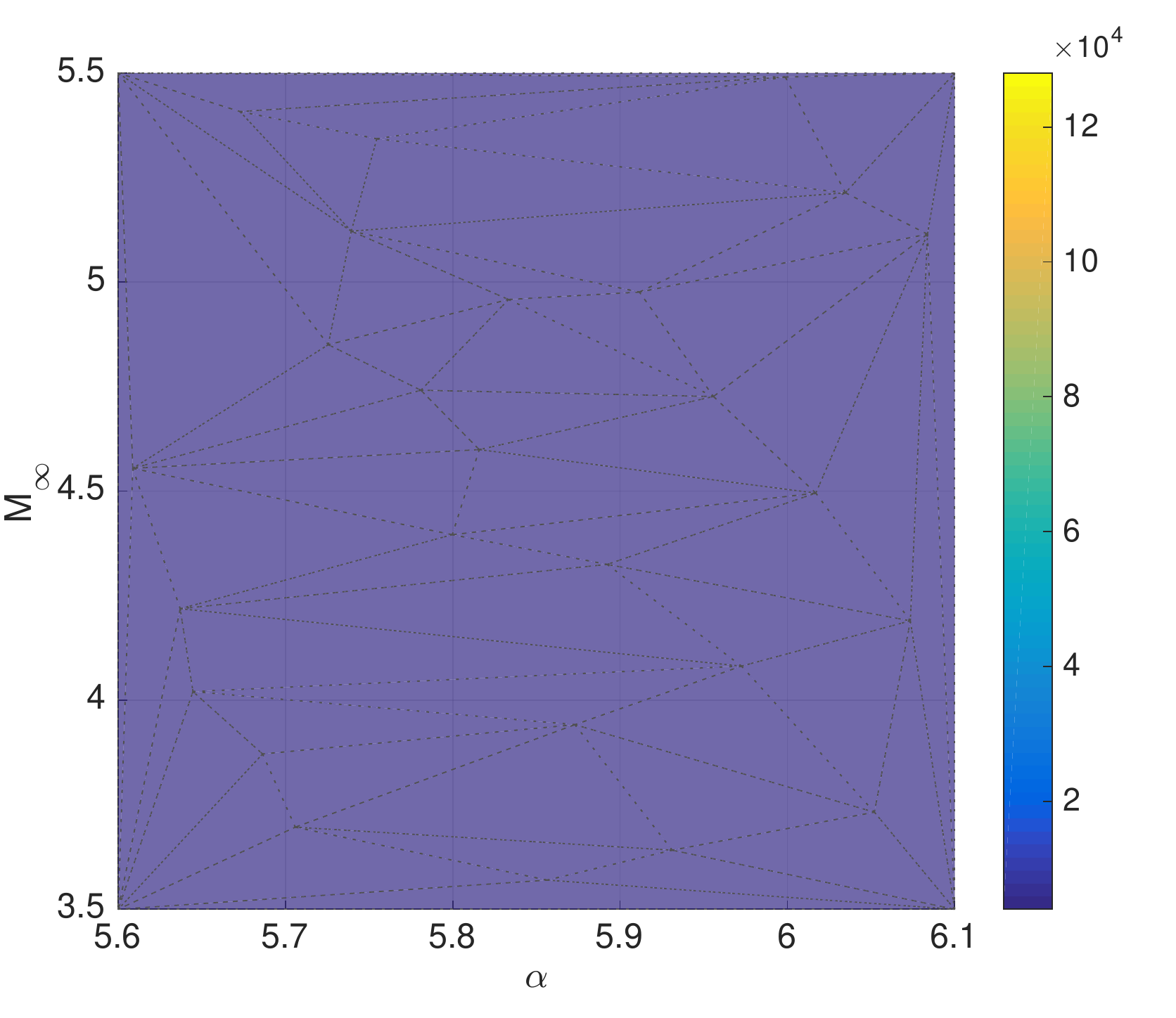} & \includegraphics[trim=2cm 7.5cm 2cm 7.5cm, clip=true, scale=0.3]{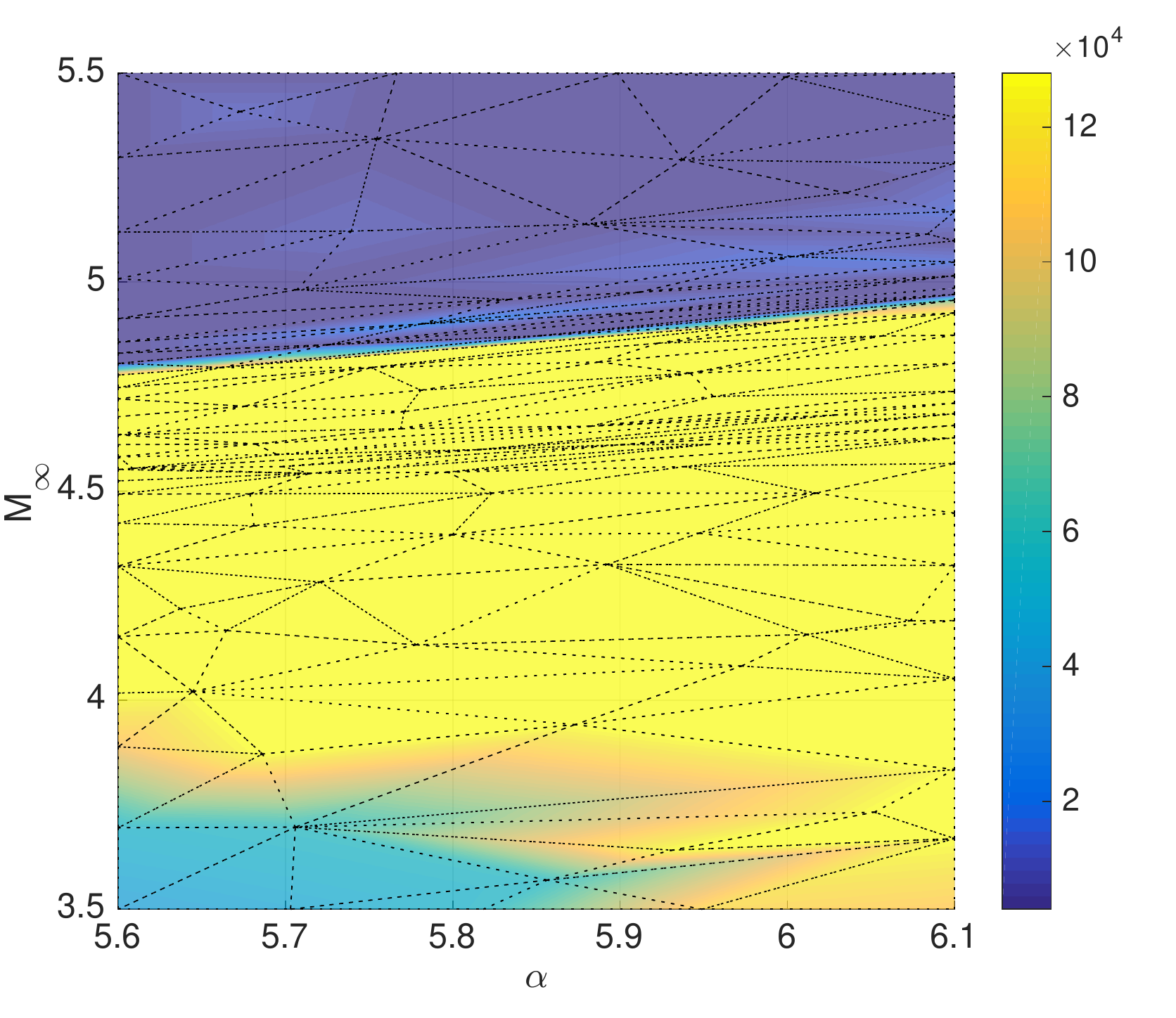} & \includegraphics[trim=2cm 7.5cm 2cm 7.5cm, clip=true, scale=0.3]{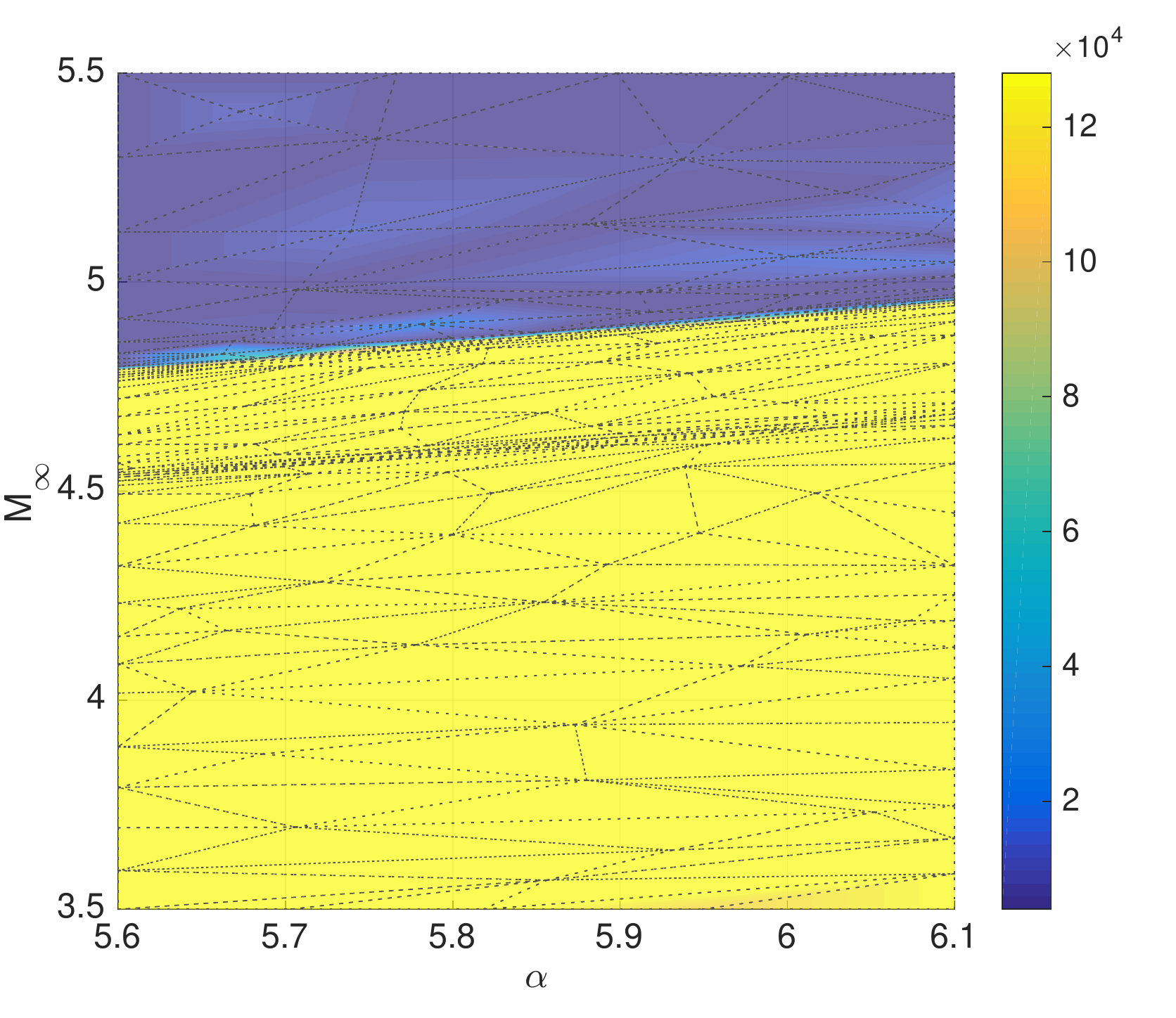}
\end{tabular}
\caption{Inlet problem: some results of the coupled adaptation strategy for error control. 
Results are displayed at the end of cycle $0$ (column (a)), cycle $2$ (b) and cycle $19$ (c).
Adapted stochastic partitions (first row) and isocontour maps are shown in the parametric space: the surface response of the approximated QoI $j(M_{\infty},\alpha)$ ($2^{\text{nd}}$ row), -- spatial discretization error maps  $\varepsilon(M_{\infty},\alpha)$ ($3^{\text{rd}}$ row) and -- corresponding optimal spatial complexity map $\Ccal_{\bx}(M_{\infty},\alpha)$ (bottom row).  }
\label{fig:inlet_MeshErrorComplexitySurface}
\end{figure}
Figure \ref{fig:inlet_MeshErrorComplexitySurface}
 shows some results at the final sub-iteration of adaptation cycle 0 (i.e. at the initial stage), cycle 2 and final cycle $n_{cycle}$.
The final design contains $N_{\bxi,n_{cycle}}=177$ simulations, averaging about $\overline{N_{\bx,n_{cycle}}}=90 000$ vertices.
We see that the integrated $C_p$ response (cf. second row) in the parametric space is quite regular except for two sharp oblique transition regions of varying $M_{\infty}$ and $\alpha$, corresponding to the impingement of the inlet first compression shock onto the region of interest $\Gamma$. The first of these regions (for lower $M_{\infty}$) corresponds to the entrance of the shock to the left of $\Gamma$, while the other one corresponds to its exit. The metric-based stochastic adaptation is able to sense these regions of poor smoothness. It gradually and anisotropically refines the approximation by adding more samples in these areas. We emphasize that the singularities are captured early in the adaptive process, starting from cycle $1$ and that the average stochastic error is quickly minimized. On the other hand, largest deterministic error contributions arise in those regions as well. They correspond to situations where the area of interest is impacted by the presence of the shock wave and require mesh refinements. Therefore, target spatial complexities are also increased for the concerned samples. However, quite early in the adaptation cycles, the deterministic error magnitude $\varepsilon$ in the lower part of the domain eventually requires complexities that exceed the chosen $\Ccal^{max}_{\bx}$ (as seen in the center and most right bottom images).
We also notice that some of the CFD computations who have reached the minimum error level (or maximum computational capacity) are kept unchanged, the data being only read through the remaining iterations until the stopping criterion is reached. This is the case for many of the computations with high Mach values. At the final stage of our example, we conclude from the results that our computational budget has been used to even down the deterministic error contribution over the parametric domain. Largest errors subsist for parameters combinations for which the compression shock sits at the left boundary of $\Gamma$. This is indeed the case where the reflected compression shock and expansion fan generated along the ramp have the strongest pressure interaction, requiring subsequent mesh adaptation effort. More detailed quantitative results are summarized in Table \ref{tab:inlet_stats}.

\begin{figure}
    \centering
    \begin{subfigure}{.495\linewidth}
        \centering
        \includegraphics[trim=2.5cm 6.5cm 2.5cm 7cm, clip=true, scale=0.45]{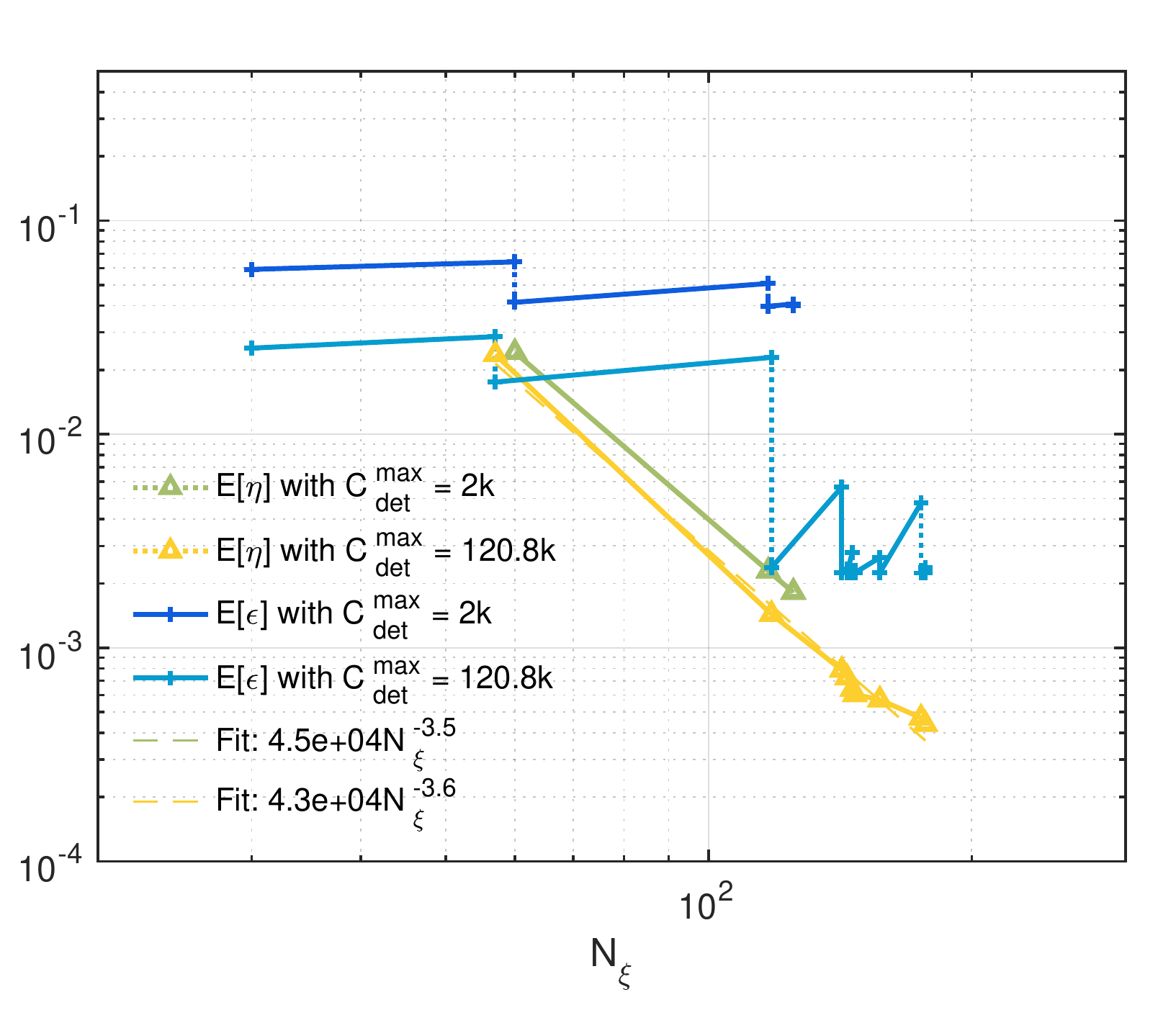}
        \caption{Fully (deterministic/stochastic) adapted schemes for different spatial refinements complexity targets $\mathcal{C}_{x}^{max}$.\\
        Results are displayed for $\mathcal{C}_{x}^{max}=120.8e3$ vs. $\Ccal^{max}_{\bx}=2e3$ complexities.}
        \label{subfig1:inlet_errorConv}
    \end{subfigure}
    \begin{subfigure}{.495\linewidth}
        \centering
        \includegraphics[trim=2.5cm 6.5cm 2.5cm 7cm, clip=true, scale=0.45]{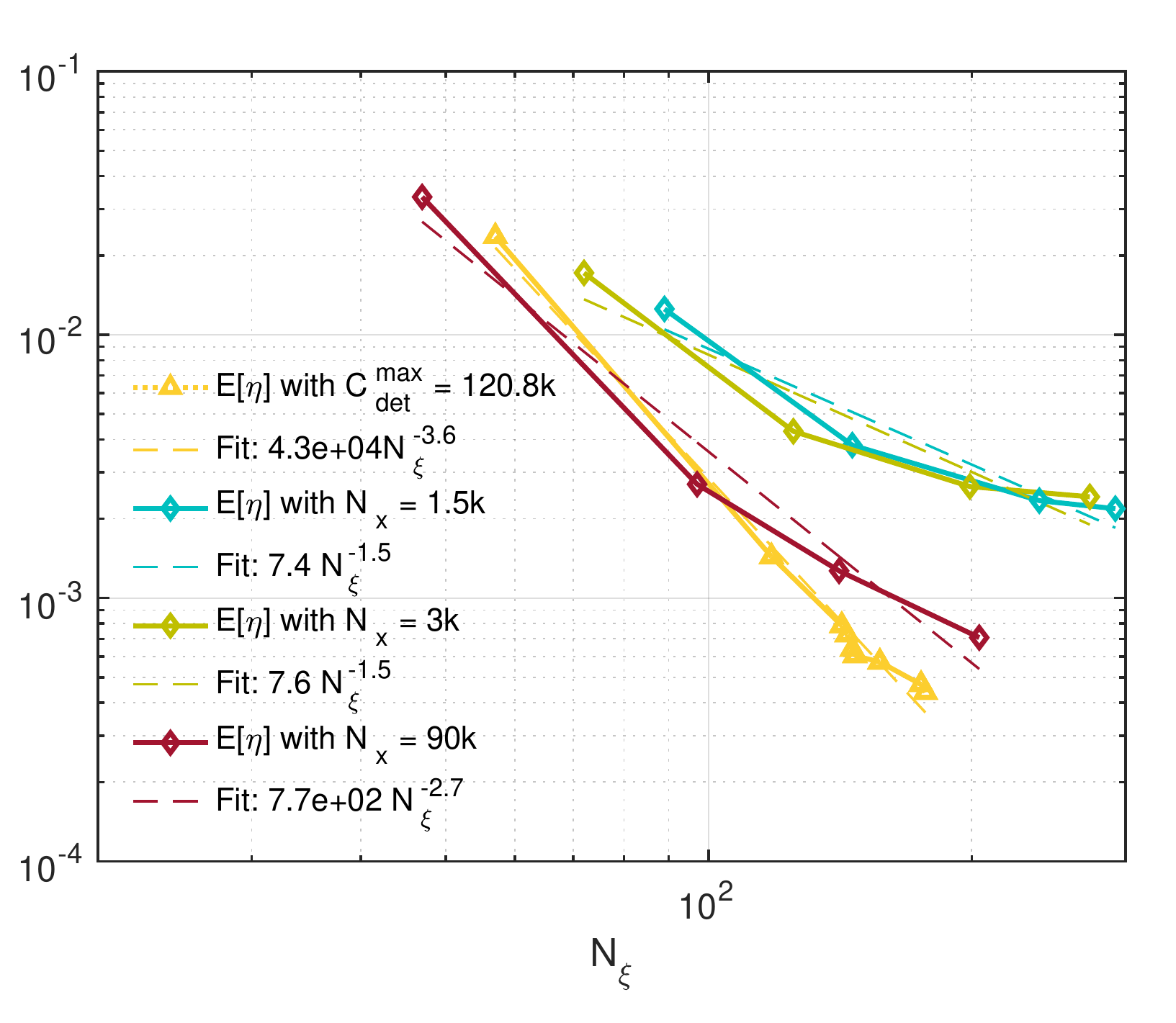}
        \caption{Fully (deterministic/stochastic) vs. partially (stochastic) adapted schemes. Stochastically adapted schemes are investigated for various isotropic spatial discretizations: $N_{\bx}=\{1500,3000,90000\}$.}
        \label{subfig2:inlet_errorConv}
    \end{subfigure}
    \caption{Inlet problem: convergence of the mean approximation errors, i.e. deterministic $\BigE [\varepsilon]$ and stochastic $\BigE [\eta]$ contributions (a) and stochastic $\BigE [\eta]$ contributions (b) respectively, vs. the total number of solution samples $N_{\bxi}$.}
\label{fig:inlet_errorConv}
\end{figure}

~\\
The evolution of the convergence of average deterministic and stochastic components of the error with increasing resolution in both parametric and physical domain are illustrated in Figure \ref{fig:inlet_errorConv}.
Figure \ref{subfig1:inlet_errorConv} compares fully adapted (deterministic/stochastic) schemes for different spatial refinements complexity targets $\mathcal{C}_{x}^{max}$, while Figure \ref{subfig2:inlet_errorConv} compares full with partially adapted (stochastic) schemes.
For ease of representation, we make the choice of representing the convergence versus the total number of solution samples $N_{\bxi}$.
The results for large $\mathcal{C}_{x}^{max}$ from Figure \ref{subfig1:inlet_errorConv} show that the error control operated in the stochastic space converges quickly at an impressive rate.
The convergence plot of the deterministic contribution counterpart of the error is more difficult to interpret as it exhibits some iterated cycles for which the error increases. This is due to the fact that the new solution samples, which are introduced in order to improve the stochastic surrogate, are first run with the default coarse finite-element mesh discretization. Despite a global second-order convergence rate vs. $N_{\bx}$ (not presented here, c.f. \cite{langenhove2017these}), we also notice that deterministic refinement eventually leads to a stagnation of the error once the maximum complexity is reached $\Ccal^{max}_{\bx}$. Therefore, it indicates that for this application and with our chosen algorithmic setup, the deterministic error dominates over the stochastic error for which the convergence rate is larger (see Table \ref{tab:inlet_stats}). This assessment is still valid when the experiment is repeated for a choice of more modest $\Ccal^{max}_{\bx}=2000$. Moreover, as expected, the average $\varepsilon$ error is larger compared to the previous results. Interestingly, the average stochastic error contribution is larger as well. This finding seems to illustrate a coupling between the errors and requires additional investigations.\\
Next, we wish to investigate if an adaptation solely carried in the stochastic space, suffices in getting satisfactory convergence. To this end, we voluntarily degenerate the spatial discretization of the flow model in the inlet: i.e. we choose to rely on fixed non-adapted meshes for all samples. Results for meshes with $N_{\bx}=\{1500,3000,90000\}$ vertices are presented in Figure \ref{subfig2:inlet_errorConv} and compared to the fully adapted approach. It is very clear that the poor discretization of the retained isotropic meshes  is not capable of accurately capturing the flow features that affect the QoI. They each introduce a certain level of discretization error that eventually impinges on the convergence of the QoI stochastic error. It seems indeed that a too coarse spatial resolution will induce errors and spurious oscillations in the QoI response surface that may in turn mislead the stochastic refinement and induce unnecessary and irrelevant computational burden, leading to a stagnation in the decrease of the stochastic contribution to the approximation error.\\ In conclusion, adaptations in both spaces are needed in this framework. The goal-oriented adaptation in the deterministic space is particularly key, especially for computational budgets with low to moderate mesh cardinalities. Coupling of the adaptive schemes is also crucial and may offer unexpected computational savings due to the interplay between the successive refinements. This is the topic of ongoing investigations.

\begin{figure}[!h]
\centering
\begin{tabular}{cc}
\includegraphics[scale=0.415]{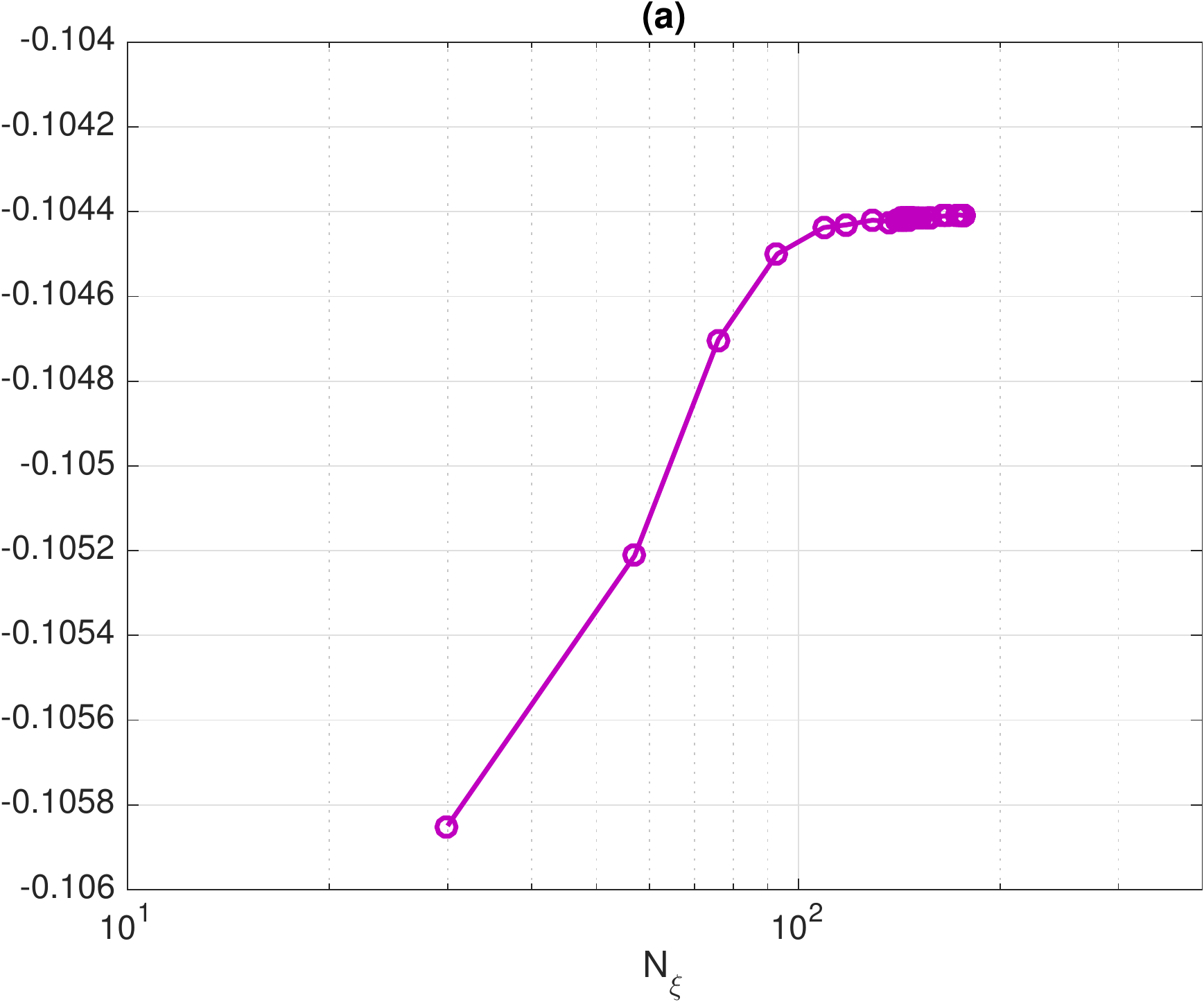} &
\includegraphics[scale=0.415]{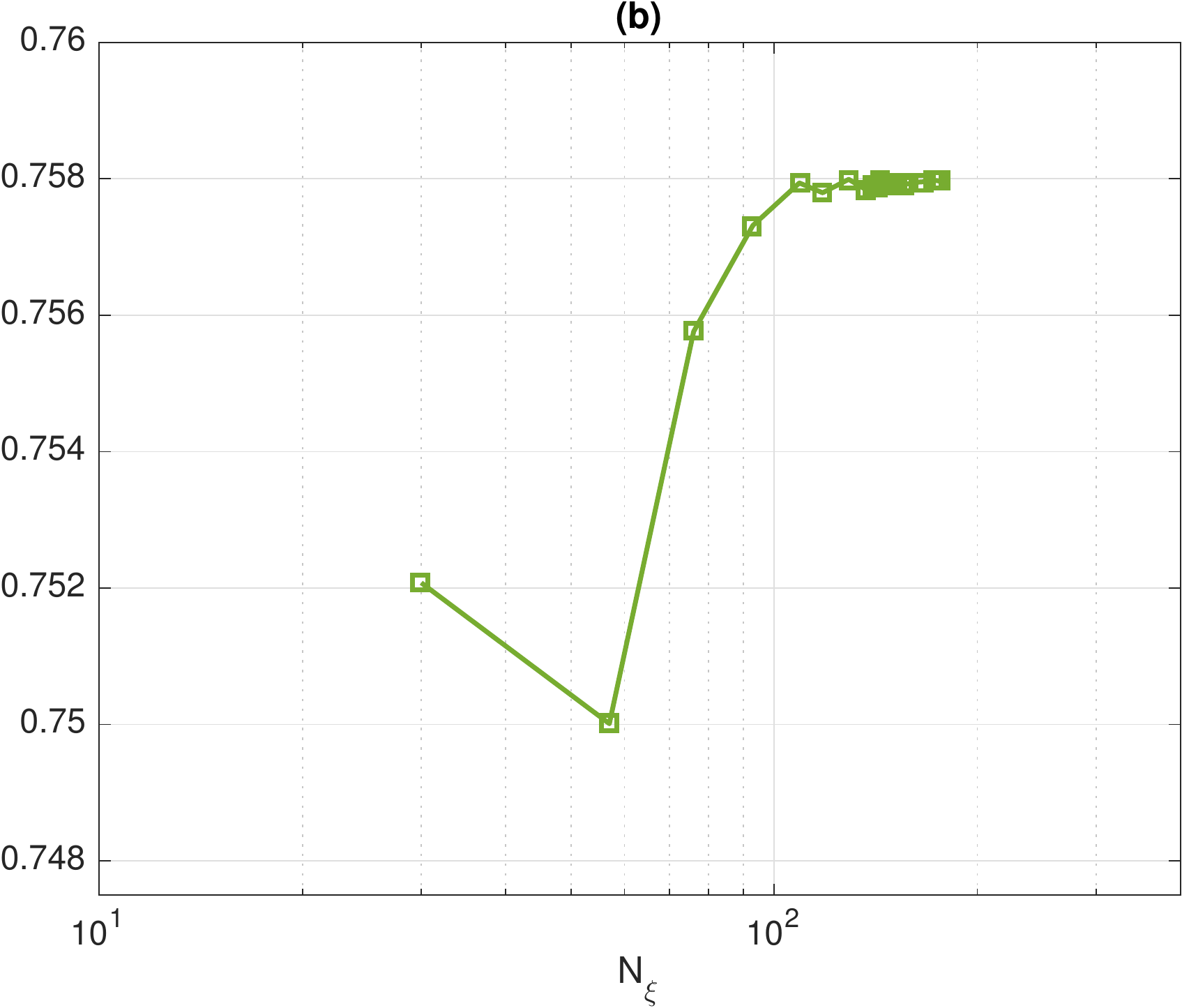} 
 \end{tabular}
\caption{Inlet problem: convergence of the QoI statistics vs. the total number of solution samples $N_{\bxi}$; mean $\overline{j}$ (a) and the coefficient of variation $\text{Cov}_j$ (b).}
\label{fig:inlet_statisticsConv}
\end{figure}

~\\
One may also monitor the convergence of the QoI. Although the exact statistics of $C_p$ for the parametric ranges and distributions considered are not available for our problem we can still visualize the evolution of low-order statistics along the subsequent adaptations of the approximated QoI. The results displayed in Figure \ref{fig:inlet_statisticsConv} hint at some type of numerical convergence for the mean and the coefficient of variation. The large values of the coefficient of variation confirm the large variability of the integrated $C_p$ depending on the changing flow features. \\
As a final note in this subsection, we emphasize the fact that the proposed adaptation algorithms may be improved, potentially speeding up the convergence. For instance, one could avoid using the default mesh complexity target for the new samples introduced. Instead, one may interpolate complexity from the response of the new sample neighbors.

\begin{table} 
    \footnotesize
    \arraycolsep=0.5pt
    \medmuskip=0.3mu
\centering
\begin{tabular}{|l | c | c || c | c | c | c |}
\hline 
Adaptation & \multicolumn{2}{c||}{Discretizations}  & \multicolumn{2}{c|}{QoI statistics} &  \multicolumn{2}{c|}{Error estimations} \\ 
\hline 
Cycle nr. & $N_{\bxi}$ & $\overline{N_{\bx}}$  & $\mathbb{E}[Cp]$ & $\text{Var}[Cp]$ & $\BigE[{\varepsilon}]$ & $\BigE[{\eta}]$\\ 
 \hline 
 \hline 
 $0$ & $30$ & $4000$ & $-0.10584943$ & $0.00633727$ & $0.0253$ & \\
 \hline
 $1$ & $57$ & $4000$ & $-0.10521272$ & $0.00622762$ & $0.0286$ & $0.0236$\\
 \hline
 $2$ & $57$ & $6446$ & $-0.10521254$ & $0.00622689$ & $0.0175$ & $0.0236$\\
 \hline
 $3$ & $118$ & $4885$ & $-0.10443541$ & $0.00626330$ & $0.0229$ & $0.0014$\\
 \hline
 $4$ & $118$ & $82373$ & $-0.10443116$ & $0.00626273$ & $0.0024$ & $0.0014$\\
 \hline
 $5$ & $142$ & $70471$ & $-0.10441939$ & $0.00626304$ & $0.0057$ & $0.0008$\\
 \hline
 $6$ & $142$ & $88849$ & $-0.10441761$ & $0.00626286$ & $0.0023$ & $0.0008$\\
 \hline
 $7$ & $144$ & $88708$ & $-0.10441652$ & $0.00626417$ & $0.0023$ & $0.0007$\\
 \hline
 $8$ & $144$ & $89294$ & $-0.10441650$ & $0.00626417$ & $0.0023$ & $0.0007$\\
 \hline
 $\vdots$ & $\vdots$ & $\vdots$ & $\vdots$ & $\vdots$ & $\vdots$ & $\vdots$\\
 \hline
 $18$ & $177$ & $89351$ & $-0.10441013$ & $0.00626303$ & $0.0024$ & $0.0004$\\
 \hline
 $19$ & $177$ & $90096$ & $-0.10440860$ & $0.00626287$ & $0.0022$ & $0.0004$\\
 \hline
\end{tabular}
\caption{Inlet problem: adaptation statistics and estimated average errors of the coupled error control strategy with $\mathcal{C}_{\bx}^{max}=128 000$. Not all adaptation cycles are displayed.}
\label{tab:inlet_stats}
\end{table}

\section{Conclusion}
This work  proposes a new extension of a deterministic goal-oriented anisotropic mesh adaptation, used to control the discretization/approximation error on a scalar-output functional, to the case of system models bearing a parametric stochastic component. The proposed method is particularly relevant to nonlinear systems such as complex compressible flows with few uncertain parametric dimensions but leading to output quantities of interest (QoI) presenting low regularity and isotropy (or even discontinuities) with respect to the parameters. These situations are for instance encountered when sharp solution gradients become very sensitive to small changes in probable parameters values. Ultimately, the goal is to control the two sources of numerical errors acting on the prediction of the QoI statistics due to deterministic and stochastic discretizations/approximations.\\
The original and the new derivations both rely on a continuous framework to model a mesh/partition of their respective space and its elements, which exploits the duality between Riemannian metric spaces and discrete meshes. 
Besides, the formulations require the minimization of {error estimates of the numerical approximation errors,} that is then recast in the continuous mesh framework. {From the error estimate an optimal metric field is deduced that governs the anisotropic mesh adaptation}. This is realized by solving an optimization problem where one seeks an optimal mesh minimizing the error under some constraint on the cardinality of the targeted discretizations.
On the deterministic side, a ``goal-oriented" approach is followed minimizing {an {\em a priori} error estimate of the {\em system solution}}, thanks to the use of an adjoint solution of the dual problem, acting as a weight that optimally distribute the degrees of freedom
for the specific target.   
The analysis is deployed for a particular numerical scheme, here used to solve 2D compressible Euler system of equations.\\
On the stochastic side, the error estimation {directly relates to the {\em QoI} through interpolation error and} is similar to the one used for deterministic ``feature-based" mesh adaptation {found in the literature}. 
Indeed, the proposed error estimate is independent of the considered {model} and relies only on "geometric" information instead. This is because of the absence of any differential operators along the parametric dimensions.
In this case, there is no need for adjoint solution and only the QoI hessian is used. However the error is weighted by the parameters probability measures. This makes very natural the integration of complex and/or correlated measures in the adaptation. \\

In order to take full advantage of the existing framework and numerical tools developed to manage the main components involved in the adaptation, i.e. the flow solver, the error estimate, the solution interpolation, and the mesh generator, the choice is made to rely on a robust and flexible stochastic collocation with {\em linear} simplex elements, allowing unstructured $h-$type adaptivity and guaranteeing  a Local Extremum Diminishing property. This makes easier the adaptation of the numerical tools to the stochastic framework.\\
The original formulation was shown to be very efficient for goal-oriented adaptation when working with error estimations controlled in $L^p$-norms. In particular, it was shown that anisotropic mesh adaptations achieve  second-order convergence in these norms even for nonlinear systems with discontinuous responses.
For this work, the estimation of both errors due to deterministic and stochastic discretizations are derived in $L^1$-norm. Inherently, this guarantees a convergence in the mean of the stochastic surrogate and therefore ensures that the mean value of the approximated QoI does convergence to the exact QoI mean. Ongoing work is pursued in order to extend the stochastic error estimation and control to $L^2$-norm. \\
For the numerical examples treated in this paper, we have shown a {convergence in the mean} for multivariate nonlinear {discontinuous} functionals of the stochastic space and we have at least recovered a rate of convergence of $\Ocal \left (N_{\bxi}^{-2 /d_{\bxi}} \right )$. The algebraic convergence of this approach is obviously  strongly impacted by the dimensionality of the stochastic problem. For moderate number of random dimensions ($d_{\bxi}>4$) the approximation becomes less efficient than Monte-Carlo like sampling methods. Moreover, it was reported in the literature that Delaunay type triangulation becomes poorly conditioned for a number of random dimensions superior to eight \cite{edeling2016simplex}. Nevertheless, in our opinion, this first attempt establishes an interesting proof-of-concept.
Moreover, the method is robust: for a given functional with a fixed number of parameters, a change in the nature and/or the regularity of the probability measure of the parameters does not affect the convergence rate, nor the choice of the refinement step size, and flexible: capable of handling non-hypercube probability spaces.
~\\

This work also proposes a way of reducing the {\em total} approximation error in the QoI by combining the aforementioned metric-based anisotropic adaptations. As it is critical to identify which approximation space contributes more to this error, algorithmic developments are presented to tackle the quantification and the adaptive adjustment of the error components in the deterministic and stochastic approximation spaces. 
The capability of the proposed approach is tested on various problems including a supersonic scramjet inlet subject to geometrical and operational parametric uncertainties. It is demonstrated to generate highly anisotropic unstructured meshes that accurately capture compressible discontinuous flow features impacting pressure-related quantities of interest, while balancing computations and refinements in both spaces. From a practical point of view, prescribing discretizations minimal size is not required but given a computational budget, maximal number of simulations and affordable highest mesh cardinality are to be defined.\\
For long-term perspectives, the error control strategy may benefit from: -- the derivation of error estimators for higher-order approximation schemes \cite{alex2018} and curved meshes, -- the extension of metric-based adaptation to unsteady problems potentially including moving geometries \cite{barral} and -- the development of mesh-adaptation for turbulent Navier-Stokes equations. 
More specifically, several potential improvements have been identified throughout this work. First and easiest, the stochastic error estimation may be derived and adapted in the $L^2$-norm and this should therefore guarantee a convergence in the variance of the stochastic surrogate. 
Another aspect is $h/p$ adaptation that should bring computational savings in this framework. A possible path would be the one of Yano and Darmofal \cite{YanoJCP2012} who proposed a Riemannian-based optimization to control an error estimate for (Discontinuous Galerkin) numerical schemes of any order $p$ across unstructured meshes. One may consider an extension of their work in which the mesh is parameterized by not just a metric field but also a polynomial order field. This would obviously require some strategy to measure the error sensitivity to the solution order. Another more straightforward but less optimal choice would be to first perform an $h$ adaptation of the partition as proposed in this work and then to adapt {\em a posteriori} the order of the stochastic approximation in each simplex depending on the approximated solution regularity and the availability of higher-order interpolation stencils of vertices of surrounding elements, as proposed in \cite{witteveen2010simplex}.

\section*{Acknowledgments}
The authors are thankful to our late colleague, Dr. Ir. Jeroen A.S. Witteveen, for valuable early discussions related to some of the aspects of this project.

\bibliographystyle{plain}  
\bibliography{references_these}

\appendix

\section{Simplex-Stochastic Collocation Approximation}\label{ssc_elem}

In the following, we provide a brief outline of a baseline version of the Simplex (elements) Stochastic Collation (SSC) approximation, that was greatly developed and extended by Witteveen {\em et al.} in the context of robust adaptive uncertainty quantification, see for instance \cite{witteveen2009adaptive,witteveen2012refinement,witteveen2013subcell, witteveen2013simplex}.

The baseline version of the method is straightforward. It relies on a piecewise multivariate polynomial approximation of a QoI $j$ that depends on a vector of random parameters $\bxi \in \Xi \subset \BigR^{d_{\bxi}}$, i.e. $j(\bxi)$. The SSC method 
discretizes the parameter space into a tessellation of $N_{\Xi}^{elem}$ simplices.
On a linear Lagrange element $\Xi_{(i)}$ the continuous representation of the response surface is obtained by an interpolation on $N_q$ points; we define the interpolation operator on element $\Xi_{(i)}$:
\begin{align*}
    \mathcal{I}_{\Xi_{(i)}}j(\bxi) = \sum_{k=1}^{N_q} j(\bxi_{(k)}) b_k(\bxi),
\end{align*}
where $b_k$ are the barycentric coordinates of $\bxi$ locally in $\Xi_{(i)}$.
The global interpolant on mesh $\mathcal{H}_{\xi}$ can then be written as
\begin{align*}
    \mathcal{I}_{\mathcal{H}_xi} j(\bxi)|_{\Xi_{(i)}} = \mathcal{I}_{\Xi_{(i)}} j(\bxi).
\end{align*}

The tessellation of the parametric domain decomposes any integral over the parameter space into a summation of integrals over $N_{\Xi}^{elem}$ simplices. Let $\Omega_{\xi,(i)}$ be the $i^{th}$ such element and $\Xi_{(i)}$ its image in the parameter space. The integral over each element is approximated by a quadrature, e.g. Newton-Cotes (NC) quadrature, using $N_q$ quadrature points, while the quadrature weights must account for the local probability measure. For simplicity, we assume that $N_q$ is the same for each element.\\
Hence the expectation will be approximated as:
\begin{align}
    \begin{split}
        \mathbb{E}[j(\bxi)] &= \int_{\Xi} j(\bxi) \rho_{\bxi} \mathrm{d} \bxi\\
                            &= \sum_{i=1}^{N_{\Xi}^{elem}} \int_{{\Xi}_{(i)}} j(\bxi) \rho_{\bxi} \mathrm{d}\bxi \approx \sum_{i=1}^{N_{\Xi}^{elem}} \sum_{k=1}^{N_q} c_{(i,k)} j(\bxi_{(i,k)}),
    \end{split}
    \label{appox_exp}
\end{align}
where $\bxi_{(i,k)}$ and $c_{(i,k)}$ are the parameter value and the quadrature weight corresponding to the $k^{\text{th}}$ quadrature point in the $i^{\text{th}}$ element. For a given pdf and tesselation, quadrature weights are computed once and for all as integrals of Lagrange basis:
\begin{align}
    c_{(i,k)} = \int_{\Xi_{(i)}} L_{(i,k)}(\xi_{1}, \xi_{2}, \ldots, \xi_{d_{\bxi}}) \rho_{\bxi}(\xi_{1}, \xi_{2}, \ldots, \xi_{d_{\bxi}}) \mathrm{d}\bxi \  \  \text{ with } i=1, \ldots, N_{\Xi}^{elem},
    \label{n-c_wgt_def}
\end{align}
where $L_{(i,k)}(\bxi)$ is the Lagrange polynomial in the $i^{\text{th}}$ element corresponding to the $k^{\text{th}}$ quadrature point evaluated at $\bxi$. In this work, first degree NC quadrature rule is used. This corresponds to one quadrature/sample point at each vertex of the $d_{{\xi}}-$simplex. In that case $L_{(i,k)}$ are linear Lagrange polynomial on element $\Xi_{(i)}$, taking value $1$ on vertex $k$ and $0$ on all other vertices. \\
As proposed in \cite{witteveen2009adaptive}, a higher degree quadrature with standard weights is used in order to compute the integral in (\ref{n-c_wgt_def}). Let $q_{(i,l)}$ be the standard NC quadrature weight from a quadrature rule containing $N_{q_{sub}}$ points, then the weights $c_{(i,k)}$ can be computed as
\begin{align}
    c_{(i,k)} \approx \sum_{l=1}^{N_{q_{sub}}} q_{(i,l)} L_{(i,k)}(\bxi_{(i,l)}) \rho_{\bxi}(\bxi_{(i,l)}).
    \label{c_def}
\end{align}

NC quadrature rules up to degree $8$ have been implemented in $d_{\bxi}=2$ and up to degree $6$ in $d_{\bxi}=3$. 
Moreover, these computations may be facilitated and made more efficient by mapping whichever simplex to a reference element. 
Do note that these higher degree quadratures are not needed to integrate the linear Lagrange polynomials but merely account for the complexity introduced by the local probability density functions.
Unless mentioned otherwise, we used a $5^{\text{th}}$ degree NC quadrature for $2d$ problems with smooth pdf, and a $3^{\text{rd}}$ degree NC quadrature in $3d$ cases with smooth pdf.\\ 

The non-intrusive SSC method is used in our framework as a robust and efficient stochastic collocation approximation. It is based on an unstructured tessellation rather than relying on the more common tensor product structure and allows for anisotropic refinements. Moreover, it is capable of handling non-hypercube probability spaces \cite{witteveen2012simplex}.
Note that our approximation differs slightly form the SSC method proposed in \cite{witteveen2009adaptive}.
In \cite{witteveen2009adaptive} a second degree NC is used resulting in a piecewise quadratic approximation of the stochastic response. Furthermore, in order to avoid unphysical oscillations near singularities, the original SSC method splits quadratic elements into smaller first degree elements when an extremum is detected that is not located at one of the quadrature points. In later publications other refinement criteria were investigated \cite{witteveen2012refinement} as well as more elaborate schemes in order to obtain truly discontinuous representations \cite{witteveen2013subcell, witteveen2013simplex}.
In this work, only linear elements are used and we rely only on $h-$adaptivity to control the error. The $h-$adaptivity used in this work is also very different from the one applied in the SSC method which splits existing elements according to a predefined patterns whereas the $h-$adaptivity in this paper uses a metric-based approach allowing for far greater flexibility.

\end{document}